\documentclass[18pt,reqno]{amsart}
\usepackage{amsmath,amssymb,amsthm,hyperref}
\usepackage{stmaryrd}
\usepackage{enumerate}
\usepackage{graphicx}
\usepackage{float}
\usepackage{tikz}
\usepackage{url}
\usepackage{mathrsfs}
\usepackage{mathtools}
\usepackage{subcaption} 
\usepackage{geometry}
\usepackage[ backend=biber, 
style=numeric, 
sorting=nyt, 
giveninits=true, 
doi=false, 
url=false, 
isbn=false ]{biblatex} 
\renewbibmacro{in:}{} \DeclareFieldFormat[article]{pages}{#1} 

 \newtheorem{theorem}{Theorem}[section]
 
 \newtheorem{corollary}[theorem]{Corollary}
 
 \newtheorem{lemma}[theorem]{Lemma}
 \newtheorem{proposition}[theorem]{Proposition}
 \newtheorem{example}[theorem]{Example}
 
 \theoremstyle{definition}
 \newtheorem{definition}[theorem]{Definition}
 \newtheorem{remark}[theorem]{Remark}

\newcommand{\N}{\mathbb{N}}
\newcommand{\Z}{\mathbb{Z}}
\newcommand{\R}{\mathbb{R}}
\newcommand{\D}{\mathbb{D}}
\newcommand{\Q}{\operatorname{Q}}

\newcommand{\FF}{\mathcal{F}}

\newcommand{\T}{\mathbb{T}}

\newcommand{\E}{\mathbb{E}}

\newcommand{\dint}{{\rm d}}
\newcommand{\dTV}{\operatorname{d}_{\rm TV}}
\newcommand{\dKOL}{\operatorname{d}_{\rm Kol}}
\newcommand{\dWass}{\operatorname{d}_{\rm W}}
\newcommand{\HH}{\mathcal{H}}

\newcommand{\I}{\mathbb{I}}
\renewcommand{\P}{\mathbb{P}}
\newcommand{\dHH}{\operatorname{d}_{\HH}}

\newcommand{\WW}{\mathcal{W}}

\renewcommand{\dint}{{\rm d}}

\renewcommand{\L}{\operatorname{L}}
\newcommand{\Id}{\operatorname{Id}}
\newcommand{\RR}{\operatorname{R}}
\newcommand{\Pt}{\operatorname{P}}

\newcommand{\Ga}{\Gamma}
\newcommand{\Var}{\operatorname{Var}}

\newcommand{\Tj}{\operatorname{T}}
\renewcommand{\S}{\operatorname{S}}

\newcommand{\LL}{\mathcal{L}}

\newcommand{\Cpol}[1]{\mathscr{C}_{\text{pol}}^{#1}(\mathbb{R})}
\newcommand{\Cpold}[2]{\mathscr{C}_{\text{pol},\sigma}^{#1}(\mathbb{R}^{#2})}

\newcommand{\projgeq}[1]{\operatorname{Proj}_{\geq #1}}

\newcommand{\J}{\operatorname{J}}

\newcommand{\rank}{\operatorname{rank}}
\newcommand{\Fvect}{\underline{F}}
\DeclareMathOperator{\prob}{\mathbf{P}}

\DeclareMathOperator{\Dom}{\operatorname{Dom}}

 \DeclarePairedDelimiter{\norm}{\lVert}{\rVert}
 \numberwithin{equation}{section}
\title{Edgeworth expansion on sums of Wiener chaos}
\author{Paul Mansanarez}
\address{
	Paul Mansanarez, Nantes Universit\'e/Universit\'e libre de Bruxelles, France/Belgium. E-mail:
	paul.mansanarez@ulb.be}
\begin{document}

\maketitle

\begin{abstract}
    In this paper, we study Edgeworth expansions for random varaibles belonging to a fixed finite sum of Wiener chaoses. Building on the work of Mansanarez, Poly and Swan \cite{MPS2025.EE} for a single Wiener chaos, we investigate the extension of their approach to functionals $F$ living in a finite sum of Wiener chaoses $\WW_{\leq d}$. Within the Malliavin–Stein framework, we bound the total variation distance between the law of $F$ and its Edgeworth approximation of order $2dm$ by $(\sum_{p,q}^d \Var\Gamma(F_p,F_q))^{m/2}(\Var\Gamma(F,\L^{-1}F))^{1/2}$, where $F_p$ denotes the projection of $F$ onto the $p$-th Wiener chaos. The presence of several chaotic components introduces additional interactions, since these components are orthogonal but generally dependent. Since our proof rely on the integration-by-parts formalism and the spectral properties of the Ornstein-Uhlenbeck operator, this result extends naturally to the broader framework of Markov chaoses.

\end{abstract}

\section{Introduction}

Edgeworth expansions provide a systematic refinement of Gaussian approximation by identifying correction terms that describe deviations from the limiting normal distribution. For nonlinear functionals of Gaussian fields, obtaining such expansions requires quantitative estimates that account for the underlying dependence structure. Wiener chaos offers a natural framework for this problem, through the interplay between Gaussian analysis, Malliavin calculus, and the spectral properties of the Ornstein–Uhlenbeck generator.

\medskip

The starting point of the present work is the article \cite{MPS2025.EE}. For a random variable $F$ with variance one, in the $p$-th Wiener chaos $\WW_p$, $p\geq 2$, they consider the signed measure
\begin{align*}
    \boldsymbol{\gamma}_{F,M}(\dint x) = \frac{e^{-\frac{x^2}{2}}}{\sqrt{2\pi}}\left( 1+\sum_{k=3}^{M} \frac{\mathbb E[H_k(F)] }{k!}H_k(x) \right)\,\dint x, \qquad M\geq3,
\end{align*}
where $H_k$ denotes the $k$-th Hermite polynomial. They establish
\begin{align*}
    \dTV\left(\prob_F,\boldsymbol{\gamma}_{F,M}\right) \leq C_{p,m} \bigl ( \Var\Gamma(F,F) \bigl )^{\frac{m+1}{2}} \, ,
\end{align*}
where $\prob_F$ denotes the law of $F$, $\Gamma(F,F)$ the \emph{carré-du-champ} operator of $F$ and $M = 4m-1$, $m\geq 1$. Their approach combines Malliavin–Stein arguments with spectral analysis on the \emph{Ornstein-Uhlenbeck} operator in dimension one.

\medskip

This result naturally raises the question of extending the approach to random variables belonging to a finite sum of Wiener chaoses. We therefore consider centered random variables of the form
\begin{align*}
    F=\sum_{p=1}^{d} F_p \, ,\qquad \mathbb{E}[F^2]=1 \, ,
\end{align*}
where $d\geq 2$ is fixed, $F_p \in \WW_p$ for all $p$. This setting accommodates polynomial Gaussian functionals with several nonzero chaotic projections, allowing different orders to contribute simultaneously to the fluctuations.
\medskip

The passage to a finite sum changes the spectral structure of the problem and introduces a genuine difficulty. For a random variable in a single chaos, the identity $LF=-pF$ directly relates the Ornstein–Uhlenbeck generator $\L$ to the variable itself, an identity that considerably simplifies the integration-by-parts formulas underlying Malliavin--Stein approximation (see, e.g., \cite{NP2012}). For a finite sum of chaoses, this scalar eigenfunction relation is no longer available, one has instead
\begin{align*}
    \L F=-\sum_{p=1}^{d} pF_p\, ,
\end{align*}
so the operator $\L$ acts differently on each component. Furthermore, although distinct chaotic components are orthogonal, they are generally dependent, and interactions between the components enter both the correction coefficients and the remainder. These observations suggest that extending the argument requires a careful treatment of the full chaotic decomposition, and that establishing an expansion in this setting requires more than applying the fixed-chaos result separately to each component. 

\medskip

The purpose of this paper is to investigate this extension, building directly on the method of \cite{MPS2025.EE}. Our objective is to understand how their expansion procedure can be adapted from a single Wiener chaos to a finite sum of chaoses. The main challenge is to identify quantities that control the resulting approximation error. 

\subsection{Asymptotic normality for functionals of a Gaussian field}

The Gaussian approximation of nonlinear functionals of Gaussian fields has been extensively studied through the combination of Malliavin calculus and Stein’s method. An important starting point is the celebrated \emph{Fourth Moment Theorem} of Nualart and Peccati \cite{NP2005}, which states that a sequence of random variables $F_n$ in a fixed Wiener chaos $\WW_p$, with unit variance, converges in distribution to the standard Gaussian distribution if and only if its fourth cumulant $\kappa_4(F_n)$ converges to $0$. Soon after, in \cite{NO2008}, Nualart and Ortiz-Latorre show that the latter condition is equivalent to
\begin{align*}
    \Var\Gamma(F,F) \xrightarrow[n\to +\infty]{} 0\, ,
\end{align*}
where $\Gamma(F,F)$ is the \emph{carré-du-champ} operator of $F$, which can also be expressed as $\langle DF,DF\rangle_{\mathfrak{H}}$, where $D$ denotes the Malliavin derivative operator (see \cite{NP2012} for more details), therefore drawing a interesting connection with Malliavin calculus. Then, in \cite{NP2009b} Nourdin and Peccati developed a quantitative approach known as the \emph{Malliavin--Stein} method, or \emph{Nourdin--Peccati} analysis, by combining Malliavin calculus with Stein’s characterization of the normal distribution. For a centered random variable $F\in\mathbb{D}^{1,2}$ with unit variance and $N\sim\mathcal N(0,1)$, a basic estimate is
\begin{align}\label{intro:sm_bound}
    \dWass(F,N) \leq \bigl ( \Var\Gamma(F,\L^{-1} F) \bigr )^{\frac{1}{2}} \, ,
\end{align}
where $\L^{-1}$ is the \emph{pseudo-inverse} of $\L$ and $\dWass$ denotes the Wasserstein-$1$ distance. For $F$ in the $p$-th Wiener chaos, the latter bound \eqref{intro:sm_bound} reduces to $p^{-1} (\Var\Gamma(F,F))^{1/2}$. This identity leads to explicit approximation bounds in terms of the square root of the fourth cumulant, equivalent quantity. Quantitative refinements of \eqref{intro:sm_bound} for finite sums of Wiener chaoses have been investigated in \cite{BCKBS2025,DDMN2026,EsSebaiy2026,FT2016}. These results, however, either impose additional structural restrictions or yield bounds of a different form. 

These estimates also motivate a more precise description of the approximation error. Edgeworth expansions address this problem by adding explicit correction terms to the Gaussian law, with coefficients determined by moments or cumulants. In their work \cite{NP2009a} on exact Berry–Esseen asymptotics, Nourdin and Peccati showed that Malliavin–Stein arguments can identify the leading term of the error under suitable joint convergence assumptions. Namely, they prove under appropriate conditions, the following local one-term Edgeworth expansion
\begin{align*}
    \P(F_n \leq z) = \Phi(z) - \frac{\E[H_3(F_n)]}{6} \Phi^{(3)}(z) + o_z(\varphi(n)) \, , \qquad z \in \R\, ,
\end{align*}
where $(F_n)$ is a sequence of elements of $\D^{1,2}$ and $\varphi(n) := \bigl( \Var\Gamma(F_n,\L^{-1} F_n) \bigr )^{\frac{1}{2}}$ (see also \cite{Campese2013,KP2018}). This results show that $\bigl( \Var\Gamma(F_n,\L^{-1} F_n) \bigr )^{\frac{1}{2}}$ should be the correct bound to use to quantify an Edgeworth expansion for functionals of a Gaussian field, in particular elements of a finite sum of Wiener chaoses. For more details on Edgeworth expansions and the Nourdin--Peccati analysis, we refer to \cite{MPS2025.EE,NP2012}. 

\subsection{Main result}
Our main result provides higher-order Gaussian approximations for random variables belonging to a finite sum of Wiener chaoses. In contrast to the homogeneous setting, several chaotic components may contribute simultaneously, and the approximation must account for their interactions. The result is formulated directly in total variation and requires no decomposition into independent summands.
Let $F$ be a centred element of $\WW_{\leq d}$ such that $\E[F]=0$ and $\E[F^2]=1$. For $p \in \{1,\ldots,d\}$, denote by $F_p$ the orthogonal projection of $F$ onto $\WW_p$. For an integer $M>0$, define the signed measure $\boldsymbol{\gamma}_{F,M}$ with Radon--Nikodym derivative
\begin{align*}
x \longmapsto \frac{e^{-x^2/2}}{\sqrt{2\pi}} \left ( 1 + \underset{k=3}{\overset{M}{\sum}} \, \frac{\E[H_k(F)]}{k!}H_k(x) \right )\,,
\end{align*}
where $H_k$ denotes the Hermite polynomial of order $k$ (see Section \ref{sec.polynomials} for a precise definition). The coefficients of this approximation are explicit functions of the moments, or equivalently of the cumulants, of $F$. For instance, writing $\kappa_j(F)$ for the $j$-th cumulant, the normalization of $F$ gives 
\begin{align*}
    \E[H_3(F)]=\kappa_3(F) \, , \qquad \E[H_4(F)]=\kappa_4(F) \, , \qquad \E[H_5(F)]=\kappa_5(F)\, ,
\end{align*}  
and $\E[H_6(F)] = \kappa_6(F)+10\,\kappa_3(F)^2$. Thus the correction incorporates both individual cumulants and their interactions. The measure $\boldsymbol{\gamma}_{F,M}$ has total mass one, but need not be nonnegative. Moreover, Hermite orthogonality implies the exact moment-matching property 
\[ \int_{\mathbb R}P(x)\, \boldsymbol{\gamma}_{F,M}(dx) = \E[P(F)]\,, \] 
for every polynomial $P$ of degree at most $M$, since $H_0,\ldots,H_M$ form a basis of the space of polynomials of degree at most $M$. The main question is whether this explicit moment-matching correction also approximates the law of $F$ with a quantitatively improved error. To state our first result, we introduce
\begin{align*}
    \Lambda(F) := \sum_{1\leq p\leq q\leq d} \Var\Gamma(F_p,F_q) \, .
\end{align*}
The quantity $\Lambda(F)$ records fluctuations of the diagonal and mixed carr\'e-du-champ terms associated with the chaotic decomposition. In this article, we prove the following extension of the result of \cite{MPS2025.EE}.
\begin{theorem}\label{thm:EE_Wiener_sums}
Let $m$ be a positive integer. Then there exists $C_{d,m}>0$ such that for all $F\in \WW_{\leq d}$ with $\E[F]=0$ and $\E[F^2]=1$,
    \begin{align}\label{eq:maintheorem}
    \dTV(\prob_F,\boldsymbol{\gamma}_{F,M}) & \leq C_{d,m} \, \Lambda(F)^{\frac{m}{2}}\bigl ( \Var\Gamma(F,\L^{-1} F) \bigr )^{\frac{1}{2}}\, ,
\end{align}
where $M = 2dm$.
\end{theorem}

Theorem \ref{thm:EE_Wiener_sums} extends the result of \cite{MPS2025.EE} from random variables belonging to a single Wiener chaos to random variables belonging to a finite sum of Wiener chaoses. As mentioned before, although distinct chaotic projections are orthogonal, they are generally dependent. These interactions are reflected explicitly in the terms $\Gamma(F_p,F_q)$ appearing in $\Lambda(F)$. The estimate also separates two aspects of the approximation. The factor $(\Var\Gamma(F,\L^{-1}F))^{1/2}$ is the discrepancy entering the usual Malliavin--Stein bound for Gaussian approximation, whereas $\Lambda(F)^{m/2}$ quantifies the additional gain obtained by including higher-order corrections. In particular, whenever $\Lambda(F_n)\to0$, the corrected approximation has an error that is smaller than the scale $(\Var\Gamma(F,\L^{-1}F))^{1/2}$. Importantly, the constants depend only on the maximal chaos order $d$ and the expansion order $m$.

\begin{remark}
The quantity $\Lambda(F)$ measures the extent to which the chaotic component of $F$ fail to exhibit a jointly Gaussian structure. In the recent work of Herry, Malicet and Poly \cite{HMP2024+}, they prove that, for a normalized sequence $(F_n)_{n\geq1}$ in $\WW_{\leq d}$, asymptotic normality of the sum forces the vector $(F_{n,1},\ldots, F_{n,d})$ of chaotic projections to become asymptotically Gaussian. More precisely, writing
\begin{align*}
    \sigma_{n,p}^2=\E[F_{n,p}^2] \,, \qquad D_n=\operatorname{diag} (\sigma_{n,1}^2,\ldots,\sigma_{n,d}^2)
\end{align*}
the law of $(F_{n,1},\ldots, F_{n,d})$ approaches the centred Gaussian law with covariance matrix $D_n$. Together with the fixed-chaos criteria for Gaussian convergence and the corresponding contraction estimates, it yields 
\begin{align*}
    F_n\xrightarrow{\mathrm{law}}\mathcal N(0,1) \quad\Longrightarrow\quad \Lambda(F_n)\xrightarrow[n\to+\infty]{} 0 \, .
\end{align*}
Thus, the interaction term $\Lambda(F)$ in the bound in \eqref{eq:maintheorem} vanishes in the central regime.
\end{remark}
The next proposition gives a quantitative comparison of $\Lambda(F)$ in terms of the discrepancy $\Var\Gamma(F,\L^{-1}F)$, in the setting of the Wiener space.
\begin{proposition}\label{intro:bound-by-var}
    Let $d\geq 2$ be an integer. There exists a constant $C_d>0$ such that, for every centered $F\in \WW_{\leq d}$ with $\E[F^2]=1$,
    \begin{align*}
    \Lambda(F) \leq C_d \,\bigl ( \Var\Gamma(F,\L^{-1} F) \bigr )^{\frac{1}{2d-2}}\,.
\end{align*}
\end{proposition}
For an element $F$ of a single Wiener chaos, the identity $-\L^{-1}F=\tfrac{1}{p}F$ reduces the standard Malliavin--Stein discrepancy $\Var\Gamma(F,\L^{-1}F)$ to a multiple of the variance of $\Gamma(F,F)$. For $F$ in a finite sum, however, one has $-\L^{-1}F=\sum_{q=1}^d\tfrac{1}{q}F_q$, and hence 
\begin{align*}
    \Gamma(F,-\L^{-1}F) = \sum_{p,q=1}^d\frac{1}{q} \, \Gamma(F_p,F_q)
\end{align*}
The term $\Var\Gamma(F,\L^{-1}F)$ therefore involves a weighted sum of mixed interactions, in which cancellations may occur. Controlling each componentwise fluctuation through this single quantity is consequently a substantive additional step. In particular, the proposition shows that, within a fixed finite sum of Wiener chaoses, convergence of $\Var\Gamma(F,\L^{-1}F)$ to zero forces the fluctuations of every componentwise term $\Gamma(F_p,F_q)$ to vanish, with an explicit quantitative bound. Combining the comparison from Proposition \ref{intro:bound-by-var} with Theorem \ref{thm:EE_Wiener_sums} yields an expansion whose remainder term is solely controlled in terms of $\Var\Gamma(F,\L^{-1}F)$.
\begin{corollary}\label{cor:EE_Wiener_sums}
Let $m$ be a positive integer. Then there exists $C_{d,m}>0$ such that for all $F\in \WW_{\leq d}$ with $\E[F]=0$ and $\E[F^2]=1$,
    \begin{align*}
    \dTV(\prob_F,\boldsymbol{\gamma}_{F,M}) & \leq C_{d,m} \bigl ( \Var\Gamma(F,\L^{-1} F) \bigr )^{\frac{m+1}{2}}\, ,
\end{align*}
where $M = 4d(d-1)m$.
\end{corollary}

The proof of Theorem \ref{thm:EE_Wiener_sums} uses the spectral and diffusion properties of the Ornstein--Uhlenbeck generator, rather than a representation in a fixed finite number of Gaussian coordinates. This suggests an extension to suitable diffusive symmetric Markov generators, under the corresponding spectral and algebraic assumptions. As in \cite[Section 7]{MPS2025.EE}, this perspective builds on the insightful article \cite{Ledoux2012} and the subsequent developments \cite{ACP2014,AMMP2016,CNPP2016}. The appropriate analogue of Theorem \ref{thm:EE_Wiener_sums} can be formulated in that framework by retaining the componentwise interaction term. We refer to \cite[Section 7]{MPS2025.EE} for the relevant Markov-chaos structure. By contrast, Proposition \ref{intro:bound-by-var} is established here specifically on the Wiener space. Consequently, the reduction to a single scalar discrepancy in Corollary \ref{cor:EE_Wiener_sums} should not be asserted for general Markov chaoses without an additional comparison estimate.

\subsection{Structure of the paper}

The paper is structured as follows. In Section \ref{sec.prelim} , we provide notations and remind useful results that will be used throughout the article. In Section \ref{sec.overview}, we outline the main strategy of the proof of Theorem \ref{thm:EE_Wiener_sums}. In Section \ref{sec-FO}, we establish the theorem in the case $m=1$ and develop the key ingredients needed to address the general case.
In Section \ref{sec-HO}, we prove Theorem \ref{thm:EE_Wiener_sums} in full generality. Finally, in Section \ref{sec.shenanigans}, we prove Proposition \ref{intro:bound-by-var} and compare the different quantities at play on specific settings and examples.

\section{Preliminaries}\label{sec.prelim}
In this section, we provide the necessary definitions, notations, and preliminary lemmas required for the proof of the main theorem. An important part of the result are easily find in the literature, hence we provide proofs for the sake of clarity. Otherwise, we refer the reader to standard references.

\subsection{Notations}

Let us start with some notations.

\begin{itemize}
\item The set of nonnegative integers is denoted by $\N$, and the set of positive integers by $\N^*$. For $a,b$ integers such that $a\leq b$, we will denote by $\llbracket a, b \rrbracket$ the set $\{a, a+1, \ldots,b\}$. For a real number $\alpha$, we will also denote by $\lceil \alpha \rceil$ (resp. $\lfloor \alpha \rfloor$) the ceiling (resp. floor) function evaluated at $\alpha$ that is, the smallest (resp. biggest) integer $k$ such that $k-1 < x \leq k$ (resp. $k\leq x <k+1$).

\item We denote by $\mathscr{C}_b^0(\R)$ the set of continuous functions $h:\R \rightarrow \R$ that are bounded.

\item For every $m$ positive integer, we write $f^{(m)}$ the $m$-th derivative of $f$. For $m=1,2$, we will also use the classical notations $f', f''$ for the first and second derivatives. If $f$ has two variables $(t,x)$, with $t \in \R$ and $x \in \R^d$, we will denote by $\partial_tf$ the first partial derivative with respect to the (time) variable $t$. We will denote also, for $p,q\in \llbracket 1,d \rrbracket$, $\partial_{p} f$, $\partial_p^2 f$, $\partial_{pq}^2 f$ respectively the first and second partial derivative of $f$ with respect to the $p$-th coordinate $x_p$, and the second partial derivative of $f$ with respect to $x_p$ and $x_q$. Finally, if $\alpha$ is an element of $\N^d$, we will denote by $\partial^\alpha_x$ the partial derivative $\partial_{1}^{\alpha_1} \cdots \partial_{d}^{\alpha_d}$, where by convention $\partial_{p}^{0} = \Id$.

\item For a non-negative integer $m$, and $\sigma = (\sigma_1,\ldots, \sigma_d) \in ]0,1]^d$, we denote by $\Cpold{m}{d}$ the set of $m$-times differentiable functions $f$ such that, for all $\alpha \in \N^d$ such that $|\alpha| := \sum_{i=1}^d \alpha_i = m$, $\partial_x^\alpha f$ is continuous, and for all $\beta \in \N^d$ such that $\beta|\leq m$, $\partial_x^\beta f$ has polynomial growth, that is, there exist $C_\beta>0$ and $d_\beta \in \N$ verifying $\forall x \in \R^d, \; |\partial_x^\beta f (x)| \leq C_\beta(1+\|I_\sigma^{-1}x\|^{d_\beta})$, where $I_\sigma$ is the diagonal matrice $\operatorname{diag}(\sigma_1,\ldots, \sigma_d)$. The space $\Cpold{m}{d}$ is naturally endowed with a structure of algebra over $\R$. Finally, we define $\Cpold{\infty}{d} \coloneq \underset{m\in\N}{\bigcap} \Cpold{m}{d}$. When $d=1$ and $\sigma = 1$, we will denote simply $\Cpold{m}{1} = \Cpol{m}$.

\item The identity operator on a vector space $E$ is denoted by $\Id_E$, or $\Id$ when the associated vector space is clear.

\item If $J$ is a linear operator from a vector space $E$ into itself, for $m \in \N$, we will denote by $J^m$ the composition $J\circ \cdots \circ J$ where $J$ appears $m$ times. By convention $J^0$ is the identity operator $\Id_E$. Furthermore, if $P = \prod_{k=1}^m (X-\alpha_k)$ is a polynomial, with $\alpha_1, \cdots, \alpha_m$ real numbers, and for all $k\in \{1,\ldots,m\}$, $(J-\alpha_k\Id_E)$ is invertible, we will denote by $P(J)^{-1}$ the composition $ (J-\alpha_1\Id_E)^{-1} \circ \cdots \circ (J-\alpha_m\Id_E)^{-1}$. Using the classical linear algebra notation, we will write $\alpha \Id_E = \alpha$, for all scalar $\alpha$ and we will also sometimes denote multiplicatively $J(u) = Ju $, for all $u \in E$.

\item Given a probability space $(X,\mathcal{A},\nu)$ and $q\geq1$, we will denote by $L^q(\nu)$ the vector space of $\mathcal{A}$-measurable functions $f:X \rightarrow \R$ such that 
\[\norm{f}_q^q \coloneq \int_X |f(x)|^q \, \nu(\dint x) <+\infty.\]

\item We will denote by $\gamma$ the standard Gaussian distribution on $\R$, that is, the measure with density $x \mapsto e^{-x^2/2}/\sqrt{2\pi}$ with respect to the Lebesgue measure on $\R$. Additionally, we will denote by $N$ a random variable with distribution $\gamma$. Since $\gamma$ admits moments of all orders, we have the inclusion $\Cpol{m} \subseteq \underset{q\geq 1}{\bigcap} L^q(\gamma)$, for all $m \in \N$.

\item For $A$ a measurable set, we denote by $\I_A$ the characteristic function (or indicator function) of the set $A$.

\end{itemize}
\subsection{Analysis on $L^2(\R^d,\gamma_\sigma)$}\label{s:gaussian_analysis_multid}
Let $\sigma_1, \ldots, \sigma_d$ be positive real numbers such that $\sigma^2_1+\cdots +\sigma_d^2 = 1$. Denote by $\sigma$ the vector $(\sigma_1, \ldots, \sigma_d)$. We consider the probability space $(\R^d, \mathscr{B}(\R^d), \gamma_\sigma)$, where $\gamma_\sigma$ is the centered Gaussian measure with covariance matrix $\operatorname{diag}(\sigma_1^2, \ldots, \sigma_d^2)$. It is the tensor product of probability spaces
\begin{align*}
\underset{p=1}{\overset{d}{\bigotimes}} \, (\R, \mathscr{B}(\R), \gamma_i)
\end{align*}
where $\gamma_i$ is the probability measure of $\mathcal{N}(0,\sigma_i^2)$. The probability measure $\gamma_\sigma$ is the invariant and reversible measure of this weighted anisotropic Ornstein-Uhlenbeck operator $\LL_\sigma$, given on twice differentiable function $f:\R^d \rightarrow \R$ by
\begin{align}\label{eq:definition_LL_sigma}
\forall x \in \R^d, \quad \LL_\sigma(f)(x) = \sum_{p=1}^d p \left( \sigma_p^2 \partial_{x_p}^2 f(x) - x_p \partial_{x_p} f(x) \right ).
\end{align}
\begin{remark}
Taking $d=1$ and $\sigma_1=1$, we recover the usual Ornstein–Uhlenbeck operator $\LL$ on the real line.
\end{remark}
In the rest of this section, we will be considering product of diagonal matrices, keeping in mind that in restriction to that space, the product is commutative. 

\subsubsection{Multivariate Hermite polynomials}\label{sec.polynomials}
The operator $\LL_\sigma$ is self-adjoint. An orthogonal basis of $L^2(\R^d,\gamma_\sigma)$ is given by the \emph{generalized multivariate} Hermite polynomials, defined by
\begin{align}\label{eq:definition_hermite_polynomial_d}
H_\alpha^\sigma(x) \coloneq \underset{p=1}{\overset{d}{\prod}} H_{\alpha_p} \left ( \frac{x_p}{\sigma_p} \right )\,, \quad \alpha = (\alpha_1, \ldots, \alpha_d) \in \N^d\,, 
\end{align}
where $(H_k)_{k\geq 0}$ are the univariate Hermite polynomials, that can be defined as
\begin{align*}
    H_k(x) \coloneq (-1)^k e^{x^2/2} \frac{\dint^k}{\dint x^k} e^{-x^2/2} \,, \quad k \in \N \, .
\end{align*}
Above, each $H_\alpha^\sigma$ has eigenvalue 
\begin{align}
-\lambda(\alpha) \coloneq -\sum_{p=1}^{d} p\, \alpha_p\,.
\end{align}
For every $\phi \in L^{2}(\gamma_\sigma)$, one has the expansion
\begin{align*}
\phi = \sum_{\alpha \in \N^d} \, c_\alpha^\sigma(\phi) H_\alpha^\sigma(x)\,, \quad c_\alpha^\sigma(\phi) := \frac{1}{\alpha!} \, \E\left[\phi(N_\sigma )H_\alpha^{\sigma}(N_\sigma)\right ]\,,
\end{align*}
where $N_\sigma$ is a Gaussian vector with mean zero and covariance matrix $\operatorname{diag}(\sigma_1^2,\ldots,\sigma^2_d)$.
The notion of Hermite rank extends as follows.
\begin{definition}\label{def:hermite-rank}
Let $\phi$ be an element of $L^2(\gamma_\sigma)$. The \emph{Hermite rank} of $\phi$ is defined as
\begin{align}
\operatorname{rank}_ \sigma(\phi) \coloneq \min \bigl\{ |\alpha| \;;\; \alpha \in \N^d, \;c_\alpha^\sigma(\phi) \neq 0 \bigr\}\,.
\end{align}
When $d=1$, we simply denote the Hermite rank of $\phi$ by $\rank(\phi)$.
\end{definition}
\noindent The following lemma follows directly from the preceding definition combined with the fact that the family $(H_{\alpha}^\sigma)_{|\alpha|\leq r-1}$ forms a basis of the space of polynomials of degree at most $r-1$.
\begin{lemma}\label{lemma:hermite_rank_d}
Let $\phi \in L^2(\gamma_\sigma)$ with Hermite rank $r > 0$. Then for every polynomial $Q \in \R[X_1, \ldots, X_d]$ of degree at most $r-1$, one has
\[
\E\left [Q(N_\sigma)\phi(N_\sigma)\right ] = 0\,.
\]
\end{lemma}
\noindent In what follows, we interpret a univariate function evaluated at a sum as a multivariate function in the following way.
\begin{definition}\label{definition:phi_oplus}
Let $\varphi$ be an element of $L^2(\gamma)$. We denote by $[\varphi]^\oplus$ the function $\R^d \rightarrow \R$ defined by
\begin{align*}
\forall (x_1,\ldots,x_d) \in \R^d, \quad [\varphi]^\oplus(x_1,\ldots,x_d) = \varphi \left ( x_1+\cdots+x_d \right)\,.
\end{align*}
In particular, when $\varphi$ is continuous and bounded, so is $[\varphi]^\oplus$.
\end{definition}

\begin{remark}
    Let $\varphi$ be an element of $\Cpol{m}$. In the rest of the article, we will use, without reminding it, the fact that $[\varphi]^\oplus \in \Cpold{m}{d}$.
\end{remark}

\noindent By considering the generating function of the Hermite polynomials, one obtains the following multinomial formula for Hermite polynomials.

\begin{lemma}\label{lemma:multinomial_hermite}
Let $k$ be an element of $\N$. Then one has
\begin{align}
H_k \left (x_1+\cdots +x_d \right )= \sum_{\underset{|\alpha|=k}{\alpha \in \N^d}} \, \sigma^\alpha \binom{k}{\alpha} H_\alpha^{\sigma}(x)\,,
\end{align}
where $\binom{k}{\alpha} \coloneq \frac{k!}{\alpha_1 ! \cdots \alpha_d!}$ and $\sigma^\alpha \coloneq \sigma_1^{\alpha_1} \cdots \sigma_d^{\alpha_d}$.
\end{lemma}

\noindent We will be interested in the relation between the Hermite rank of a function $\varphi:\R \rightarrow \R$ and the Hermite rank of its transformation $[\varphi]^\oplus$. Lemma \ref{lemma:multinomial_hermite} yields the following.

\begin{lemma}\label{lemma:egalite-hermite-rank}
Let $\varphi$ be an element of $L^2(\gamma)$. Then $[\varphi]^\oplus \in L^2(\gamma_\sigma)$ and
\begin{align*}
\rank_\sigma\left ([\varphi]^\oplus \right ) = \rank(\varphi)\,.
\end{align*}
\end{lemma}
\begin{proof}
Let $\alpha$ be an element of $\N^d$. Expanding $\varphi$ with respect to the Hermite polynomials, and using Lemma \ref{lemma:multinomial_hermite}, one gets
\begin{align*}
c_\alpha\left ([\varphi]^\oplus \right) & = \frac{1}{\alpha !}\,\E \left [ [\varphi]^\oplus( N_\sigma)H_\alpha^\sigma(N_\sigma) \right ] = \sum_{n = 0}^{+\infty} \frac{c_n(\varphi)}{\alpha !}\,\E \left [ [H_n]^\oplus( N_\sigma)H_\alpha^\sigma(N_\sigma) \right ] \\
& = \sum_{n = 0}^{+\infty} \,\sum_{\underset{|\beta|=n}{\beta \in \N^d}}\,\sigma^\beta \binom{n}{\beta} \frac{c_n(\varphi)}{\alpha !}\,\E \left [ H_\beta^\sigma( N_\sigma)H_\alpha^\sigma(N_\sigma) \right ] = \sigma^\sigma \binom{|\alpha|}{\alpha} c_{|\alpha|}(\varphi)\,,
\end{align*}
where we used the orthogonality in $L^2(\gamma_\sigma)$ of the polynomials $H_\beta^\sigma$ and $H_\alpha^\sigma$ for $\beta \neq \alpha$. That concludes the proof.
\end{proof}

\subsubsection{Hermite rank and projection}

We have the following property of conservation of Hermite rank.
\begin{lemma}\label{lemma:conv_rank_L_d}
    Let $\phi$ be an element of $\Cpold{2}{d}$ such that $\rank_\sigma(\phi)=r$. Then, for every $q > -r$, 
    \begin{align*}
        \rank_\sigma \bigl ((\LL_\sigma-q\Id)\phi\bigr ) \geq r \, .
    \end{align*}
\end{lemma}
\begin{proof}
    Let $\alpha$ be an element of $\N^d$. Then, by symmetry of $\LL_\sigma$, one has
    \begin{align*}
        \alpha !\,c_\alpha^\sigma \bigl ( (\LL_\sigma-q\Id)\phi \bigr ) & =\E[H_\alpha^\sigma(N_\sigma) (\LL_\sigma-q\Id) \phi(N_\sigma)] =\E[(\LL_\sigma-q\Id) H_\alpha^\sigma(N_\sigma) \phi(N_\sigma)] \\
        &= -(\lambda(\alpha)+q)\alpha !\, c_\alpha^\sigma (\phi)\,.
    \end{align*}
    The result follows since $\lambda(\alpha) + q \geq r+q>0$.
\end{proof}

\begin{definition}\label{definition:proj_d}
For a positive integer $m$ and $\phi$ an element of $L^{2}(\gamma_\sigma)$, we define the function $\projgeq{m}^\sigma(\phi):\R^d \rightarrow \R$ by  
\begin{align}\label{eq:proj_d}
     \projgeq{m}^\sigma(\phi) := \phi - \sum_{\alpha \in \Lambda_{m}} c_\alpha^\sigma(\phi) H_\alpha^\sigma\,,
\end{align}
where $\Lambda_m \coloneq \bigl \{ \alpha \in \N^d \; / \; |\alpha|\leq m-1 \bigr \}$. It is the projection of $\phi$ on the orthogonal space spanned by the family $(H_\alpha^\sigma)_{\alpha \in \Lambda_m}$.
\end{definition}

\noindent The operator $\projgeq{m}^{\sigma}$ enjoys the following continuity and stability property.  

\begin{lemma}\label{lemma:continuity_proj_d}
Let $m \geq 1$ be fixed. Consider a sequence $(\phi_n)_{n \geq 1}$ of functions in $\Cpold{0}{d}$ converging pointwise to a function $\phi$, and assume that there exist a constant $C>0$ and an integer $d \geq 0$ such that
\[
\forall n \in \N, \ \forall x \in \R^d, \quad |\phi_n(x)| \leq C(1+\|I_\sigma^{-1}x\|^l)\,.
\]
Then the sequence $\bigl( \projgeq{m}^\sigma(\phi_n) \bigr)_n$ converges pointwise to $\projgeq{m}^\sigma(\phi)$.
\end{lemma}

\begin{proof}
Fix $x \in \R^d$. By~\eqref{eq:proj_d}, one has
\[
\projgeq{m}(\phi_n)(x) = \phi_n(x) - \sum_{\alpha \in \Lambda_{m}} c_\alpha^\sigma(\phi_n) H_\alpha^\sigma(x)\,.
\]
By assumption, for every $\alpha \in \Lambda_{m}$,
\begin{align}\label{aux:continuity_proj_d}
    \forall n \in \N, \quad 
    |\phi_n(N_\sigma) H_\alpha^\sigma(N_\sigma)| \leq C (1+\|I_\sigma^{-1}N_\sigma\|^l) |H_\alpha^\sigma(N_\sigma)| \quad \text{a.s.}\,,
\end{align}
and the right-hand side of~\eqref{aux:continuity_proj_d} is integrable since $N_\sigma$ has finite moments of all orders.  
Moreover, $(\phi_n(N_\sigma))_n$ converges almost surely to $\phi(N_\sigma)$, so by the dominated convergence theorem,
\begin{align*}
c_\alpha^\sigma(\phi_n) = \frac{1}{\alpha!} \, \E[\phi_n(N_\sigma)H_\alpha^\sigma(N_\sigma)] \;\xrightarrow[n\to+\infty]{}\; \frac{1}{\alpha!}\, \E[\phi(N_\sigma)H_\alpha^\sigma(N_\sigma)] = c_\alpha^\sigma(\phi) \, .
\end{align*}
Therefore, for every $x \in \R^d$, 
\begin{align*}
    \projgeq{m}^\sigma(\phi_n)(x) \;\xrightarrow[n\to+\infty]{}\; \projgeq{m}^\sigma(\phi)(x)\,,
\end{align*}
which completes the proof.
\end{proof}

\begin{lemma}\label{lemma:estimates_proj_polynomial}
    Let $m$ be a positive integer. Let $\phi$ be an element of $\Cpold{0}{d}$ verifying
    \begin{align*}
        \forall x \in \R^d, \quad |\phi(x)| \leq A(1+\|I_\sigma^{-1}x\|^l)
    \end{align*}
    for some $A>0$ and some nonnegative integer $l$. Then there exists $B_{m,d}>0$ such that
    \begin{align*}
        \forall x \in \R^d, \quad |\projgeq{m}^\sigma(\phi)(x)| \leq AB_{m,l}\left (1+\|I_\sigma^{-1}x\|^{\max(l,m-1)} \right )\,,
    \end{align*}
    where $B_{m,d}$ does not depend on $\phi$ nor $\sigma$. Furthermore, for all $k \in \N$, if $\phi \in \Cpold{k}{d}$ then $\projgeq{m}^\sigma(\phi)\in \Cpold{k}{d}$. 
\end{lemma}
\begin{proof}
    Let $x$ be an element of $\R^d$. Again, by definition, one has
    \begin{align*}
        \projgeq{m}^\sigma(\phi)(x) = \phi(x) - \sum_{\alpha \in \Lambda_m} \, c_\alpha^\sigma(\phi_n)H_\alpha^\sigma(x)\,.
    \end{align*}
    In particular, for all $k\in \N$, if $\phi \in \Cpold{k}{d}$ then $\projgeq{m}^\sigma(\phi)\in \Cpold{k}{d}$, because $\projgeq{m}^\sigma(\phi)$ is a sum of $\phi$ with a polynomial.
    \begin{align*}
        \bigl |\projgeq{m}^\sigma(\phi)(x) \bigr | & \leq |\phi(x)| + \sum_{\alpha \in \Lambda_m} \, |c_\alpha^\sigma(\phi_n)||H_\alpha^\sigma(x)|\\
        & \leq A(1+\|I_\sigma^{-1}x\|^l) + \sum_{\alpha \in \Lambda_m} \, \frac{A}{\alpha!}\,\E[(1+\|I_\sigma^{-1}N_\sigma \|^l)|H_\alpha^\sigma(N_\sigma)|]|H_\alpha^\sigma(x)|\,,
    \end{align*}
    which yields the result, since $H_\alpha^\sigma$ is a polynomial of degree $|\alpha|$ and $I_\sigma^{-1}N_\sigma$, $H_\alpha^\sigma(N_\sigma)$ does not depend on $\sigma$.
\end{proof}
\begin{remark}\label{rmk:estimates_polynomial_d}
Above, we used the fact that a polynomial is a function with at most polynomial growth in the sense of our definition. The --generalized-- inequality of arithmetic and geometric means yields, for $\alpha \neq 0$,
\begin{align*}
|x^\alpha | \leq \left (\frac{\alpha_1x_1^2+\cdots+\alpha_dx_d^2}{|\alpha|} \right )^{|\alpha|/2} \leq \left ( \frac{\max(\alpha_1, \ldots, \alpha_d)}{|\alpha|} \right )^{|\alpha|/2}\|x\|^{|\alpha|}\,.
\end{align*}
Then, using the triangular inequality and the fact that for $0\leq n \leq q$ one has 
\[\forall x \in \R^d, \quad \|x\|^n \leq 1+\|x\|^q,\]
we can bound $|P(x)|$ by $C\left (1+\|x\|^{\deg(P)} \right )$ where the constant $C$ depends on the coefficients of $P$ and its degree.
\end{remark}
\begin{lemma}\label{lemma:conv_and_hermiterank_d}
    Let $(\phi_n)_n$ be a sequence of elements of $\Cpold{0}{d}$ of Hermite rank $r$, converging pointwise to a function $\phi\in L^2(\gamma_\sigma)$, and verifying
    \begin{align*}
        \forall n \in \N, \;\forall x \in \R^d, \quad |\phi_n(x)| \leq A(1+\|I_\sigma^{-1}x\|^l)\,,
    \end{align*}
    for some $A>0$ and some nonnegative integer $l$ which do not depend on $n$. Then $\rank_\sigma(\phi) \geq r$.
\end{lemma}
\begin{proof}
    Let $\alpha \in \N^d$. One has 
    \begin{align}\label{aux1:conv_and_hermiterank_d}
        c_\alpha^\sigma(\phi_n) \xrightarrow[n\to+\infty]{} c_\alpha^\sigma(\phi)\,.
    \end{align}
    Indeed, by assumptions, one has $\phi_n(N_\sigma) \xrightarrow[n\to+\infty]{} \phi(N_\sigma) \; \text{ a.e.}$, and 
    \begin{align}\label{aux2:conv_and_hermiterank_d}
        \forall n \in \N, \quad |H_\alpha^\sigma(N_\sigma)\phi_n(N_\sigma)| \leq A(1+\|I_\sigma^{-1}N_\sigma\|^d)|H_\alpha^\sigma(N_\sigma)|\,,
    \end{align}
    where the right-hand side of \eqref{aux2:conv_and_hermiterank_d} is integrable, since $N_\sigma$ has moments of all order and $H_\alpha^\sigma$ is a polynomial. By Lebesgue's dominated convergence theorem, one gets 
    \begin{align*}
        \E[H_\alpha^\sigma(N_\sigma)\phi_n(N_\sigma)] \underset{n\rightarrow +\infty}{\longrightarrow} \E[H_\alpha^\sigma(N_\sigma)\phi(N_\sigma)]\,.
    \end{align*}
    that is \eqref{aux1:conv_and_hermiterank_d}. To conclude, one has to notice that, in light of \eqref{aux1:conv_and_hermiterank_d}, $c_\alpha^\sigma(\phi)$ is zero when $c_\alpha^\sigma(\phi_n)$ is.
\end{proof}

\subsubsection{The modified Ornstein-Uhlenbeck semigroup}\label{ss:semigroup_d}
As in dimension one, the operator $\LL_\sigma$ is the infinitesimal generator of a semigroup $(\Pt_t^\sigma)_{t\geq0}$ on $L^2(\gamma_\sigma)$, and similarly it verifies the following Mehler's formula.

\begin{proposition}\label{proposition:mehler_formula_d}
Let $\phi$ be an element of  $L^2(\gamma_\sigma)$. For every $t\geq0$ and $x\in \R^d$, one has 
\begin{align}\label{eq:mehler_formula_d}
    \Pt_t^\sigma \phi(x) 
    & \coloneq \int_{\R^d} \phi \left(U_t x + M_ty\right) \, \gamma_\sigma(\dint y) = \E \left[ \phi \left(U_t x+ M_t N_\sigma \right) \right]\,,
\end{align}
where $U_t \coloneq \operatorname{diag} \left ( e^{-t}, \ldots, e^{-dt} \right )$ and $M_t \coloneq \operatorname{diag} \left ( \left (1-e^{-2t}\right )^{1/2}, \ldots, \left (1-e^{-2dt} \right )^{1/2} \right )$.
\end{proposition}
\noindent We refer to \cite{BGL2013} for more details. As in the one-dimensional case, we have the following.
\begin{proposition}
The family $(\Pt^\sigma_t)_{t\geq0}$ is a strongly continuous semigroup on $L^2(\gamma_\sigma)$ with infinitesimal generator $\LL_\sigma$. For every $\phi \in L^2(\gamma_\sigma)$, one has:
\begin{enumerate}
    \item $\displaystyle \lim_{t \to 0} \Pt_t^\sigma \phi = \phi$ ;
    \item $\displaystyle \lim_{t \to +\infty} \Pt^\sigma_t \phi = \int_{\R^d} \phi(x) \, \gamma_\sigma(\dint x)$ ;
    \item if $\phi$ is bounded, then for all $t \geq 0$, $\Pt^\sigma_t\phi$ is bounded and 
    \[
    \left \|\Pt_t^\sigma \phi \right \|_\infty \leq \left \|\phi \right \|_\infty \,;
    \]
    \item for all $x \in \R^d$ and $t>0$, 
    \begin{align}\label{eq:kolmogorov}
        \partial_t \Pt^\sigma_t \phi(x) = \LL_\sigma(\Pt^\sigma_t \phi)(x)\,;
    \end{align}
    \item for all $t \geq 0$, 
    \[
    \Pt^\sigma_t \phi = \sum_{\alpha \in \N^d}e^{-\lambda(\alpha)t} c_\alpha^\sigma (\phi) H_\alpha^\sigma\,,
    \]
    where the series converges in $L^2(\gamma_\sigma)$ ;
    \item The operator $\LL_\sigma$ has domain
    \begin{align*}
    \Dom(\LL_\sigma) \coloneq \left \{ f \in L^2(\gamma_\sigma) \;\big  / \; \underset{t \rightarrow 0}{\lim} \, \frac{\Pt_t^\sigma f - f}{t} \text{ exists in } L^{2}(\gamma_\sigma) \right \}\,.
    \end{align*}
\end{enumerate}
\end{proposition}

\noindent Let $V$ be a non-zero polynomial and let $\phi$ be an element of $\Cpol{d}$. We are interested in quantities of the form
\begin{align*}
    \Pt^{\sigma,V}_t \phi( x) \coloneq \E \left [ V(I_\sigma^{-1} N_\sigma) \phi \left (U_t x+M_t N_\sigma \right ) \right ], \quad x \in \R^d\,,
\end{align*}
where $U_t=\operatorname{diag}\left (e^{-t},\ldots, e^{-dt} \right )$, $M_t = \operatorname{diag}\left (\left (1-e^{-2t} \right)^{1/2},\ldots, \left (1-e^{-2dt} \right)^{1/2} \right )$.
\begin{lemma}\label{lemma:continuity_PtV}
Let $(\phi_n)_{n\in \N}$ be a sequence of elements of $\Cpold{0}{d}$ converging pointwise to a function $\phi$. Suppose that there exists $C>0$ and $l \in \N$ such that
\begin{align*}
\forall n \in \N, \; \forall x \in \R^d, \quad |\phi_n(x)| \leq C \left ( 1 + \|I_\sigma^{-1} x \|^l \right ).
\end{align*}
Then, for every $t \in [0,+\infty[$, the sequence $(\Pt_t^{\sigma,V} \phi_n)_{n\in\N}$ converges pointwise to $\Pt^{\sigma,V} \phi$.
\end{lemma}
\begin{proof}
Let $x$ be an element of $\R^d$. Since $(\phi_n)$ converges pointwise, one has 
\begin{align*}
V \left (I_\sigma^{-1} x \right ) \phi_n \left ( U_t x + M_t N_\sigma \right ) \xrightarrow[n\rightarrow +\infty]{} V \left (I_\sigma^{-1} x \right ) \phi \left ( U_t x + M_t N_\sigma \right ).
\end{align*}
Furthermore, one has the domination
\begin{align*}
\left | V \left (I_\sigma^{-1} x \right ) \phi_n \left ( U_t x + M_t N_\sigma \right ) \right | \leq C\left |V \left (I_\sigma^{-1} N_\sigma \right ) \right | \left ( 1+ \| I_\sigma^{-1}U_tx+I_\sigma^{-1}N_\sigma \|^l \right ),
\end{align*}
hence the result, by Lebesgue's dominated convergence theorem.
\end{proof}
\begin{proposition}\label{prop:derive_PtV}
   Let $\beta$ be an element of $\N^d$ such that $|\beta| = \sum_{p=1}^d \beta_p >0$. Then
\begin{align}\label{eq:derive_PtV}
    \partial_x^{\beta}  \Pt^{\sigma,V}_t \phi(x) =  \frac{(-1)^{|\beta|}e^{-\lambda(\beta)t}}{\prod_{p=1}^d \sigma_p^{\beta_p}\left ( 1-e^{-2pt}\right )^{\beta_p/2}}\Pt_t^{\sigma,A_\beta(V)}\phi(x)\,,
\end{align}
for every $x \in \R^d$, where 
\begin{align}\label{def:AbetaV}
    A_\beta(V)(u) \coloneq \sum_{\alpha\leq\beta} (-1)^{|\alpha|}\binom{\beta}{\alpha} \partial^{\beta-\alpha}V \left ( u \right )H_\alpha^\sigma \left ( I_\sigma u \right )\,.
\end{align}
\end{proposition}
\begin{proof}
    Let $x$ be an element of $\R^d$. Using the change of variable $u=M_t y$, one has
    \begin{align*}
        \Pt^{\sigma,V}_t \phi(x) & = \int_{\R^d} V(I_\sigma^{-1}y)\phi \left (U_tx+M_ty \right ) \, \gamma_\sigma( \dint y) \\
        & = \int_{\R^d} V\left ( I_\sigma^{-1}M_t^{-1}u \right )\phi \left (U_tx+u \right ) \frac{e^{-\frac{1}{2} \langle u,  I^{-2}_\sigma M_t^{-2}u \rangle} }{(2\pi)^{d/2}|\det(M_tI_\sigma)|}\, \dint u
       = \tilde{\phi} \ast g_{t,V} \left ( -U_t x \right )\,,
    \end{align*}
    where $\tilde{\phi}(u) \coloneq \phi(-u)$, and $g_{t,V}(u) \coloneq V\left ( I_\sigma^{-1}M_t^{-1}u \right )\rho_t^\sigma(u)$ with $\rho_{t}^\sigma(u) \coloneq \frac{e^{-\frac{1}{2} \langle u,I_\sigma^{-2} M^{-2}_tu \rangle }}{(2\pi)^{d/2}|\det(M_tI_\sigma)|}$. The function $\rho_t^\sigma$ is in $\mathscr{C}^\infty(\R^d)$ and one has
    \begin{align*}
        \forall \alpha \in \N^d, \; \forall u \in \R^d,  \quad \partial^\alpha_u \rho_t^\sigma(u) = (-1)^{|\alpha|} \left( \prod_{p=1}^d \sigma_p^{-\alpha_p}M_{t,p}^{-\alpha_p} \right )H_\alpha^\sigma \left ( M_t^{-1}u \right )\rho_t^\sigma(u)\,,
    \end{align*}
    where $M_{t,p}$ is the diagonal coefficient of index $(p,p)$ of the matrix $M_t$, that is $M_{t,p} =\left (1-e^{-pt}\right)^{1/2}$. Therefore, $g_{t,V}$ is in $\mathscr{C}^\infty(\R^d)$, and using Leibniz's formula, for all $u \in \R^d$ one gets
    \begin{align}
    \partial_x^\beta g_{t,V}(u) & = \left( \prod_{p=1}^d \sigma_p^{-\beta_p} M_{t,p}^{-\beta_p} \right ) \left ( \sum_{\alpha\leq\beta} (-1)^{|\alpha|} \binom{\beta}{\alpha} \partial^{\beta-\alpha}V \left ( I_\sigma^{-1}M_t^{-1}u \right )H_\alpha^\sigma \left ( M_t^{-1}u \right ) \right )\rho^\sigma_t(u) \\
    &  = \left( \prod_{p=1}^d \sigma_p^{-\beta_p} M_{t,p}^{-\beta_p} \right ) A_\beta^\sigma(V) \left ( I_\sigma^{-1}M^{-1}_t u \right )\rho_t^\sigma(u)\,,
    \end{align}
    where $A_\beta^\sigma (V)(u) \coloneq \sum_{\alpha\leq\beta} (-1)^{|\alpha|}\binom{\beta}{\alpha} \partial^{\beta-\alpha}V \left ( u \right )H_\alpha^\sigma \left ( I_\sigma u \right )$. By property of the convolution product, we deduce that $\Pt_t^{\sigma,V} \phi$ is infinitely differentiable as well on $\R^d$. Furthermore, one has
    \begin{align*}
        \partial_x^\beta \Pt^{\sigma,V}_t \phi(x) & = (-1)^{|\beta|}\left ( \prod_{p=1}^d U_{t,p}^{\beta_p} \right ) \tilde{\phi} \ast \left ( \partial_x^\beta g_{t,V} \right ) \left( -U_t x \right ) \\
        & = \frac{(-1)^{|\beta|}}{\sigma^\beta} \left ( \prod_{p=1}^d \frac{U_{t,p}^{\beta_p}}{M_{t,p}^{\beta_p}} \right )\int_{\R^d} \phi \left ( U_t x+u \right ) A_\beta^\sigma(V) \left (I_\sigma^{-1}M_t^{-1}u \right ) \rho_t^\sigma(u) \, \dint u\,.
    \end{align*}
Using again the change of variable $y = M_t^{-1}u$ leads to
    \begin{align*}
        \partial_x^\beta \Pt^{\sigma,V}_t \phi(x) & = \frac{(-1)^{|\beta|}}{\sigma^\beta} \left ( \prod_{p=1}^d \frac{U_{t,p}^{\beta_p}}{M_{t,p}^{\beta_p}} \right )\int_{\R^d} \phi \left ( U_t x+M_t y \right ) A_\beta^\sigma(V) \left (I_\sigma^{-1}y \right ) \, \gamma_\sigma(\dint y) \\
        & = \frac{(-1)^{|\beta|}}{\sigma^\beta} \left ( \prod_{p=1}^d \frac{U_{t,p}^{\beta_p}}{M_{t,p}^{\beta_p}} \right ) \E\left [ A_\beta^\sigma(V)(I_\sigma^{-1}N_\sigma) \phi \left (U_t x + M_t N_\sigma\right ) \right ] \\
        &= \frac{(-1)^{|\beta|}}{\sigma^\beta} \frac{e^{-\lambda(\beta)t}}{\prod_{p=1}^d \left ( 1-e^{-2pt}\right )^{\beta_p/2}}\Pt_t^{\sigma,A_\beta^\sigma(V)}\phi(x)\,.
    \end{align*}
\end{proof}

\begin{remark}
Proposition \ref{prop:derive_PtV} will mostly be used in the case $V=1$, in which case $A_\beta(V)(u) = (-1)^{|\beta|} H_\beta^\sigma(I_\sigma u)$.
\end{remark}

\begin{lemma}\label{lemma:estimates_PtV_pol}
    Let $V$ be a non-zero polynomial, let $\beta$ be an element of $\N^d$ and let $\phi$ be an element of $\Cpol{0}$ such that 
    \begin{align*}
        \forall x \in \R^d, \quad |\phi(x)| \leq C\left (1+\|I_\sigma^{-1}x\|^l \right )\,,
    \end{align*}
    where $C>0$ and $l\in\N$.
    Then, there exists $B_{l,V}>0$ such that
    \begin{align}\label{ineq:estimates_PtV_pol}
     \forall x \in \R^d, \quad\left |\Pt^{\sigma,V}_t\phi(x) \right | & \leq CB_{l,V} \left (1+\|I_\sigma^{-1}x\|^l\right )\,,
    \end{align}
    where $B_{V,l}$ does not depend on $\sigma$.
\end{lemma}
\begin{proof}
    Let $x$ be an element of $\R^d$. Using the triangular inequality, one has
    \begin{align*}
        \forall t \in \R_+^*, \quad \left |\Pt^{\sigma,V}_t\phi(x) \right | & \leq \E \left [ \left |V(I_\sigma^{-1}N_\sigma) \right |\left |\phi \left( U_t x + M_t N_\sigma \right) \right | \right ] \\
        & \leq C \, \E \left [ \left | V(I_\sigma^{-1}N_\sigma) \right | \right ]+ C \, \E \left [ \left |V(I_\sigma^{-1}N_\sigma) \right |\left \| U_tI_\sigma^{-1} x + M_tI_\sigma^{-1} N_\sigma \right \|^l  \right ] \\
        & \leq  C \, \E \left [ \left | V(I_\sigma^{-1}N_\sigma) \right | \right ]+ C \, \E \left [ \left |V(I_\sigma^{-1}N_\sigma) \right |2^{l-1}\left (\| U_t I_\sigma^{-1}x \|^l+ \|M_t I_\sigma^{-1}N_\sigma \|^l  \right )\right ] \\
        & \leq C B_{V,l}\left (1+\|I_\sigma^{-1}x\|^l \right )\,,
    \end{align*}
    where $B_{l,V} \coloneq 2^{l-1} \,\E \left [ \left | V(I_\sigma^{-1}N_\sigma ) \right | \left ( 1 + \|I_\sigma^{-1}N_\sigma\|^l \right ) \right ]$, since 
    \begin{align*}
        \forall t \in \R_+^*, \quad \|U_t\| \leq \underset{1\leq p\leq d}{\max} \,  e^{-pt} \leq 1\,, \quad \|M_t\|\leq \underset{1\leq p\leq d}{\max} \, \left (1-e^{-2pt} \right )^{1/2} \leq 1\,.
    \end{align*}
    The fact that $B_{l,V}$ does not depend on $\sigma$ follows from the fact that $I_\sigma^{-1}N_\sigma$ does not.
\end{proof}

\begin{proposition}\label{proposition:controle_semigroup_hermite_d}
    Let $V$ be a non-zero polynomial and let $\phi$ be an element of $\Cpold{0}{d}$ of Hermite rank $r$ such that 
    \begin{align*}
        \forall x \in \R^d, \quad |\phi(x)| \leq C\left (1+\|I^{-1}_\sigma x\|^l \right ),
    \end{align*}
    where $C>$ and $l\in\N$. Then there exists $B_{l,V,r}>0$ such that
    \begin{align}\label{eq:controle_semigroup_hermite_d}
        \left | \Pt_t^{\sigma,V} \phi(x) -\sum_{|\alpha|\leq r-1} c_\alpha \left (\Pt_t^{\sigma,V} \phi \right ) H^\sigma_\alpha( x)\right | \leq CB_{l,r,V}e^{-rt}\left ( 1 + \|I_\sigma^{-1} x\|^{l+r} \right ),
    \end{align}
    for every $x\in \R^d$. Furthermore, $B_{l,r,V}$ does not depend on $\sigma$.
\end{proposition}
\begin{proof}
Fix $t \geq 1$ and let $\beta$ be an element of $\N^d$ such that $|\beta|=r$. Consider the function 
\begin{align*}
    g :x \mapsto \Pt_t^{\sigma,V} \phi( x) -\sum_{|\alpha|\leq r-1} c_\alpha \left (\Pt_t^{\sigma,V} \phi \right ) H^\sigma_\alpha( x)\,.
\end{align*}
Then $\forall x \in \R^d, \quad \partial_x^\beta g(x) = \partial_x^\beta \left ( \Pt_t^{\sigma,V}\phi  \right )(x)$, and $g$ has Hermite rank greater than $r-1$.

Now, fix $\alpha \in \N^d$ such that $|\alpha|=r-1$. Let $x \in \R^d$ and let $y \in \R^d$, $y \neq x$. Consider the smooth function $f:s \mapsto \partial_x^\alpha g \bigl (y + s(x-y) \bigr  )$.  By the mean value theorem, there exists $s \in ]0,1[$ such that $f(1)-f(0) = f'(s)$, that is
\begin{align}\label{eq:mean_value_thm_partial^alpha_g}
    \partial_x^\alpha g(x)-\partial_x^\alpha g(y) &= \langle \nabla\left ( \partial_x^\alpha g \right) \bigl( y +s(x-y) \bigr ) ,x-y \rangle  = \sum_{p=1}^d \,(x_p-y_p)\partial_x^{\alpha+e_p}g\bigl( y +s(x-y) \bigr )\,.
\end{align}
Integrating \eqref{eq:mean_value_thm_partial^alpha_g} over $y$ with respect to the measure $\gamma_\sigma$, one gets
\begin{align}\label{eq:MV_thm_partial^alpha_g_integrated}
    \partial_x^\alpha g(x) &= \sum_{p=1}^d \,\int_{\R^d} \,(x_p-y_p)\partial_x^{\alpha+e_p}g\bigl( y +s(x-y) \bigr ) \, \gamma_\sigma(\dint y),
\end{align}
using that $\int_{\R^d} \partial_x^\alpha g(y) \, \gamma_\sigma(\dint y) = \int_{\R^d}  H_\alpha^\sigma(y) g(y) \, \gamma_\sigma(\dint y)=0$ since $g$ has Hermite rank greater than $r-1$ and $|\alpha|=r-1$. Substituting \eqref{eq:derive_PtV} in \eqref{eq:MV_thm_partial^alpha_g_integrated}, one gets
\begin{eqnarray}\label{eq:MV_thm_partial^alpha_g_integrated2}
    && \partial_x^\alpha g(x) = \sum_{p=1}^d \frac{x_p}{\sigma_p}\frac{(-1)^{r-1}e^{-\lambda (\alpha+e_p )t}}{\sigma^\alpha \prod_{q=1}^d \left (1-e^{-2qt} \right )^{(\alpha_q+\delta_{pq})/2}}\int_{\R^d}  \Pt_t^{\sigma,A^\sigma_{\alpha+e_p}(V)}\phi\bigl( y +s(x-y) \bigr ) \gamma_\sigma(\dint y) \nonumber \\
    &&\qquad \qquad + \sum_{p=1}^d \frac{(-1)^{r}e^{-\lambda (\alpha+e_p )t}}{\sigma^\alpha \prod_{q=1}^d \left (1-e^{-2qt} \right )^{(\alpha_q+\delta_{pq})/2}}\int_{\R^d}  \frac{y_p}{\sigma_p}\Pt_t^{\sigma,A^\sigma_{\alpha+e_p}(V)}\phi\bigl( y +s(x-y) \bigr ) \gamma_\sigma(\dint y)\,. 
\end{eqnarray}
On one hand, using Lemma \ref{lemma:estimates_PtV_pol},  there exist $B_{l,A^\sigma_{\alpha+e_1}(V)}, \ldots, B_{l,A^\sigma_{\alpha+e_d}(V)} >0$ such that 
\begin{align}\label{ineq:estimates_pol_PtV_proof_hermite_aux1}
    \forall p \in \llbracket 1,d \rrbracket, \forall w \in \R^d, \quad \left | \Pt_t^{\sigma, A^\sigma_{\alpha+e_p}(V)} \phi(w) \right | \leq CB_{l,A^\sigma_{\alpha+e_p}(V)} \left (1+\|I_\sigma^{-1} w\|^l \right )\,.
\end{align}
On the other hand, since $t\geq 1$, one has 
\[\forall p \in \llbracket 1,d \rrbracket, \quad \left (1-e^{-pt}\right )^{(\alpha_p+\delta_{pq})/2} \geq 2^{-(\alpha_p+\delta_{pq})/2}\,,\]
and since also $\lambda(\alpha+e_p) = \lambda(\alpha)+\lambda(e_p)\geq r-1+p \geq r$, one gets
\begin{align}\label{ineq:estimates_denom_derive_Pt}
    \frac{e^{-\lambda (\alpha+e_p )t}}{\sigma^\alpha \prod_{q=1}^d \left (1-e^{-2qt} \right )^{(\alpha_q+\delta_{pq})/2}} \leq \frac{2^{r/2}e^{-rt}}{\sigma^\alpha}\,.
\end{align}
Combining \eqref{ineq:estimates_pol_PtV_proof_hermite_aux1} and \eqref{ineq:estimates_denom_derive_Pt} in \eqref{eq:MV_thm_partial^alpha_g_integrated2}, one obtains
\begin{eqnarray}
    && \left |\partial_x^\alpha g(x) \right | \leq \sum_{p=1}^d \frac{|x_p|}{\sigma_p}\frac{2^{r/2}e^{-rt}}{\sigma^\alpha}\int_{\R^d}  \left | \Pt_t^{\sigma,A^\sigma_{\alpha+e_p}(V)}\phi\bigl( y +s(x-y) \bigr )\right | \gamma_\sigma(\dint y) \nonumber \\
    &&+ \sum_{p=1}^d \frac{2^{r/2}e^{-rt}}{\sigma^\alpha}\int_{\R^d}  \frac{|y_p|}{\sigma_p}\left |\Pt_t^{\sigma,A^\sigma_{\alpha+e_p}(V)}\phi\bigl( y +s(x-y) \bigr ) \right |\gamma_\sigma(\dint y)  \\
    && \leq CB_{l,\alpha,V}\frac{2^{r/2}e^{-rt}}{\sigma^\alpha}\sum_{p=1}^d\int_{\R^d} \left ( \frac{|x_p|}{\sigma_p}+\frac{|y_p|}{\sigma_p} \right )\left ( 1 +  \|  I^{-1}_\sigma y+sI^{-1}_\sigma(x-y) \|^l \right )  \gamma_\sigma(\dint y)\,, \nonumber
\end{eqnarray}
where $B_{l,\alpha,V} \coloneq \underset{1\leq p\leq d}{\max} \, B_{l,A^\sigma_{\alpha+e_p}(V)}$. Note that, with the exact constant obtained in the proof of Lemma \ref{lemma:estimates_PtV_pol}, each $B_{l,A^\sigma_{\alpha+e_p}(V)}$ does not depend on $\sigma$. Indeed, one has 
\begin{align*}
B_{l,A^\sigma_{\alpha+e_p}(V)} = 2^{l-1} \,\E \left [ \left | A^\sigma_{\alpha+e_p}(V)(I_\sigma^{-1}N_\sigma ) \right | \left ( 1 + \|I_\sigma^{-1}N_\sigma\|^l \right ) \right ]\,,
\end{align*}
where
\begin{align}\label{eq:A(V)_Nsigma}
A^\sigma_{\alpha+e_p}(V)(I_\sigma^{-1}N_\sigma) = \sum_{\beta \leq\alpha+e_p} \binom{\alpha+e_p}{\beta} \partial^{\alpha+e_p-\beta}V \left ( I_\sigma^{-1}N_\sigma \right )H_\beta^\sigma \left ( N_\sigma \right )\,.
\end{align}
Using the definition \eqref{eq:definition_hermite_polynomial_d}, one can see that each $H_\beta^\sigma \left ( N_\sigma \right )$ does not depend on $\sigma$. Then, using again the fact that $I_\sigma^{-1}N_\sigma$ does not depend on $\sigma$, we get that $B_{l,A^\sigma_{\alpha+e_p}(V)}$ does not depend on $\sigma$ either.
Using the change of variables $I_\sigma u = y$ in the integrals above, one gets
\begin{eqnarray*}
&& \left |\partial_x^\alpha g(x) \right | \leq CB_{l,\alpha,V}\frac{2^{r/2}e^{-rt}}{\sigma^\alpha}\sum_{p=1}^d\int_{\R^d} \left ( \frac{|x_p|}{\sigma_p}+|u_p| \right )\left ( 1 + \| u+s(I^{-1}_\sigma x- u)\|^l \right ) \, \gamma_d(\dint u)\,,
\end{eqnarray*}
where $\gamma_d$ denote the centered Gaussian measure on $\R^d$ with covariance matrix $I_d$. Using the fact that $0<s<1$, and $\|\sigma\| = 1$, we obtain that there exists a constant $A_{l,\alpha,V} >0$, which does not depend on $\sigma$, such that
\begin{align*}
    \left |\partial_x^\alpha g(x) \right | \leq CA_{l,r,V} \frac{e^{-rt}}{\sigma^\alpha} \left ( 1 + \|I_\sigma^{-1}x\|^l\sum_{p=1}^d \frac{|x_p|}{\sigma_p} \right )\,,
\end{align*}
for every $x \in \R^d$. The argument can be iterated.
\end{proof}
\begin{lemma}
    Let $\phi$ be an element of $\Cpold{0}{d}$ and let $V$ be a nonzero polynomial. Then 
    \begin{align*}
        \forall t \in \R_+, \quad \rank_\sigma( \Pt_t^{\sigma,V}\phi) \geq \rank_\sigma(\phi)-\deg(V) \, .
    \end{align*}
\end{lemma}
\begin{proof}
Let us write $r=\operatorname{rank}_{\sigma}(\phi)$ and $k=\deg (V)$. If $r\leq k$, the asserted lower bound is trivial, so assume $r>k$. Fix a multi-index $\alpha\in\mathbb{N}^d$ such that
$|\alpha|<r-k$. Let us show that  $c_\alpha^\sigma(P_t^{\sigma,V}\phi)=0$. We have
\begin{align*}
    \alpha!\,c_\alpha^\sigma(P_t^{\sigma,V}\varphi) & = \mathbb{E}\left[ H^\sigma_\alpha(N_\sigma) P_t^{\sigma,V}\varphi(N_\sigma) \right ] = \mathbb{E}\left [ H^\sigma_\alpha(N_\sigma)V(I_\sigma^{-1}Y_\sigma)
        \phi(U_tN_\sigma+M_tY_\sigma) \right],
\end{align*}
where $Y_\sigma$ is an independent copy of $N_\sigma$. Now, consider the Gaussian vectors
\begin{align*}
    Z := U_t N_\sigma + M_t Y_\sigma, \quad W := M_tN_\sigma-U_tY_\sigma \, .
\end{align*}
The matrices $U_t$ and $M_t$ are symmetric, commute, and satisfy
$U_t^2+M_t^2=I_d$. Consequently, we have
\begin{align*}
    \mathbb{E}[ZZ^\top] = \mathbb{E}[WW^\top] = I_\sigma^2, \quad \E[ZW^\top] = U_tM_t-M_tU_t = 0 \, .
\end{align*}
Thus $Z$ and $W$ are independent standard Gaussian vectors,  with common law $\gamma_\sigma$. Hence
\begin{align*}
    \alpha!\,c_\alpha^\sigma(P_t^{\sigma,V}\varphi) & = \mathbb{E}\left [ H^\sigma_\alpha(U_tN_\sigma+M_tY_\sigma)V(I_\sigma^{-1}(M_tN_\sigma -U_tY_\sigma))
        \phi(N_\sigma) \right],
\end{align*}
Using the independence of $N_\sigma$ and $Y_\sigma$, we obtain
\begin{align*}
    \alpha!\,c_\alpha^\sigma
    \bigl(P_t^{\sigma,V}\phi\bigr) = \mathbb{E}\left [\phi(N_\sigma)Q_{\alpha,t}^\sigma(N_\sigma)\right],
\end{align*}
where
\begin{align}\label{aux:integral}
    Q_{\alpha,t}^\sigma(z) := \int_{\mathbb{R}^d} H_\alpha^\sigma(U_tz+M_ty) V(I_\sigma^{-1}(M_tz-U_ty)) \,\gamma_\sigma(\dint y) \,.
\end{align}
The integrand in \eqref{aux:integral} is a polynomial in $(z,y)$ of total degree at most $|\alpha|+k$. Indeed, composition with a linear map does not increase polynomial degree, and the two factors have degrees at most $|\alpha|$ and $k$, respectively. Integrating with respect to $y$ replaces each monomial in $y$ by its Gaussian moment and does not increase the degree in $z$.
Consequently, $Q_{\alpha,t}^\sigma$ is a polynomial satisfying
\begin{align*}
    \deg (Q_{\alpha,t}^\sigma) \leq |\alpha|+k<r\, .
\end{align*}
Thus, by Lemma \ref{lemma:hermite_rank_d}, it follows that
\begin{align*}
    \alpha!\,c_\alpha^\sigma
    \bigl(P_t^{\sigma,V}\phi\bigr) = \mathbb{E}\left [\phi(N_\sigma)Q_{\alpha,t}^\sigma(N_\sigma)\right]=0 \, .
\end{align*}
This concludes the proof.
\end{proof}
\begin{remark}
    In particular, when $\phi$ verifies $\rank_\sigma(\phi) \geq r + \deg(V) $, for $V \neq 0$, above in \eqref{eq:controle_semigroup_hermite_d}, one has 
    \begin{align}\label{eq:controle_semigroup_hermite_final_d}
         \forall t \geq 1, \forall x \in \R^d, \quad \left|\Pt_t^{\sigma,V}\phi(x)  \right | \leq e^{-rt}CB_{l,r,V} \left (1+\|I_\sigma^{-1} x\|^{l+r} \right )\,.
    \end{align}
\end{remark}

We end this section on the Ornstein-Uhlenbeck semigroup with a commutation relation.

\begin{lemma}\label{lemma:commutation_semigroup_d}
Let $\phi$ be an element of $\Cpold{1}{d}$. Then, for all $t \in \R_+$, $\Pt_t \phi$ is differentiable and 
\begin{align*}
\forall p \in \llbracket 1, d \rrbracket, \quad \partial_p \Pt_t^\sigma \phi = e^{-tp} \Pt_t^\sigma \partial_p \phi\,.
\end{align*}
\end{lemma}

\subsubsection{Resolvent operator}\label{ss:resolvent_d}

The operator $\LL_\sigma$ has spectrum the set of nonpositive integers $\Z_-$. Its resolvent operator is defined by 
\[
\RR_\sigma(\xi ) := (\LL_\sigma-\xi \Id)^{-1}\,,
\]
which is a bounded operator on $L^2(\gamma_\sigma)$ for every $\xi \in \R \setminus \Z_-$. As it is the case for $\LL$, for $\xi \in \Z_-$, the operator $\LL_\sigma$ can also be inverted on suitable subspaces of $L^2(\gamma_\sigma)$ in the same way. For further background on resolvent operators of $\LL_\sigma$, we refer to \cite{BGL2013} and \cite{EN1999}. The next proposition collects the key properties of the inversion of $\LL_\sigma$, 
We include a proof for completeness.

\begin{proposition}\label{proposition:resolvent_d}
Let $\xi \in \R$ and let $\phi \in \Cpold{0}{d}$ such that $\rank_\sigma(\phi) > -\xi$. Then the function 
\begin{align*}
u: \begin{pmatrix}
\, \R^d & \longrightarrow & \R \, \\
\, x & \longmapsto & -\int_0^{+\infty} e^{-\xi t} \Pt^\sigma_t \phi(x)\,\dint t\,
\end{pmatrix}\,,
\end{align*}
is well-defined and continuous on $\R^d$. Furthermore, it has the following properties:
\begin{enumerate}[(i)]
\item $u \in \Dom (\LL_\sigma)$ and for every $x \in \R^d$, the limit
\begin{align}
\underset{s \rightarrow 0}{\lim}\, \frac{\Pt_s^\sigma u(x) -u(x)}{s}
\end{align}
exists. We denote this limit by $\LL_\sigma u(x)$.
\item for every $x \in \R^d$, $\left (\LL_\sigma-\xi \Id \right )u(x) = \phi(x)$. We will denoted by $\RR_\sigma(\xi) \phi$ the function $u$.
\end{enumerate}
\end{proposition}
\begin{proof}
Denote $r$ the quantity $\rank_\sigma(\phi)$. Consider the function $u:\R^d \rightarrow \R$ defined by 
\begin{align*}
\forall x \in \R^d, \quad u(x) \coloneq -\int_0^{+\infty} e^{-\xi t} \Pt^\sigma_t \phi(x)\,\dint t\,.
\end{align*}
By Lemma \ref{lemma:estimates_PtV_pol} and Proposition \ref{proposition:controle_semigroup_hermite_d}, since $\phi \in \Cpold{0}{d}$ and $r > -\xi$, there exists a constant $B = B_{\phi,r} >0$ and an integer $l\geq 0$ such that
\begin{align}\label{aux1:resolvent}
\left |\Pt_t^\sigma \phi(x) \right | \leq B \left ( 1 + \|I_\sigma^{-1} x\|^l \right )\I_{[0,1]}(t) + B e^{-rt}\left ( 1 + \|I^{-1}_\sigma x\|^{l+r} \right )\I_{[1,+\infty[}(t)\,,
\end{align}
for every $x \in \R^d$ and $t \in \R_+$. It follows that $u$ is well-defined and continuous on $\R^d$. Let $s>0$. One has, by Fubini,
\begin{align*}
\Pt_s^\sigma u(x) & = -\int_0^{+\infty} e^{-\xi t} \Pt_{t+s}^\sigma \phi(x) \, \dint t = -e^{\xi s} \int_{s}^{+\infty} e^{-\xi t} \Pt_t^\sigma \phi(x) \, \dint t \\
& = e^{\xi s} \left ( \int_{0}^{s} e^{-\xi t}\Pt_t^\sigma \phi(x) \, \dint t+u(x) \right )\,.
\end{align*}
Then, one has
\begin{align*}
\frac{\Pt_s^\sigma u(x) -u(x)}{s} = \frac{e^{\xi s}}{s} \int_{0}^{s} e^{-\xi t}\Pt_t^\sigma \phi(x) \, \dint t+\frac{e^{\xi s}-1}{s} u(x)\,,
\end{align*}
and taking the limit when $s$ goes to $0$, one gets 
\begin{align*}
\forall x \in \R^d, \quad \underset{s \rightarrow 0}{\lim}\, \frac{\Pt_s^\sigma u(x) -u(x)}{s} = \phi(x)+\xi u(x)\,.
\end{align*}
Hence $u \in \Dom \LL_\sigma$ and $\forall x \in \R^d, \quad \LL_\sigma u(x) = \phi(x)+\xi u(x)$, that is, 
\[\left ( \LL_\sigma-\xi \Id\right )u(x) = \phi(x)\,.\]
\end{proof}

\begin{proposition}
Let $\xi \in \R$ and let $\phi \in \Cpold{0}{d}$ such that $r=\rank_\sigma(\phi) > -\xi$. Then one has 
\begin{align}\label{eq:resolvent_spectral_d}
\RR_\sigma(\xi) \phi = -\underset{\alpha\in\N^d, \, |\alpha| \geq r}{\sum} \, \frac{c_\alpha^\sigma(\phi)}{\lambda(\alpha)+\xi}H_\alpha^\sigma \, ,
\end{align}
where the serie in \eqref{eq:resolvent_spectral_d} converges in $L^2(\gamma_\sigma)$.
\end{proposition}

\noindent A consequence of the latter property is the next lemma about the conservation of Hermite rank.

\begin{lemma}\label{lemma:conservation_rank_d}
  Let $\phi \in L^2(\gamma_\sigma)$ with Hermite rank $r$. Then $\RR_\sigma(\xi)\phi$ has Hermite rank $r$, for every $\xi>-r$.
\end{lemma}
\begin{proof}
    Let $\xi\in \bigl (\R \setminus \Z \bigr ) \, \cup ]-r,0]$. Then, using equality \eqref{eq:resolvent_spectral_d}, one deduces 
    \begin{align*}
    \forall \alpha \in \N^d, \quad c_\alpha^\sigma \left ( \RR_\sigma(\phi) \right ) = -\frac{c_\alpha^\sigma(\phi)}{\lambda(\alpha)+\xi}\,,
    \end{align*}
    hence the result.
\end{proof}

\begin{corollary}\label{corollary:prod_resolvent_conservation_rank_d}
    Let $\xi_1, \ldots, \xi_l$ be real numbers, and denote by $Q$ the rational fraction $\frac{1}{(X-\xi_1)\cdots (X-\xi_l)}$. For all $\phi \in \Cpold{0}{d}$ such that $\rank_\sigma(\phi)> -\underset{1\leq i \leq l}{\min} \xi_i$, 
   \[\rank_\sigma \left (Q(\LL_\sigma)\phi \right ) = \rank_\sigma \left ( \phi \right ) \, .\]
\end{corollary}
\begin{proof}
    By definition, $Q(\LL_\sigma) = \RR_\sigma(\xi_1) \cdots \RR_\sigma(\xi_l)$. By applying Lemma \ref{lemma:conservation_rank_d} $r$ times, one gets the result wanted.
\end{proof}
\noindent The operator $\RR_\sigma(\xi)$ satisfies a commutation relation.
\begin{proposition}\label{proposition:resolvent_commutation_d}
    Let $\xi$ be an element of $\R$. For every $\phi$ in $\Cpold{1}{d}$ such that  $\rank_\sigma(\phi) > -\xi$, $\RR_\sigma(\xi)\phi$ is differentiable and
\begin{align}\label{eq:resolvent_commutation}
     \forall p \in \llbracket 1,d\rrbracket, \quad \partial_{p} \bigl ( \RR_\sigma(\xi) \phi \bigr )& =\RR_\sigma(\xi+p) \, \partial_{p} \phi\,.
\end{align}
\end{proposition}
\begin{proof}
    Let $\phi$ be an element of $\Cpold{1}{d}$ such that $r = \rank_\sigma(\phi)>-\xi$, and fix $p\in\llbracket 1,d\rrbracket$. Then $\partial_p \phi \in \Cpold{0}{d}$ and $\partial_p \phi$ has Hermite rank $r-1$. Using Proposition \ref{proposition:controle_semigroup_hermite_d}, there exist a constant $B = B_{\partial_p \phi, r}>0$ and $l\in \N$ such that 
\begin{align*}
\left |\Pt_t^\sigma (\partial_p \phi)(x) \right | \leq B \left ( 1 + \|I_\sigma^{-1} x\|^l \right )\I_{[0,1]}(t) + B e^{-(r-1)t}\left ( 1 + \|I^{-1}_\sigma x\|^{l+r} \right )\I_{[1,+\infty[}(t)\,,
\end{align*}
for every $x \in \R^d$ and $t\in \R_+$. Therefore, using Lemma \ref{lemma:commutation_semigroup_d}, we get
\begin{align}\label{aux1:resolvent_commutation}
\left |\partial_p\Pt_t^\sigma \phi(x) \right | \leq B e^{-pt}\left ( 1 + \|I_\sigma^{-1} x\|^l \right )\I_{[0,1]}(t) + B e^{-(r-1+p)t}\left ( 1 + \|I^{-1}_\sigma x\|^{l+r} \right )\I_{[1,+\infty[}(t)\,.
\end{align}
Hence, using Proposition \ref{proposition:resolvent_d} and by differentiating under the integral sign, we obtain 
\begin{align*}
\partial_{p} \bigl ( \RR_\sigma(\xi) \phi \bigr )(x)&  = -\int_0^{+\infty} e^{-\xi t} \partial_p \Pt_t^\sigma \phi(x) \, \dint t = -\int_0^{+\infty} e^{-(\xi+p) t} \Pt_t^\sigma (\partial_p\phi)(x) \, \dint t \\
& =\RR_\sigma(\xi+p) \partial_p \phi(x)\,,
\end{align*}
for every $x\in \R^d$.
\end{proof}
\noindent We end the section with two lemmas: an estimate and a continuity property for resolvent operators of $\LL_\sigma$.
\begin{lemma}\label{lemma:estimates_resolvent_polynomial_d}
Let $\xi \in \R$. Let $\phi$ be a continuous function such that $\rank_\sigma(\phi)\geq r >-\xi$, verifying  $\forall x \in \R^d, \;|\phi(x)| \leq C\left (1+\|I_\sigma^{-1} x \|^l\right )$, for some $C>0$ and $l\in \N$. Then, there exists $K_{r,l,\xi} >0$ such that
\[\forall x \in \R^d, \quad
|\RR_\sigma(\xi) \phi(x) | \leq  CK_{r,l,\xi}\left (1+\|I_\sigma^{-1} x\|^{l+r} \right )\,.
\]
\end{lemma}
\begin{proof}
Let $x$ be an element of $\R^d$. Using Proposition \ref{proposition:resolvent_d}, one has
\begin{align*}
 \RR_\sigma(\xi)\phi(x) &= - \int_{0}^{+\infty} e^{-\xi t} \Pt^\sigma_{t} \phi(x) \, \dint t =-\int_{0}^{1} e^{-\xi t} \Pt_{t}^\sigma \phi(x) \, \dint t -\int_{1}^{+\infty} e^{-\xi t} \Pt^\sigma_{t} \phi(x) \, \dint t\,.
\end{align*}
Hence,  using Lemma \ref{lemma:estimates_PtV_pol} and \eqref{eq:controle_semigroup_hermite_final_d}, there exist $B_{l,1}>0$ and $B_{l,r,1}>0$ such that
\begin{align*}
| \RR_\sigma(\xi) \phi(x)|  &\leq C \left (B_{l,1}\int_0^1 e^{-\xi t} \, \dint t +B_{l,r,1}\int_{1}^{+\infty}  e^{-(\xi+r) t}\, \dint t \right )\left (1+\|I_\sigma^{-1}x\|^{l+r}\right ) \\
& \leq  CK_{r,l,\xi} \left (1+\|I_\sigma^{-1}x\|^{l+r} \right )\,,
\end{align*} 
where $K_{r,l,\xi} \coloneq   B_{l,1}\int_0^1 e^{-\xi t} \, \dint t +\frac{e^{-\xi-r}B_{l,r,1}}{\alpha+r}$.
\end{proof}

\begin{remark}
    In particular, from Propositions \ref{proposition:resolvent_commutation_d} and \ref{lemma:estimates_resolvent_polynomial_d}, one deduces that for every $\xi \in \R$, every integer $l>0$ and every $\phi \in \Cpold{l}{d}$ such that $\rank_\sigma(\phi) >-\xi$, $\RR_\sigma(\xi) \phi$ belongs to $\Cpold{l}{d}$.
\end{remark}

\begin{lemma}\label{lemma:continuity_resolvent_d}
    Let $\xi$ be an element of $\R$, and consider a sequence $(\phi_n)$ in $\Cpold{0}{d}$ such that $\rank_\sigma(\phi_n) > -\alpha$, converging pointwise to a function $\phi$. Suppose that there exist a constant $C$ and a nonnegative integer $l$ such that
    \[\forall n \in \N, \; \forall x \in \R^d, \quad |\phi_n(x)| \leq C\left (1+\|I_\sigma^{-1} x\|^l \right ).\]
    Then the sequence $(\RR_\sigma(\xi)\phi_n)_n$ converges pointwise to $\RR_\sigma(\xi)\phi$.
\end{lemma}
\begin{proof}
Fix $x$ in $\R$. Using Proposition \ref{proposition:resolvent_d} and then \eqref{eq:mehler_formula_d}, for all $n$ in $\N$, we have
\begin{align*}
    \RR_\sigma(\xi) \phi_n(x) & = \int_0^{+\infty} e^{-\xi t} \Pt^\sigma_t \phi_n(x) \, \dint t =  \int_0^{+\infty} e^{-\xi t} \, \E \left [\phi_n \left (U_tx+M_tN_\sigma \right ) \right  ] \, \dint t\,,
\end{align*}
where $N_\sigma \sim \mathcal{N}_d(0,I_\sigma^2)$. By assumption, one has 
\[\forall n \in \N, \; |\phi_n(xe^{-t}+\sqrt{1-e^{-2t}}N)| \leq C\left (1+\|U_tI^{-1}_\sigma x+M_tI^{-1}_\sigma N_\sigma\|^l \right )\,.\]
Since $N_\sigma$ admits moments of all order, Lebesgue's dominated convergence theorem yields
\begin{align*}
    \E \left [\phi_n \left (U_tx+M_tN_\sigma \right ) \right  ]\xrightarrow[n\to+\infty]{} 
   \E \left [\phi \left (U_tx+M_tN_\sigma \right ) \right  ] = \Pt^\sigma_t \phi(x)\,.
\end{align*}
Then, by \eqref{eq:controle_semigroup_hermite_final_d}, one has the domination 
\begin{align*}
    |e^{-\alpha t} \Pt_t \phi_n(x)| & \leq e^{-(\xi+r)t}CB_{l,r,V}(1+\|I_\sigma^{-1} x\|^{l+r})\,,
\end{align*}
where $r = \underset{n}{\min} \rank_\sigma(\phi_n)$ and $B_{r,l}$ is positive. Since $\xi+r >0$, the dominated convergence theorem leads to the result, by taking the limit inside the integral.
\end{proof}

\subsection{Analysis on the Wiener space}\label{s:analysis.wiener}

We work on the following countable product of probability spaces which we call a \emph{Wiener space}:
\begin{equation}\label{def:wiener-space}
  (\Omega,\mathscr{F},\prob) \coloneq \left(\mathbb{R},\mathscr{B}(\mathbb{R}),\gamma\right)^{\mathbb{N}}.
\end{equation}
We define the \emph{coordinate functions}
\begin{equation}\label{eq:coordinate.functions}
N_i \coloneq
\begin{cases}
  \Omega & \longrightarrow  \mathbb{R},
  \\ (x_{0}, x_{1},\cdots) & \longmapsto  x_i\,.
\end{cases}
\end{equation}
By construction, the $N_i$'s are independent random variables on $(\Omega,\mathscr{F},\prob)$, with common law the standard Gaussian distribution.

\subsubsection{Wiener chaoses}

  For $m \in \mathbb{N}$, the \emph{$m$-th Wiener chaos} $\mathcal{W}_{m}$ is defined as the $L^{2}(\prob)$-closure of the linear span of functions of the form $\prod_{i=0}^\infty H_{k_i}(x_i)$ where the $k_{i}$'s satisfy $\sum_{i=0}^\infty k_i=m$.
The above product is finite since $H_0(x)=1$ and only finitely many integers $(k_i)_{i\ge 0}$ are non zero.
For $m = 0$, we find that $\mathcal{W}_{0}$ is the linear space of constant functions.
Importantly, Wiener chaoses provide the following  decomposition.
\begin{theorem}[It\^o-Wiener decomposition]
    Let $m\neq n$ be two nonnegative integers. The linear subspaces $\WW_m$ and $\WW_n$ are orthogonal. Furthermore, we have the direct sum
    \begin{align*}
L^2(\prob)&=\bigoplus_{m=0}^\infty \mathcal{W}_m\,.
    \end{align*}
\end{theorem}
\noindent We sometimes work in a finite sum of Wiener chaoses.
Accordingly, let us define
\begin{equation*}
  \mathcal{W}_{\leq m} \coloneq \bigoplus_{k=0}^m \mathcal{W}_k\,.
\end{equation*}
Note that this exposition of the Wiener chaoses is only one among several. Another important way of introducing them is through the so-called Wiener-It\^o multiple integrals. We refer to \cite{NP2012, Nua2006} and Section \ref{sec.contraction} below for that approach.

  We often use that on $\mathcal{W}_{\leq m}$ all the $L^{p}(\prob)$-norms are equivalent.
Namely, we have the following result.
\begin{proposition}\label{prop:hypercontractivity}
    Let $m \in \mathbb{N}$, and $1 < p < q < \infty$. Then, for all $F \in \WW_{\leq m}$ such that $\E[F]=0$, one has 
\begin{align}\label{eq:hypercontractivity}
  \norm{F}_{p} & \leq \norm{F}_{q} \leq \left ( \frac{q-1}{p-1} \right )^{m/2} \norm{F}_{p}\,.
\end{align}
\end{proposition}

\noindent This fact is a well-known consequence of hypercontractivity estimates; see, for instance, \cite[Corollary~2.8.14]{NP2012}.

\subsubsection{$\L$ and $\Gamma$ operators}

We now introduce the operators from Malliavin calculus that will be used throughout this article.  
For a broader introduction and detailed proofs of the properties below, we refer to the textbooks \cite{BH1991,NP2012,Nua2006}.

\begin{definition}
The \emph{Ornstein--Uhlenbeck operator} $\L$ is defined by  
\begin{align}\label{eq:def_L}
    \L F \coloneq -\sum_{p=0}^{+\infty} p \, J_p(F)\,,
\end{align}
on the domain
\[
\operatorname{Dom} \L \coloneq \Bigl\{ F \in L^2(\prob) \;\big|\; \sum_{p=0}^{+\infty} \,p^2 \, \|J_p(F)\|_2^2 < +\infty \Bigr\},
\]
where $J_p(F)$ denotes the orthogonal projection of $F$ onto the $p$-th Wiener chaos $\WW_p$.
\end{definition} 

\begin{definition}
    The \emph{pseudo-inverse} $\L^{-1}$ of $\L$ is defined by
    \begin{align*}
        \forall F \in L^2(\prob), \quad \L^{-1} F \coloneq -\sum_{p=1}^{+\infty} \frac{1}{p} J_p(F) \, .
    \end{align*}
    For any $F \in L^2(\prob)$, one has $\L^{-1}F \in \Dom \L$ and $\L \L^{-1} F = F -\E[F]$.
\end{definition}

\begin{definition}
The \emph{carr\'e-du-champ operator} (or \emph{square field operator}) $\Ga$ associated with $\L$ is defined, for $\phi,\psi \in \mathscr{C}_{pol}^2(\R^m)$, $m \in \N^*$, by  
\begin{align}\label{eq:def_Gamma}
    \Gamma(\phi,\psi) \coloneq \sum_{k=1}^{m} \partial_{x_k}\phi \;\partial_{x_k}\psi\,.
\end{align}
It extends by closability to the domain
\[
\operatorname{Dom} \Ga \coloneq \Bigl\{ H \in L^2(\prob) \;\big|\; \sum_{p=0}^{+\infty} p \,\|J_p(H)\|_2^2 < +\infty \Bigr\}\,.
\]
In particular, for any $m \in \N$, one has $\WW_{\leq m} \subset \operatorname{Dom} \L \subseteq \operatorname{Dom} \Ga$.
\end{definition}

\noindent The following result records some classical properties.

\begin{proposition}\label{prop:LGonsmooth}
\begin{enumerate}[(i)]
    \item For all $p \in \N$, \quad $\WW_p = \operatorname{Ker}(\L+p\Id)$.
    \item For all $m \in \N^*$ and all $\phi \in \mathscr{C}_{pol}^2(\R^m)$,  
    \[
    \L \phi = \Delta \phi - \sum_{k=1}^m x_k \,\partial_{x_k}\phi\,.
    \]
\end{enumerate}
\end{proposition}

\noindent The next proposition summarizes some fundamental identities involving $\L$ and $\Ga$, which will be heavily used in the sequel.

\begin{proposition}\label{prop:LGamma}
Let $m \in \N$. For all $F,G \in \WW_{\leq m}$, the following hold:
\begin{enumerate}[(a)]
    \item \emph{Integration by parts:}  
    \begin{align}\label{eq:LIPP}
        \E[\Ga(F,G)] = -\E[G \L F] = -\E[F \L G]\,.
    \end{align}
    \item \emph{Diffusion property:} for any $\varphi \in \Cpol{2}$,  
    \begin{align}\label{eq:diffusion}
        \L \varphi(F) = \varphi''(F)\,\Ga(F,F) + \varphi'(F)\,\L F\,.
    \end{align}
    Equivalently, $\Ga$ acts as a derivation in the sense that  
    \[
    \Gamma(\varphi(F),F) = \varphi'(F)\,\Gamma(F,F)\,.
    \]
    \item \emph{Spectral stability:}  
    \begin{align}\label{eq:spectral_stability}
        FG \in \WW_{\leq 2m}.
    \end{align}
    In particular, if $m \geq 1$, then
    \begin{align}\label{eq:gamma_stability}
        \Ga(F,F) \in \WW_{\leq 2m-2}\,.
    \end{align}
\end{enumerate}
\end{proposition}
\noindent For brevity, we will often write $\Ga(F) \coloneq \Ga(F,F)$.  
Proofs of Propositions~\ref{prop:LGonsmooth} and~\ref{prop:LGamma} can be found in \cite{BH1991}.

As well as in the finite dimensional case, we define also a semigroup.

\begin{definition}
    The \emph{Ornstein-Uhlenbeck semigroup} $(P_t)_{t\geq 0}$ is defined by
    \begin{align*}
        \forall F \in L^{2}(\prob), \quad P_t(F) \coloneq \sum_{p=0}^{+\infty} e^{-pt} J_p(F) \, ,
    \end{align*}
    for every $t\geq 0$.
\end{definition}
It enjoys similar properties as $(\Pt_t)_{t\geq 0}$. We refer the reader to \cite[Section 2.8]{NP2012} for more details.

\subsubsection{Wiener chaoses and contractions}\label{sec.contraction}

We briefly explain how the Gaussian product space considered here fits into
the formalism of isonormal Gaussian processes used in \cite{NP2012}, in order to use some of the results used through the use of contraction. We refer to \cite{Janson1997,NP2012,Nua2006} for more details.
Recall the coordinate maps $N_k$ defined in \eqref{eq:coordinate.functions}. The
sequence $(N_k)_{k\geq 1}$ then consists of independent standard Gaussian
random variables.
Let $\mathfrak{H}=\ell^2(\N)$, endowed with its usual inner product, and for $h=(h_k)_{k\geq 1}\in\mathfrak{H}$,
define
\begin{align*}
    X(h) := \sum_{k= 1}^{+\infty} h_k N_k \, , 
\end{align*}
where the series converges in $L^2(\Omega)$. The family
\begin{align*}
    X = \bigl\{X(h):h\in\mathfrak{H}\bigr\}
\end{align*}
is an \emph{isonormal} Gaussian process over $\mathfrak{H}$. Indeed, it is a centered Gaussian family satisfying
\begin{align*}
   \forall (h,g)\in \mathfrak{H}^2, \quad \E\bigl[X(h)X(g)\bigr] = \langle h,g\rangle_{\mathfrak{H}}\,.
\end{align*}
Moreover, the Gaussian space generated by $X$ coincides with that generated by the coordinate functions $(N_k)_{k\geq 1}$. Thus, the Wiener space $\Omega$ can be identified with the Gaussian probability space associated with an isonormal Gaussian process over $\ell^2(\N)$. 

For every integer $q\geq 1$, we denote by $\mathfrak{H}^{\odot q}$ the $q$-th symmetric tensor power of $\mathfrak{H}$ and by $I_q(f)$ the multiple Wiener--It\^o integral of a kernel $f\in\mathfrak{H}^{\odot q}$. The previously defined $q$-th Wiener chaos $\WW_q$ can then be represented as 
\begin{align*} 
\WW_q = \left\{ I_q(f) : f\in\mathfrak{H}^{\odot q} \right\}\,. \end{align*} 
Thus, every element of $\WW_q$ admits a representation as a multiple Wiener--It\^o integral of order $q$ with a symmetric kernel in $\mathfrak{H}^{\odot q}$.

\subsubsection{Estimates through contractions}

Let $p,q$ be positive integers, let $f\in\mathfrak{H}^{\odot p}$ and $g\in\mathfrak{H}^{\odot q}$, and let $r\in\{0,\ldots,p\wedge q\}$. We denote by 
\begin{align*} 
f\otimes_r g \in \mathfrak{H}^{\otimes(p+q-2r)} 
\end{align*} 
the contraction of order $r$ of $f$ and $g$, and by $f\widetilde{\otimes}_r g$ its symmetrization. By convention, \begin{align*} 
f\otimes_0 g = f\otimes g \, , 
\end{align*} 
while, when $p=q$,
\begin{align*} 
f\otimes_p g = \langle f,g\rangle_{\mathfrak{H}^{\otimes p}}\,. 
\end{align*} 
We refer to \cite[Appendix~B.4]{NP2012} for the precise definition of contractions and their main properties. We now recall several results from \cite{NP2012} that will be used below. 
\begin{proposition}[Product formula]\label{prop:product_formula} 
Let $p,q$ be positive integers and let $F=I_p(f)$ and $G=I_q(g)$, where $f\in\mathfrak{H}^{\odot p}$ and $g\in\mathfrak{H}^{\odot q}$. Then 
\begin{align}\label{eq:product-formula}
FG = \sum_{r=0}^{p\wedge q} r! \binom{p}{r} \binom{q}{r} I_{p+q-2r} \left( f\widetilde{\otimes}_r g \right)\, . 
\end{align} 
\end{proposition} 
\noindent As a consequence of the product formula, the fourth cumulant of an element of a fixed Wiener chaos can be expressed in terms of contractions of its kernel. 
\begin{lemma}[Lemma 5.2.4, \cite{NP2012}]\label{lem:fourth_cumulant_contractions} 
Let $F$ be an element of $\WW_q$, with $q\geq 2$, and write $F=I_q(f)$, where $f\in\mathfrak{H}^{\odot q}$. Then 
\begin{align*} 
\E[F^4] - 3\E[F^2]^2 &= (q!)^2 \sum_{r=1}^{q-1} \binom{q}{r}^2 \Bigg( \left\|f\otimes_r f\right\|^2_{ \mathfrak{H}^{\otimes(2q-2r)} } + \binom{2q-2r}{q-r} \left\|f\widetilde{\otimes}_r f\right\|^2_{ \mathfrak{H}^{\otimes(2q-2r)} } \Bigg)\, . 
\end{align*} 
\end{lemma}
\noindent The following estimate controls the mixed terms associated with two random variables belonging to different Wiener chaoses. 
\begin{lemma}[Lemma 6.2.1, \cite{NP2012}]\label{lem:mixed_gamma_contractions} 
Let $p,q$ be positive integers and let $F=I_p(f)$ and $G=I_q(g)$, where $f\in\mathfrak{H}^{\odot p}$ and $g\in\mathfrak{H}^{\odot q}$. If $p<q$, then 
\begin{align*} 
\frac{1}{q^2} \E\left[ \Gamma(F,G)^2 \right] &\leq (p!)^2 \binom{q-1}{p-1}^2 (q-p)! \left\|f\right\|^2_{\mathfrak{H}^{\otimes p}} \left\|g\otimes_{q-p}g\right\|_{ \mathfrak{H}^{\otimes 2p} } \\ &\quad + \frac{p^2}{2} \sum_{r=1}^{p-1} (r-1)!^2 \binom{p-1}{r-1}^2 \binom{q-1}{r-1}^2 (p+q-2r)! \\ &\qquad\qquad\times \left( \left\|f\otimes_{p-r}f\right\|^2_{ \mathfrak{H}^{\otimes 2r} } + \left\|g\otimes_{q-r}g\right\|^2_{ \mathfrak{H}^{\otimes 2r} } \right)\, .
\end{align*} 
\end{lemma} 
\noindent The preceding estimate yields the following more convenient bound. 
\begin{lemma}\label{lemma:estimates_variance_pq} 
Let $p,q$ be positive integers, let $F$ be an element of $\WW_p$, and let $G$ be an element of $\WW_q$. If $p<q$, then there exists a constant $C_{p,q}>0$, depending only on $p$ and $q$, such that 
\begin{align*} 
\E\left[ \Gamma(F,G)^2 \right] \leq C_{p,q} \Bigg( \Var\bigl(\Gamma(F)\bigr) + \Var\bigl(\Gamma(G)\bigr) + \E[F^2] \Bigl( \Var\bigl(\Gamma(G)\bigr) \Bigr)^{1/2} \Bigg)\, . 
\end{align*} 
\end{lemma}

\subsection{Stein's method on Wiener chaos}\label{ss:steinsmethod}

The classical Stein's method starts from the following fundamental observation.

\begin{lemma}[Stein's lemma, 1972]
Let $\FF$ be the set of absolutely continuous functions $f \in L^2(\gamma)$ such that $f' \in L^1(\gamma)$, and let $X$ be a random variable. Then 
\begin{align}\label{eq:stein_charac}
    X \sim \mathcal{N}(0,1) 
    \quad \Longleftrightarrow \quad 
    \forall f \in \FF, \quad \E[f'(X)-Xf(X)]=0\,.
\end{align}
\end{lemma}

This characterization allows one to quantify the proximity of a random variable $X$ to the standard normal distribution by controlling the quantity $\E[f'(X)-Xf(X)]$.  
More precisely, let $\HH$ be a suitable class of measurable functions $h:\R \to \R$ (see \cite[Definition~C.1.2]{NP2012}) such that 
\begin{align}\label{eq:def_dHH}
    \dHH(F,G) \coloneq \sup_{h \in \HH} \, \bigl|\E[h(F)-h(G)]\bigr|
\end{align}
defines a distance on the space of Borel probability measures on $\R$.  
Important examples include:
\begin{enumerate}[(a)]
    \item $\HH = \{\I_A : A \in \mathscr{B}(\R)\}$, yielding the \emph{total variation} distance $\dHH = \dTV$;
    \item $\HH = \{\I_{(-\infty,z]} : z \in \R\}$, yielding the \emph{Kolmogorov} distance $\dHH = \dKOL$;
    \item $\HH = \operatorname{Lip}_1(\R)$, the set of $1$-Lipschitz functions, yielding the \emph{Wasserstein-1} distance $\dHH = \dWass$ on measures with finite first moment.
\end{enumerate}
The distance~\eqref{eq:def_dHH} can also be extended to bounded signed measures (see \cite{Dud2002}).  

\medskip

The core of Stein's method proceeds as follows.  
Fix $h \in \HH$. The \emph{Stein equation} associated with $h$ is the differential equation
\begin{align}\label{eq:stein_equation}
    \phi'(x) - x \phi(x) = h(x) - \E[h(N)]\,, \quad x \in \R\,,
\end{align}
where $N \sim \mathcal{N}(0,1)$.  
A particular solution $\phi_h$ is given by
\begin{align} \label{eq:stein_solll}
    \phi_h(x) \coloneq \sqrt{2\pi}\, e^{x^2/2} \int_{-\infty}^x \bigl(h(y)-\E[h(N)]\bigr)\,\gamma(\dint y)\,, \quad x \in \R\,.
\end{align}
For any random variable $F$, plugging $x=F$ into \eqref{eq:stein_equation} and taking expectations yields
\begin{align}\label{eq:auxh}
    \E[h(F)-h(N)] = \E[\phi'_h(F)-F\phi_h(F)].
\end{align}
Therefore,
\begin{align}\label{eq:boundingbystein}
    \dHH(F,N) 
    = \sup_{h \in \HH} \,|\E[h(F)-h(N)]|
    = \sup_{h \in \HH} \,|\E[\phi'_h(F)-F\phi_h(F)]|\,.
\end{align}

The essence of Stein's method is  to control the right-hand side of~\eqref{eq:boundingbystein}.  
To this end, bounds on the solution $\phi_h$ from \eqref{eq:stein_solll} to the Stein equation are crucial.  
In our setting, we are primarily interested in the total variation distance, i.e., $\dHH$ with $\HH$ given by the set of indicator functions of measurable subsets of $\R$.  
Since such functions are not regular, this creates technical difficulties.  
However, as shown in \cite{NP2015} and \cite[Proposition~C.3.5]{NP2012}, one may use a density argument to obtain
\begin{align}\label{eq:dTV_continuous}
    \dTV(F,G) = \frac{1}{2} \sup_{\substack{h \in \mathscr{C}_b^0(\R)\\ \|h\|_\infty \leq 1}} \,|\E[h(F)-h(G)]|\,,
\end{align}
for all random variables $F,G$.  
This allows one to improve the regularity of the associated solution $\phi_h$.  

The next proposition summarizes some useful properties of $\phi_h$ in this context (see \cite{NP2012} for proofs).

\begin{proposition}\label{prop:proprietes_solstein}
Let $h:\R \to \R$ be continuous and bounded, let $p > 0$, and let $F \in \WW_p$. Then:
\begin{enumerate}[(i)]
    \item $\phi'_h$ is continuous ;
    \item $\|\phi_h\|_\infty \leq \sqrt{\tfrac{\pi}{2}}\,\|h-\E[h(N)]\|_\infty$,  
    \quad $\|\phi'_h\|_\infty \leq 2\,\|h-\E[h(N)]\|_\infty$ ;
    \item $\E[h(F)-h(N)] = \tfrac{1}{p}\,\E[\phi'_h(F)\,(p-\Ga(F))]$.
\end{enumerate}
\end{proposition}

\section{Overview of proof}\label{sec.overview}

In this section, we fix $d>1$ and $F \in \WW_{\leq d}$ such that $\E[F]=0$ and $\E[F^2]=1$.

\medskip

We extend the proof of the corresponding result of non-asymptotic Edgeworth expansion for random variables $F$ belonging to a fixed Wiener chaos $\WW_p$ by deriving the multivariate counterpart of the expansions developed in \cite{MPS2025.EE}. 

Although our main focus is on the multivariate setting, the starting point of the proof of Theorem~\ref{thm:EE_Wiener_sums} is the same as that of \cite[Theorem 1.2]{MPS2025.EE}. Namely, consider the following representation of the total variation distance between the law of $F$ and the signed measure $\boldsymbol{\gamma}_{F,M}$, given in \eqref{eq:dTV_continuous}:
\begin{align}\label{eq:dTV_continuous_bis}
    \dTV(\prob_F,\boldsymbol{\gamma}_{F,M}) = \frac{1}{2} \sup_{\substack{h \in \mathscr{C}_b^0(\R)\\ \|h\|_\infty \leq 1}} \,\left | \int_\R h(x) \, \prob_F(\dint x) - \int_\R h(x) \, \boldsymbol{\gamma}_{F,M}( \dint x) \right |\,.
\end{align}
We can rewrite the right-hand side of \eqref{eq:dTV_continuous_bis} the following way.
\begin{lemma}\label{lemma:rpz-dTV-EE}
    Let $M$ be a nonnegative integer. Then
    \begin{align}\label{eq:rpz-dTV-EE}
         \dTV(\prob_F,\boldsymbol{\gamma}_{F,M}) = \frac{1}{2} \sup_{\substack{h \in \mathscr{C}_b^0(\R)\\ \|h\|_\infty \leq 1}} \E\left [ \projgeq{M-1}(\phi_h')(F) ( \Ga(F,-\L^{-1}F)-1)\right ]\,,
    \end{align}
    where $\phi_h$ is the solution to the Stein equation, defined in \eqref{eq:stein_solll}.
\end{lemma}
\begin{proof}
Denote by $g_{F,M}$ the function $1+\sum_{k=3}^{M} \tfrac{\E[H_k(F)]}{k!} H_k(x) $, and fix $h \in \mathscr{C}^0_b(\R)$ such that $\|h\|_\infty \leq 1$. We have
\begin{align}\label{dvp-integral}
    \int_\R h(x) \, \boldsymbol{\gamma}_{F,M}( \dint x) = \E[h(N)g_{F,M}(N)] = \E[h(N)]+\sum_{k=3}^{M} \frac{\E[H_k(N)h(N)]}{k!}\E[H_k(F)] \, .
\end{align}
Then, on one hand, by integrating by parts (see, e.g., \cite{BBNP2012}), we get
\begin{align*}
    \E[H_k(N)h(N)] = -k\,\E[H_{k-2}(N)\phi'_h(N)] \ ,
\end{align*}
where $\phi_h$ is the solution to the Stein equation defined in \eqref{eq:stein_solll}. On the other hand, using \eqref{eq:LIPP}, we get
\begin{eqnarray*}
    && \E [ H_k(F) (\Ga(F,-\L^{-1} F)-1)  ] = -\E  [ H_k(F)\Ga(F,\L^{-1} F)]  ]- \E [ H_k(F)  ] \\
    && =\frac{-1}{(k+1)} \E [ \Ga(H_{k+1}(F),\L^{-1} F)  ]- \E  [ H_k(F)  ] =\frac{1}{(k+1)} \E  [ H_{k+1}(F) F ]- \E [ H_k(F) ] =\frac{1}{(k+1)} \E [ H_{k+2}(F) ]\, ,
\end{eqnarray*}
hence $\E [ H_{k}(F) ]= (k-1) \, \E [ H_{k-2}(F) (\Ga(F,-\L^{-1} F)-1)  ]$. Plugging it back into \eqref{dvp-integral}, we obtain
\begin{align*}
    \int_\R h(x) \, \boldsymbol{\gamma}_{F,M}( \dint x) = \E[h(N)] - \sum_{k=1}^{M-2} \frac{\E[H_k(N)\phi'_h(N)]}{k!}\E [ H_k(F) (\Ga(F,-\L^{-1} F)-1)  ]
\end{align*}
and then
\begin{align*}
     \int_\R h(x) \, \prob_F(\dint x) - \int_\R h(x) \, \boldsymbol{\gamma}_{F,M}( \dint x) & = \E[h(F)-h(N)]- \sum_{k=1}^{M-2} \frac{\E[H_k(N)\phi'_h(N)]}{k!}\E [ H_k(F) (\Ga(F,-\L^{-1} F)-1)  ] \\
    & =  \E\left [ \left ( \phi_h'- \sum_{k=1}^{M-2} \frac{\E[H_k(N)\phi'_h(N)]}{k!}H_k \right)(F) ( \Ga(F,-\L^{-1}F)-1)\right ] \\ \\
    & =  -\E\left [ \projgeq{M-1}(\phi_h')(F) ( \Ga(F,-\L^{-1}F)-1)\right ] \, ,
\end{align*}
where $\projgeq{M-1}(\phi'_h)(x) = \phi_h'- \sum_{k=1}^{M-2} \frac{\E[H_k(N)\phi'_h(N)]}{k!}H_k$ is defined in \eqref{eq:proj_d} (take $d=1$ and $\sigma = 1$).
\end{proof}

\begin{remark}
    More generally, if $X$ is a random variable admitting moments up to order $M$, such that $\E[X^2]=1$, then it follows from the definition of $\boldsymbol{\gamma}_{X,M}$ that 
    \begin{align}\label{eq:general-rpz-EE}
    \dTV(\prob_X,\boldsymbol{\gamma}_{X,M}) = \frac{1}{2} \sup_{\substack{h \in \mathscr{C}_b^0(\R)\\ \|h\|_\infty \leq 1}} \,\left |\E[\projgeq{M+1}(h)(F)] \right |\,.
    \end{align}
   Equality~\eqref{eq:general-rpz-EE} thus reveals an interesting connection for Edgeworth expansions in integral probability metrics: adding $M+1$ correction terms to the Gaussian distribution $\gamma$ is equivalent to subtracting from $h$ the first $M+1$ terms of its expansion in the polynomial orthogonal basis of $L^2(\gamma)$.
\end{remark}
Hence, controlling the distance between $\prob_F$ and $\boldsymbol{\gamma}_{F,M}$ reduces to bounding 
\begin{align}\label{term-to-bound} 
\E\left[ \projgeq{M-1}(\phi_h')(F) \bigl(\Ga(F,-\L^{-1}F)-1\bigr) \right]  , 
\end{align} 
uniformly over $h$. Since increasing $M$ amounts to including additional correction terms, one expects the resulting bound to improve as $M$ increases. When $F$ belongs to a fixed Wiener chaos $\WW_p$, the authors of \cite{MPS2025.EE} prove that, for $M$ of the form $4m-1$, the quantity in \eqref{term-to-bound} can be bounded by 
\begin{align*} 
\bigl(\Var \Ga(F)\bigr)^{\frac{m+1}{2}} = \bigl( p^2 \Var \Ga(F,-\L^{-1}F) \bigr)^{\frac{m+1}{2}} \, . 
\end{align*} 
In the present article, we show that, for $M$ of the form $2dm$, the quantity in \eqref{term-to-bound} can be bounded by 
\begin{align*} 
\left( \sum_{1\leq p\leq q\leq d} \Var \Gamma(F_p,F_q) \right)^{\frac{m}{2}} \bigl( \Var \Ga(F,-\L^{-1}F) \bigr)^{\frac{1}{2}} \, . 
\end{align*} 
As mentioned above, however, the quantity 
\[ \sum_{1\leq p\leq q\leq d} \Var \Gamma(F_p,F_q) \] 
is not, in general, equivalent to $\Var \Ga(F,-\L^{-1}F)$, even though a bound expressed solely in terms of the latter would be desirable. This is precisely what we obtain in Corollary \ref{cor:EE_Wiener_sums} by applying the estimate from Proposition \ref{intro:bound-by-var}. However, this estimate is not optimal, as illustrated by several counterexamples presented in Section \ref{sec.shenanigans}.

\medskip

The strategy for bounding \eqref{term-to-bound} is the same as that used in \cite{MPS2025.EE}: we consider the more general term 
\begin{align}\label{term-to-bound-gen} 
\E\left [ \varphi(F) ( \Ga(F,-\L^{-1}F)-1)\right ] \, , 
\end{align} 
where $\varphi:\R\to \R$ is a smooth function with polynomial growth, together with its derivatives. More precisely, since $\varphi(F) = [\varphi]^\oplus(F_1,\ldots,F_d)$, we are naturally led to consider the more general term 
\begin{align*} 
\E\left [ \phi(F_1, \ldots, F_d) ( \Ga(F,-\L^{-1}F)-1)\right ] \, , \end{align*} 
where $\phi:\R^d \to \R$ is a smooth function with polynomial growth, together with its derivatives. 

If $\varphi$ satisfies an appropriate Hermite rank condition, we then derive an expansion of \eqref{term-to-bound-gen} as a sum of terms, each of which can be bounded by $\bigl ( \Var \Ga(F, -\L^{-1} F) \bigr )^{\frac{1}{2}}$ multiplied by a power of $\sum_{1\leq p \leq q \leq d} \Var \Gamma(F_p,F_q)$.

As mentioned above, the main difficulty lies in effectively controlling the terms arising in this expansion. This is achieved through a careful analysis of the operators involved, which are expressed as fractions evaluated at the operator $\LL_\sigma$, using Mehler's formula \eqref{eq:mehler_formula_d} for the Ornstein-Uhlenbeck semigroup $(\Pt_t^\sigma)_{t\geq 0}$.

\medskip

Following \cite{MPS2025.EE}, we split the proof of Theorem \ref{thm:EE_Wiener_sums} into two parts: the case $m=1$, in Section \ref{sec-FO}, and the general case, in Section \ref{sec-HO}. In the first part, the idea is to get more acquainted with the definition of the objects arising from the expansions. Then in the second part, we iterate the arguments used in the first part, and adapt them appropriately.

\section{Proof of Theorem \ref{thm:EE_Wiener_sums} in the case $m=1$}\label{sec-FO}

In this section, we prove Theorem \ref{thm:EE_Wiener_sums} in the case $m=1$. 

\subsection{Multivariate diffusion property and integration by parts}\label{ss:m_ipp_sum}
Let $F$ be an element of $\WW_{\leq d}$ such that $\E[F]=0$ and $\E[F^2]=1$. Write $\sigma_p^2 := \E [ F_p^2 ]$. We suppose that $\sigma_p >0$, that is, $F$ is not orthogonal to $\WW_p$. If this were the case, we could write
\[
F=\sum_{q=1, \, q\neq p}^d F_q,
\]
(or even further, omit any index $p$ such that $\sigma_p = 0$) throughout the remainder of the proof, and the argument would remain unchanged, except for the ranges of the summation indices. In particular, this allows us to write a genuine extension of \cite[Theorem 1.2]{MPS2025.EE}, which we recover by removing all indices but one. For shorter expressions, sometimes we will also denote by $\Fvect$ the vector $(F_1,\ldots,F_d)$.

\medskip

Following the ideas developed in \cite[Section 5]{MPS2025.EE}, we extend the results in a multivariate setting and derive an expansion of the term $\E\left [\varphi(F)(\Gamma(F,-\L^{-1}F)-1)\right ]$ into a sum of terms which are bounded by $\left (\sum_{1\leq p\leq q \leq d} \Var \Gamma(F_p,F_q) \right )^{\frac{1}{2}}\bigl ( \Var\Gamma(F,\L^{-1} F) \bigr )^{\frac{1}{2}} $, up to a constant depending on $d$.

\medskip

We start with the following lemma.
\begin{lemma}
    Let $\varphi$ be an element of $\Cpol{2}$. Then
    \begin{align}\label{eq:relation_diffusion_sum}
        \L\varphi(F) = \varphi''(F) \bigl( \Gamma(F)-\E \left [ \Gamma(F) \right ] \bigr ) + \LL_{\sigma} \left [ \varphi \right ]^{\oplus}(F_1, \ldots, F_d) \, .
    \end{align}
   where $\LL_\sigma$ and $[\varphi]^\oplus$ are defined in Section \ref{s:gaussian_analysis_multid}.
\end{lemma}
\begin{proof}
One has
\begin{eqnarray*}
&& \L\varphi(F) = \varphi''(F)\Gamma(F) + \varphi'(F) \L F \\
&& = \varphi''(F) \bigl ( \Gamma(F)-\E[\Gamma(F)] \bigr ) +\varphi''(F)\sum_{p=1}^d p\sigma_p^2 - \varphi'(F)\sum_{p=1}^d p F_p \\
&&= \varphi''(F) \bigl ( \Gamma(F)-\E[\Gamma(F)] \bigr )+\sum_{p=1}^d p \left ( \sigma_p^2\,\varphi''(F)- F_p\,\varphi'(F) \right )\, .
\end{eqnarray*}
Since for every $p \in \llbracket 1, d \rrbracket$, $\varphi'(F) = \partial_p [\varphi]^\oplus(F_1,\ldots, F_d)$, we get the result.
\end{proof}
\noindent In \cite{MPS2025.EE}, the next step is then to use iteratively \eqref{eq:relation_diffusion_sum}. The problem is that we wrote it for a function of one variable $\varphi$ and get in return a function of $d$ variables $\LL_\sigma [\varphi]^\oplus$. Therefore, we write \eqref{eq:relation_diffusion_sum} in a completely multivariate setting. 
\begin{lemma}\label{lemma:relation_diffusion_d}
    Let $\phi$ be an element of $\Cpold{2}{d}$. Then one has
    \begin{align}\label{eq:relation_diffusion_d}
   \left (\LL_{\sigma}\phi \right)(\Fvect) & = -\sum_{p,q =1}^d \, \partial_{pq}^2 \phi(\Fvect) \, \bigl ( \Gamma(F_p,F_q) - p \sigma_p^2\, \delta_{p,\,q} \bigr )+ \L\phi(\Fvect) \, . 
    \end{align}
\end{lemma}
\begin{proof}
Using the property of diffusion of $\L$, one gets
\begin{eqnarray*}
   &&  \L\phi(\Fvect)  = \sum_{p,q =1}^d \, \partial_{pq}^2 \phi\Fvect) \,\Gamma(F_p,F_q) + \sum_{p=1}^d \, \partial_p \phi(\Fvect)  \L F_p  \\
   && =  \sum_{p,q =1}^d \, \partial_{pq}^2 \phi(\Fvect) \, \bigl ( \Gamma(F_p,F_q) - p \sigma_p^2\, \delta_{p,\,q} \bigr )  + \sum_{p=1}^d \, p\bigl(  \sigma_p^2\,\partial_{p}^2 \phi(\Fvect)-F_p \,\partial_p \phi(\Fvect) \bigr ) \\
    && = \sum_{p,q =1}^d \, \partial_{pq}^2 \phi(\Fvect) \, \bigl ( \Gamma(F_p,F_q) - p \sigma_p^2 \, \delta_{p,\,q} \bigr ) + \left(\LL_{\sigma}\phi \right )(\Fvect) \, . 
\end{eqnarray*}
Hence the result.
\end{proof}
Now, we make use of Lemma \ref{lemma:relation_diffusion_d} to write a first integration-by-parts formula.
\begin{lemma}\label{lemma:IPPalpha}
    Let $\alpha$ be a positive integer and let $\phi$ be an element of $\Cpold{2\alpha+2}{d}$. Then, we have
    \begin{align}\label{eq:IPPalpha}
        \E \left [ \LL_{\sigma}^{\alpha} \phi(\Fvect)(\Gamma(F,-\L^{-1}F)-1) \right ] & = -\sum_{p,q=1}^d \sum_{j=0}^{\alpha-1} \,   \E\left [ \partial_{pq}^2\LL^{\alpha-1-j}_{\sigma}\phi(\Fvect) \overline{\Gamma}_{p,q} \L^{j} (\Gamma(F-\L^{-1}F)-1) \right ] \nonumber \\ & \quad + \,\E\left [\phi(\Fvect)\L^\alpha(\Gamma(F,-\L^{-1}F)-1) \right ],
    \end{align}
    where $\overline{\Gamma}_{p,q} := \Gamma(F_p,F_q)-p\sigma_p^2\, \delta_{p,\,q}$.
\end{lemma}
The proof of Lemma \ref{lemma:IPPalpha} goes verbatim as the proof of \cite[Lemma 5.2]{MPS2025.EE}, using \eqref{eq:relation_diffusion_d} and an inductive argument.

Then, we introduce the appropriate polynomial to cancel the remaining terms that are not of "order $2$", and use it to derive the next proposition.
\begin{lemma}\label{prop:Pannulateur}
    Denote by $P$ the polynomial $X(X+1) \cdots (X+2d-2)$. Then 
    \[\forall F \in \WW_{\leq d}, \quad P(\L)(\Gamma(F,-\L^{-1}F)-1)=0 \, .\]
\end{lemma}
The latter follows from \eqref{eq:gamma_stability}, by noticing that $\Gamma(F,F) \in \WW_{\leq 2d-2}$ for $F \in \WW_{\leq d}$. We will write $P = \sum_{\alpha=1}^{2d-1} a_\alpha X^\alpha$.
\begin{proposition}\label{proposition:IPP_polynomial}
Let $\phi$ be an element of $\Cpold{4d}{d}$. Then we have
    \begin{eqnarray}\label{eq:IPP_polynomial}
        \E \left [ P \left (\LL_{\sigma} \right ) (\phi)(\Fvect)(\Gamma(F-\L^{-1}F)-1) \right ]  = -\sum_{p,q=1}^{d} \sum_{j=0}^{2d-2} \, \E\left [\partial_{pq}^2 P_{j} \left (\LL_{\sigma} \right )(\phi)(\Fvect) \overline{\Gamma}_{p,q} \L^{j}(\Gamma(F,-\L^{-1}F)-1) \right ]. 
    \end{eqnarray}
     where $P_{j}(X) := \sum_{\alpha=j+1}^{2d-1} a_{\alpha} X^{\alpha-1-j}$.
\end{proposition}
\begin{proof}
   Multiplying equality \eqref{eq:IPPalpha} by $a_\alpha$ and summing for all $\alpha$ in $\{1,\ldots,2d-1\}$, one gets
    \begin{eqnarray*}
        &&\E \left [ P\left( \LL_{\sigma}\right ) (\phi)(\Fvect)(\Gamma(F-\L^{-1}F)-1) \right ]  \\ 
        && = -\sum_{p,q=1}^d\sum_{j=`0}^{2d-2} \, \E\left [ \partial_{pq}^2 P_j(\LL_{\sigma})(\phi)(\Fvect) \overline{\Gamma}_{p,q} \L^{j}(\Gamma(F,-\L^{-1}F)-1)\right ]  + \E\left [\phi(\Fvect)P\left (\L\right )(\Gamma(F,-\L^{-1}F)-1) \right ],
    \end{eqnarray*}
    by linearity of both the expectation and the differentiation.
Then, by Lemma \ref{prop:Pannulateur}, one has 
\begin{align*}
    \E\left [\phi(\Fvect)P\left (\L\right )(\Gamma(F,-\L^{-1}F)-1) \right ] = 0 \, .
\end{align*}
\end{proof}
\noindent Inverting $P(\LL_\sigma)$, one gets the desired integration-by-parts below.
\begin{proposition}\label{proposition:IPP_final_d}
Let $\phi$ be an element of $\Cpold{4d}{d}$ such that $\rank_\sigma(\phi) \geq 2d-1$. Then we have
    \begin{eqnarray}\label{eq:IPP_final}
        &&\E \left [ \phi(\Fvect)(\Gamma(F-\L^{-1}F)-1) \right ] = \sum_{p,q=1}^d \sum_{j=0}^{2d-2} \, \E\left [\Tj^{\sigma, (p,q)}_{j }(\phi)(\Fvect)\overline{\Gamma}_{p,q} \L^{j}(\Gamma(F,-\L^{-1}F)-1) \right ],
    \end{eqnarray}
     where $\Tj^{\sigma, (p,q)}_{j } (\phi):= - \partial_{pq}^2 P_{j}(\LL_{\sigma})P(\LL_{\sigma})^{-1}(\phi)$.
\end{proposition}

If $\varphi$ is an element of $\Cpol{4d}$ such that $\rank(\varphi) \geq 2d-1$, then, by virtue of Lemma \ref{lemma:egalite-hermite-rank}, applying Proposition \ref{proposition:IPP_final_d} to $\phi = \varphi^\oplus$ leads to
\begin{align*}
     \E \left [ \varphi(F)(\Gamma(F-\L^{-1}F)-1) \right ] = \sum_{p,q=1}^d \sum_{j=0}^{2d-2} \, \E\left [\Tj^{\sigma, (p,q)}_{j }[\varphi]^\oplus(\Fvect)\overline{\Gamma}_{p,q} \L^{j}(\Gamma(F,-\L^{-1}F)-1) \right ].
\end{align*}
Also, in that case, using the commutation relations $\partial_{pq}^2 \LL_\sigma = (\LL_\sigma-p-q)\partial_{pq}^2$ and \eqref{eq:resolvent_commutation} on smooth enough functions, we have 
\begin{align*}
    \Tj^{\sigma, (p,q)}_{j } [\varphi]^\oplus= - \partial_{pq}^2 P_{j}(\LL_{\sigma})P(\LL_{\sigma})^{-1} [\varphi]^\oplus = - P_j(\LL_\sigma-p-q)P(\LL_\sigma-p-q)^{-1} [\varphi'']^\oplus \, ,
\end{align*}
since $\partial_{pq}^2[\varphi]^\oplus = [\varphi'']^\oplus$. 

\subsection{Study of the $\Tj^{\sigma, (p,q)}_{j}$ and extension of \eqref{eq:IPP_final}}\label{ss:studyofTjk}
In this section, we study the analytic properties of the operators $\Tj^{\sigma, (p,q)}_{j}$, $j \in \llbracket 0, 2d-2 \rrbracket$, evaluated at functions of the form $[\varphi]^\oplus$ where $\varphi \in \Cpol{4d}$. We will still denote $\Tj^{\sigma, (p,q)}_{j}[\varphi]^\oplus$ by $\Tj^{\sigma, (p,q)}_{j}(\varphi)$ and consider $\Tj^{\sigma, (p,q)}_j$ as an operator on functions from $\R$ to $\R$. 
The two mains goals are:
\begin{itemize}
\item[$\bullet$] to extend $\Tj^{\sigma, (p,q)}_{j}$ to $\Cpol{0}$ ;
\item[$\bullet$] to obtain pointwise polynomials estimate on $\Tj^{\sigma, (p,q)}_{j}(\varphi)$, when $\varphi \in \Cpol{0}$.
\end{itemize}
This is done by an analysis of the algebraic expressions of the polynomials $P_{j}, P$, and through the use of the properties of the Ornstein-Uhlenbeck operator, its semigroup and its resolvent, available in Sections \ref{ss:semigroup_d} and \ref{ss:resolvent_d}.

\medskip

The factorization in \cite[Section 5.2]{MPS2025.EE} still holds with the necessary adjustments.
\begin{lemma}\label{lemma:rewrite_of_T_d}
    Let $j$ be an element of $\llbracket 0, 2d-2\rrbracket$ and let $p,q$ be elements of $\llbracket 1, d \rrbracket$. For all $\varphi \in \Cpol{4d}$, 
    \begin{align}\label{eq:Tj_d}
        \Tj_j^{\sigma, (p,q)}(\varphi) & = -\Q_j^{(p,q)}(\LL_\sigma) \S_\sigma^{(p,q)} (\varphi),
    \end{align}
    where $\S_\sigma^{(p,q)}(\varphi) := (\LL_{\sigma}-p-q)^{-1}[\varphi'']^{\oplus}$ and 
    \begin{align*}
        \Q_j^{(p,q)}(X) := (X-p-q)^{-j} - \underset{\alpha=1}{\overset{j}{\sum}} \, a_\alpha (X-p-q)^{\alpha-j}P ( X-p-q  )^{-1} \, .
    \end{align*}
\end{lemma}
\begin{proof}
    Observe that $X^{j+1}P_j(X) = P(X) - \sum_{\alpha=1}^{j} \, a_\alpha X^\alpha$, hence
\begin{align}\label{az}
    P_j(X-p-q)P(X-p-q)^{-1} &= (X-p-q)^{-j-1} - \underset{\alpha=1}{\overset{j}{\sum}} \, a_\alpha  (X-p-q)^{\alpha-j-1}P  ( X-p-q )^{-1} \, .
\end{align}
Define
\begin{align*}
   \Q_j^{(p,q)}(X) & := (X-p-q)^{-j} - \underset{\alpha=1}{\overset{j}{\sum}} \, a_\alpha  (X-p-q)^{\alpha-j}P  ( X-p-q )^{-1} \, ,
\end{align*}
and, for all $\varphi \in \Cpol{4d}$, $\S_\sigma^{(p,q)}(\varphi) := (\LL_\sigma-p-q)^{-1} [\varphi'']^\oplus$. Then, one has
\begin{align*}
       P_j(X-p-q)P(X-p-q)^{-1} & = \Q_j^{(p,q)}(X)(X-p-q)^{-1},
\end{align*}
and consequently, the results follows.
\end{proof}

\begin{remark}
    On $\Cpold{4d}{d}$, one can write $\Tj_{j}^{(p,q)} (\phi) = \Q_j^{(p,q)}(\LL_\sigma) (\LL_\sigma -p-q)^{-1} \partial_{pq}^2\phi$.
\end{remark}

Since $\Q_{j}^{(p,q)}(\LL_\sigma)$ is a sum of products of resolvent operators, it is well defined on functions in $\Cpold{0}{d}$, with Hermite rank greater than $2d-1$, and satisfies polynomial estimates; see Lemma \ref{lemma:estimates_resolvent_polynomial_d}. Therefore, to establish the same properties for $\Tj_{j}^{\sigma,(p,q)}$, it remains only to analyze the operator $\S_\sigma^{(p,q)}$. This is the purpose of the following key result.

\begin{proposition}\label{proposition:prolongement_S_d_1}
Let $p,q$ be elements of $\llbracket 1,d\rrbracket$. The operator $\S_\sigma^{(p,q)}$ can be extended to an operator on functions $\varphi \in \Cpol{0}$. Furthermore, if $\varphi \in \Cpol{0}$ verifies 
\begin{align*}
\forall x \in \R^d,\quad \left | [\varphi]^\oplus(x)\right | \leq C \left ( 1 +\| I_\sigma^{-1} x \|^l \right ),
\end{align*}
for some $C>0$ and $l \in \N$, then there exist $A_{l}^{(p,q)}>0$ such that
\begin{align}\label{ineq:estimate_S_d_1}
\forall x \in \R^d, \quad \left | \S_\sigma^{(p,q)}(\varphi)(x) \right | \leq CA_{l}^{(p,q)} \left ( 1 + \|I_\sigma^{-1} x \|^{l+1} \right ).
\end{align}
\end{proposition}
\begin{proof}
Let $\varphi$ be an element of $\Cpol{2}$. Using the semigroup formula for the resolvent operator in Proposition \ref{proposition:resolvent_d}, one has 
\begin{align*}
\S_\sigma^{(p,q)}(\varphi)(x) & = -\int_0^{+\infty} e^{-(p+q)t} \,\Pt^\sigma_t[\varphi'']^\oplus( x) \, \dint t \, .
\end{align*}
 Since $[\varphi']^\oplus = \partial_b[\varphi]^\oplus$ for all $b \in \llbracket 1, d\rrbracket$, one has
    \begin{align}\label{eq:astuce_prolongement}
        [\varphi'']^{\oplus}(x_1,\ldots,x_d) &= \frac{\sum_{b=1}^d b\sigma_b^2 \partial_b^2 [\varphi]^\oplus(x_1,\ldots,x_d)}{\sum_{b=1}^d b\sigma_b^2} \nonumber\\
        & = \frac{\LL_\sigma [\varphi]^\oplus(x_1,\ldots,x_d)}{\E \left [ \Gamma(F) \right ]}+\frac{\sum_{b=1}^d bx_b\partial_b  [\varphi]^\oplus(x_1,\ldots,x_d)}{\E \left [ \Gamma(F) \right ]},
    \end{align}
    for every $x \in \R^d$. Using that $\partial_t \Pt^\sigma_{t}[\varphi]^\oplus = \Pt^\sigma_{t}\LL_\sigma[\varphi]^\oplus$, and writing 
    \begin{eqnarray*}
        && \Pt^\sigma_{t} \left ( x_b\partial_b[\varphi]^\oplus \right )( x) = \E \left [ \left ( x_be^{-bt}+\sigma_b\left (1-e^{-2bt}\right )^{1/2}N_b \right )\partial_b[\varphi]^\oplus \left ( U_t x+M_t N_\sigma \right ) \right ] \\
        && = \frac{x_be^{-bt}}{\sigma_b\left ( 1-e^{-2bt} \right )^{1/2}}\E \left [ N_b[\varphi]^\oplus \left ( U_tx+M_t N_\sigma \right ) \right ] +\E \left [ H_2(N_b) [\varphi]^\oplus \left ( U_tx+M_t N_\sigma \right ) \right ],
    \end{eqnarray*}
   we obtain
    \begin{align*}
    \Pt^\sigma_t [\varphi'']^\oplus( x) & = \frac{\partial_t \Pt^\sigma_t [\varphi]^\oplus( x)}{\E \left [ \Gamma(F) \right ]}
         +\frac{1}{\E \left [ \Gamma(F) \right ]}\E \left [ \left ( \sum_{b=1}^d bH_2(N_b) \right ) [\varphi]^\oplus\left (U_tx+M_t N_\sigma \right ) \right ] \\
        & \quad+\sum_{b=1}^d \frac{bx_be^{-bt}}{\sigma_b\left ( 1-e^{-2bt} \right )^{1/2}\E \left [ \Gamma(F) \right ]}\E \left [ N_b[\varphi]^\oplus \left ( U_tx+M_t N_\sigma \right ) \right ] \\
        & = \frac{1}{\E[\Gamma(F)]} \left ( \partial_t \Pt^\sigma_t [\varphi]^\oplus( x)+ \Pt_t^{\sigma,V_{2}} [\varphi]^\oplus(x) + \sum_{b=1}^d \frac{x_b}{\sigma_b} \frac{e^{-bt}}{\left (1-e^{-2bt} \right )^{1/2}} \Pt_t^{\sigma, V_{1,b}} [\varphi]^\oplus(x) \right),
    \end{align*}
where $V_{1,b}(x_1,\ldots,x_d) \coloneq bH_{e_b}(x)$ and $V_2(x_1,\ldots,x_d) \coloneq \sum_{b=1}^d bH_{2e_b}(x)$. Using inequality \eqref{eq:controle_semigroup_hermite_final_d} and Lemma \ref{lemma:estimates_PtV_pol}, one gets that 
\begin{align}\label{aux:integrabilité}
    \frac{1}{\E[\Gamma(F)]} \left ( \partial_t \Pt^\sigma_t [\varphi]^\oplus( x)+ \Pt_t^{\sigma,V_{2}} [\varphi]^\oplus(x) + \sum_{b=1}^d \frac{x_b}{\sigma_b} \frac{e^{-bt}}{\left (1-e^{-2bt} \right )^{1/2}} \Pt_t^{\sigma, V_{1,b}} [\varphi]^\oplus(x) \right)
\end{align}
is integrable with respect to $t$ on $[0,+\infty[$, once multiplied by $e^{-pt-qt}$. After integrating by parts, we obtain
\begin{eqnarray*}
&& \S_\sigma^{(p,q)}(\varphi)(x) \\
&& = \frac{1}{\E[\Gamma(F)]} \left ( [\varphi]^\oplus(x) + (p+q)\RR_\sigma(p+q)[\varphi]^\oplus(x) + \RR_{\sigma}^{V_2}(p+q)[\varphi]^\oplus(x) + \sum_{b=1}^d \frac{x_b}{\sigma_b} \RR_{\sigma,b}^{V_{1,b}}(p+q)[\varphi]^\oplus(x) \right ),
\end{eqnarray*}
for every $x\in\R^d$, where, for $V$ polynomial, 
\begin{align*}
    \RR_\sigma^{V}(p+q)[\varphi]^\oplus(x) := -\int_0^{+\infty} e^{-(p+q)t} \, \Pt_{t}^{\sigma, V}[\varphi]^\oplus(x) \, \dint t \, ,
\end{align*}
and for $b \in \llbracket 1,d\rrbracket$,
\begin{align}
\RR_{\sigma,b}^{V}(p+q)[\varphi]^\oplus(x) \coloneq -\int_0^{+\infty} \frac{e^{-(p+q+b)t}}{\left ( 1-e^{-2bt} \right )^{1/2}} \Pt_t^{\sigma,V}[\varphi]^\oplus(x) \, \dint t \, .
\end{align}
Using again the estimate \eqref{eq:controle_semigroup_hermite_final_d}, we deduce that there exist a constant $A_l^{(p,q)}>0$ such that \eqref{ineq:estimate_S_d_1} is verified.
\end{proof}
\begin{remark}
Let us note that equality \eqref{eq:astuce_prolongement}, that describes $\left [\varphi''\right ]^\oplus = \partial_{pq}^2 [\varphi]^\oplus$ in terms of $\LL_\sigma [\varphi]^\oplus$ and $\left [ \varphi'\right ]^\oplus= \partial_b[\varphi]^\oplus$ (thus allowing us to extend $\S_\sigma^{(p,q)}$) is specific to the case where $F$ is in a finite sum of chaoses. Providing a relation between the terms
\begin{align*}
\LL_\sigma \phi, \qquad \partial_b \phi, \quad b \in \{1,\ldots,d\},
\end{align*} 
and $\partial_{pq}^2 \phi$ is more difficult -- maybe even impossible. The absence of cross derivatives in the expression of $\LL_\sigma$ on smooth functions has important consequences. In particular, extending this approach to multivariate Edgeworth expansions is likely to yield bounds in stronger, more regular distances, as it would presumably require additional regularity assumptions on $\phi_h'$, and hence on $h$.
\end{remark}

\begin{proposition}\label{proposition:continuity_S_d_1}	
Let $(\varphi_n)_{n\in\N}$ be a sequence of elements of $\Cpol{0}$ such that $\forall n \in \N, \;\rank(\varphi_n)>0$ and 
\begin{align*}
\forall n \in \N, \; \forall x \in \R^d, \quad |\varphi_n(x)| \leq C \left (1+\|I_\sigma^{-1} x \|^l \right ),
\end{align*}
for some $C>0$ and $l\in \N$. Suppose that the sequence $(\varphi_n)_{n\in\N}$ converges pointwise to a continuous function $\varphi$. Then the sequence $\left ( \S_\sigma^{(p,q)}(\varphi_n) \right )_{n\in\N}$ converges pointwise to $\S_\sigma^{(p,q)}(\varphi)$, for every $p,q$ in $\llbracket 1,d \rrbracket$.
\end{proposition}
\begin{proof}
We have seen in the proof of Proposition \ref{proposition:prolongement_S_d_1} that, for every $n$ in $\N$,
\begin{eqnarray*}
&& \S_\sigma^{(p,q)}(\varphi_n)(x) = \\
&&  \frac{1}{\E[\Gamma(F)]} \left ( [\varphi_n]^\oplus(x) + (p+q)\RR_\sigma(p+q)[\varphi_n]^\oplus(x) + \RR_{\sigma}^{V_2}(p+q)[\varphi_n]^\oplus(x) + \sum_{b=1}^d \frac{x_b}{\sigma_b} \RR_{\sigma,b}^{V_{1,b}}(p+q)[\varphi_n]^\oplus(x) \right ),
\end{eqnarray*}
Using Lemma \ref{lemma:continuity_PtV}, and the appropriate dominations given by \eqref{eq:controle_semigroup_hermite_final_d} and Lemma \ref{lemma:estimates_PtV_pol}, we deduce the result with Lebesgue's domination convergence theorem.
\end{proof}

Combining the factorization established in Lemma \ref{lemma:rewrite_of_T_d} with the properties proved in Propositions \ref{proposition:prolongement_S_d_1} and \ref{proposition:continuity_S_d_1} yields the desired results.

\begin{proposition}\label{proposition:continuity_Tjk_d}
Let $j$ be an element of $\llbracket 0, 2d-2 \rrbracket$ and let $p,q$ be elements of $\llbracket 1, d \rrbracket$. Let $(\varphi_n)_{n\in\N}$ be a sequence of elements of $\Cpol{0}$ such that $\forall n \in \N, \;\rank(\varphi_n)\geq 2d-1$ and 
\begin{align*}
\forall n \in \N, \; \forall x \in \R^d, \quad \left |[\varphi_n(x)]^\oplus \right | \leq C \left (1+\|I_\sigma^{-1} x \|^l \right ),
\end{align*}
for some $C>0$ and $l\in \N$. Suppose that the sequence $(\varphi_n)_{n\in\N}$ converges pointwise to a continuous function $\varphi$. Then the sequence $\left ( \Tj_{j}^{\sigma,(p,q)}(\varphi_n) \right )_{n\in\N}$ converges pointwise to $\Tj^{\sigma, (p,q)}_{j}(\varphi)$.
\end{proposition}
\begin{proof}
By Lemma \ref{lemma:rewrite_of_T_d}, the operator $\Tj^{\sigma, (p,q)}_{j}$ is a composition between a sum of products of resolvent operators of $\LL_\sigma$ and the operator $\S_\sigma^{(p,q)}$. Thus, combining \eqref{ineq:estimate_S_d_1} and Proposition \ref{proposition:continuity_S_d_1}, and then several applications of Lemmas \ref{lemma:estimates_resolvent_polynomial_d} and \ref{lemma:continuity_resolvent_d} lead to the result.
\end{proof}

Before stating the next proposition, we need to understand how the operator $\S_\sigma^{(p,q)}$ acts on the Hermite rank of continuous functions. To this end, we use the following classical lemma, whose proof can be found, for instance, in \cite[Section~5.3]{MPS2025.EE}.

\begin{lemma}\label{lemma:approximation_function_Hrank}
    Let $\varphi$ be an element of $\Cpol{0}$ with Hermite rank $r>0$. Suppose
    \begin{align*}
    \forall x \in \R, \quad |\varphi(x)| \leq C(1+|x|^d),
    \end{align*}
    where $C>0$ and $d$ is a non-negative integer. Then there exists a sequence $(\varphi_n)$ of functions in $\Cpol{\infty}$ with Hermite rank $\ge r$ such that:
    \begin{enumerate}[(i)]
        \item $(\varphi_n)$ converges pointwise to $\varphi$ ;
        \item there exists $M>0$ such that $\forall n \in \N, \forall x \in \R, \quad |\varphi_n(x)| \leq M(1+|x|^d)$.
    \end{enumerate}
\end{lemma}

\noindent Then, we derive the next result.

\begin{lemma}
    Let $(p,q)$ be an element of $\llbracket 1,d \rrbracket^2$ and let $\varphi$ be an element of $\Cpol{0}$. Then
    \begin{align*}
        \rank_\sigma \bigl (\S_\sigma^{(p,q)}(\varphi) \bigr ) \geq \max(0,\rank(\varphi)-2)\, .
    \end{align*}
\end{lemma}
\begin{proof}
    Suppose that $\rank(\varphi) \geq 2$, otherwise the inequality is clear. To prove the result, we approximate $\varphi$ by a sequence $(\varphi_n)_{n \geq 1}$ of functions in $\Cpol{\infty}$ given by Lemma \ref{lemma:approximation_function_Hrank}. First, let us notice that the classical Gaussian integration by parts, yields for $\psi \in \Cpol{2}$,
    \begin{align*}
       \forall n \geq 2, \quad  c_n(\psi) = \frac{1}{n(n-1)} c_{n-2}(\psi'') \, ,
    \end{align*}
    hence $\rank(\psi'') = \rank(\psi)-2$. Then, by Lemmas \ref{lemma:egalite-hermite-rank} and \ref{lemma:conservation_rank_d}, we deduce that
    \begin{align*}
        \forall n \geq 2, \quad \rank_\sigma \left ( \S_\sigma^{(p,q)}(\varphi_n) \right ) = \rank_\sigma \left (  (\LL_\sigma -p-q)^{-1} [\varphi_n'']^\oplus \right ) \geq \rank(\varphi_n)-2 \, .
    \end{align*}
    Therefore, by Propositions \ref{proposition:prolongement_S_d_1}, \ref{proposition:continuity_S_d_1}, and a straightforward application of Lebesgue's dominated convergence theorem, letting $n$ go to infinity yields
    \begin{align*}
        \rank_\sigma \left ( \S_\sigma^{(p,q)}(\varphi) \right ) \geq \rank(\varphi)-2 \, ,
    \end{align*}
    as desired.
\end{proof}

\begin{proposition}\label{proposition:estimates_Tjk_d}
Let $j$ be an element of $\llbracket 0, 2d-2 \rrbracket$ and let $p,q$ be elements of $\llbracket 1, d \rrbracket$. Let $\varphi$ be an element of $\Cpol{0}$ such that $\rank(\varphi) \geq 2d-1$, and
\begin{align*}
\forall x \in \R^d, \quad  \left | [\varphi]^\oplus(x)\right | \leq C \left ( 1 +\| I_\sigma^{-1} x \|^l \right ),
\end{align*}
for some $C>0$, $l \in \N$. Then there exist $A_{l,j}^{(p,q)}>0$ such that
\begin{align}\label{ineq:estimate_Tjk_d}
\forall x \in \R^d, \quad \left | \Tj^{\sigma,(p,q)}_{j}(\varphi)(x) \right | \leq CA_{l,j}^{(p,q)} \left ( 1 + \|I_\sigma^{-1} x \|^{l+\mu_{1,(p,q)}} \right ) \, ,
\end{align}
where $\mu_{1,(p,q)} :=1 + \sum_{k=1}^{d_1-p-q+1} k$.
\end{proposition}
\begin{proof}
Again, we leverage the factorization \eqref{eq:Tj_d}. The operator $\Q_{j}^{(p,q)}(\LL_\sigma)$ has the form
\begin{align}\label{aux:form-resolvent-Q}
    \RR_\sigma(p+q)^{j}  - \underset{\alpha=1}{\overset{j}{\sum}} \, a_\alpha \RR_\sigma(p+q)^{j-\alpha} \prod_{k=0}^{d_1} \RR_\sigma(p+q-k) \,.
\end{align}
First, by Proposition \ref{proposition:prolongement_S_d_1}, there exists $A_l^{(p,q)}>0$ such that 
\begin{align*}
    \forall x \in \R^d, \quad |\psi(x)| \leq C A_{l}^{(p,q)} \left (1 + \|I_\sigma^{-1} x \|^{l+1} \right ) \, ,
\end{align*}
where $\psi := \S_\sigma^{(p,q)}(\varphi)$. Then we apply Lemma \ref{lemma:estimates_resolvent_polynomial_d} to each term in the expression \eqref{aux:form-resolvent-Q}, applied to $\psi$. Since $p+q>0$, there exists $K_{d,l,j}^{(p,q)} >0$ such that
\begin{align}\label{aux1:estimates-T}
    \forall x \in \R^d, \quad |\RR_\sigma(p+q)^{j}\psi(x)| \leq C K_{d,l,j}^{(p,q)} \left (1 + \|I_\sigma^{-1} x \|^{l+1} \right ) \, .
\end{align}
Now, for $j\geq 1$, in the product $\prod_{k=0}^{d_1} \RR_\sigma(p+q-k)$, some of the $p+q-k$ can be negative. In that case, since $\rank_\sigma(\psi) \geq \rank(\varphi)-2 \geq d_1-p-q+1$ and by Corollary \ref{corollary:prod_resolvent_conservation_rank_d}, we deduce that there exists $\tilde{K}_{d,l,j}^{(p,q)}>0$ such that 
\begin{align}\label{aux2:estimates-T}
    \left |\left (\RR_\sigma(p+q)^{j-\alpha} \prod_{k=0}^{d_1} \RR_\sigma(p+q-k)\right )\psi(x) \right | \leq C\tilde{K}_{d,l,j} \left ( 1 + \|I_\sigma^{-1} x \|^{l+\mu_{1,(p,q)}} \right ) \, ,
\end{align}
where
\begin{align*}
    \mu_{1,(p,q)} & := 1 + \sum_{k=p+q}^{d_1} (k-p-q +1) = 1 + \sum_{k=0}^{d_1-p-q} (k+1) = 1 + \sum_{k=1}^{d_1-p-q+1} k \\
    & = \begin{cases}
        1 & \text{ if } d_1 < p+q \, ; \\
        1+\frac{(d_1-p-q+1)(d_1-p-q+2)}{2} & \text{ if } d_1 \geq p+q\,.
    \end{cases}
\end{align*}
Combining \eqref{aux1:estimates-T} and \eqref{aux2:estimates-T} together yields the result.
\end{proof}
\begin{remark}
    The exponent $\mu_{1,(p,q)}$ can be reduced by using the integral representation of a product of resolvent operators given in Proposition \ref{proposition:resolvent_d}, together with the bound \eqref{eq:controle_semigroup_hermite_final_d} to control the semigroup appearing in this representation.
\end{remark}
\noindent Finally, we extend \eqref{eq:IPP_final} to functions $\varphi \in \Cpol{0}$, with Hermite rank high enough, by an approximation argument. 

\noindent Then, the following extension follows.

\begin{proposition}\label{proposition:IPP_final_extended}
Let $\varphi$ be an element of $\Cpol{0}$ such that $\rank(\varphi)\geq 2d-1$. Then 
 \begin{align}\label{eq:IPPfinal}
         \E \left [ \varphi(F)(\Gamma(F,-\L^{-1}F)-1) \right ] = \sum_{p,q=1}^d \sum_{j=0}^{2d-2} \, \E\left [\Tj^{\sigma, (p,q)}_{j }(\varphi)(\Fvect)\overline{\Gamma}_{p,q} \L^{j}(\Gamma(F,-\L^{-1}F)-1) \right ].
 \end{align}
\end{proposition}

\begin{proof}
We approximate $\varphi$ by the sequence given by Lemma \ref{lemma:approximation_function_Hrank} and apply \eqref{eq:IPP_final} to it. Then, using Propositions \ref{proposition:continuity_Tjk_d} and \ref{proposition:estimates_Tjk_d}, we apply Lebesgue's dominated convergence theorem and we take the limit when $n$ goes to $+\infty$ to conclude. See the proof of \cite[Proposition 5.16]{MPS2025.EE} for more details.
\end{proof}

\subsection{Bounding the remainder}\label{sec.bounding}

In this section, we provide a bound on the total variation distance between $\prob_F$ and $\boldsymbol{\gamma}_{F,\,2d}$ using the identity 
\begin{align}\label{identity.to.bound}
    \dTV(\prob_F,\boldsymbol{\gamma}_{F, \, 2d}) 
    & = \frac{1}{2} \underset{\underset{\|h\|_\infty \leq 1}{h \in 
    \mathscr{C}^0(\R)}}{\sup} \,\left | \E \left [ \projgeq{2d-1}(\phi_h')(F) 
    ( \Gamma(F, -\L^{-1} F) -1)\right ]  \right | \, .
\end{align} Let us fix $h \in \mathscr{C}^0_b(\R)$ such that $\|h\|_\infty\leq 1$. Then, the Stein solution $\phi_h$ defined in \eqref{eq:stein_solll} is continuously differenatiable and verifies $\|\phi'_h\| \leq 4$. In particular, by Lemma \ref{lemma:estimates_proj_polynomial}, the function $\projgeq{2d-1}(\phi_h')$ is in $\Cpol{0}$. Thus, applying \eqref{eq:IPPfinal} to $\projgeq{2d-1}(\phi_h')$, one gets
\begin{eqnarray}\label{expansion.proj}
    && \E \left [ \projgeq{2d-1}(\phi_h')(F)(\Gamma(F,-\L^{-1}F)-1) \right ] \nonumber \\
    && = \sum_{p,q=1}^d \sum_{j=0}^{2d-2} \, \E\left [\Tj^{\sigma, (p,q)}_{j } \bigl(\projgeq{2d-1}(\phi_h')\bigr )(\Fvect)\overline{\Gamma}_{p,q} \L^{j}(\Gamma(F,-\L^{-1}F)-1) \right ] \, .
\end{eqnarray}
We proceed to bound each term
\begin{align}\label{term.to.bound}
\E\left [\Tj^{\sigma, (p,q)}_{j } \bigl(\projgeq{2d-1}(\phi_h')\bigr )(\Fvect)\overline{\Gamma}_{p,q} \L^{j}(\Gamma(F,-\L^{-1}F)-1) \right ]\, ,
\end{align}
in the expansion \eqref{expansion.proj} above, independently of $h$. We follow the same types of arguments that in \cite[Section 5.4]{MPS2025.EE}. Let us start with the following lemma (for a proof see, e.g., \cite[Lemma 5.19]{MPS2025.EE}).
\begin{lemma}\label{bounds-L}
    If $q\in\N^*$, $W\in \WW_{\le q}$, and $k\in\N$, then 
    \begin{align*}
        \E \bigl[|\L^k W|^2\bigr]^{\frac{1}{2}} \;\le\; q^{k}\,\E \bigl[|W|^2\bigr]^{\frac{1}{2}} \, .
    \end{align*}
\end{lemma}
\noindent Then, we prove this second lemma.
\begin{lemma}\label{lemma:pre-boundation}
Let $j$ be an element of $\llbracket 0, 2d-2 \rrbracket$ and let $p,q$ be elements of $\llbracket 1, d \rrbracket$. Then there exist $A_{j}^{(p,q)}>0$ such that 
\begin{align*}
\E\left [\Tj^{\sigma, (p,q)}_{j }\bigl(\projgeq{2d-1}(\phi_h')\bigr )(\Fvect)\overline{\Gamma}_{p,q} \L^{j}(\Gamma(F,-\L^{-1}F)-1) \right ] \leq A_{j}^{(p,q)} \bigl ( \Var\Gamma(F_p,F_q) \bigr )^{\frac{1}{2}} \bigl ( \Var\Gamma(F,\L^{-1} F ) \bigr )^{\frac{1}{2}}.
\end{align*}
\end{lemma}
\begin{proof}
Since $\phi_h'$ is continuous and bounded, by Lemma \ref{lemma:estimates_proj_polynomial}, there exists $B_{d} >0$ such that
\begin{align*}
    \forall x \in \R, \quad \left | \projgeq{2d-1}(\phi_h')(x) \right | \leq B_{d}\left( 1 + |x|^{2d-2} \right ) \, .
\end{align*}
Thus, by Proposition \ref{proposition:estimates_Tjk_d}, there exists $B^{(p,q)}_{d,j} >0$ such that 
\begin{eqnarray}\label{aux:bounding-3}
    && \left | \E\left [\Tj^{\sigma, (p,q)}_{j }\bigl(\projgeq{2d-1}(\phi_h')\bigr )(\Fvect)\overline{\Gamma}_{p,q} \L^{j}(\Gamma(F,-\L^{-1}F)-1) \right ] \right | \nonumber \\
    && \leq B^{(p,q)}_{d,j} \, \E \left [ \left (1+\| I_\sigma^{-1} \Fvect \|^{2d-2+\mu_1^{(p,q)}}\right ) |\overline{\Gamma}_{p,q} \L^{j}(\Gamma(F,-\L^{-1}F)-1) |\right ] \nonumber\\
    && \leq B^{(p,q)}_{d,j} \, \E\left [ |\overline{\Gamma}_{p,q}|^{2}\right ]^{\frac{1}{2}}\, \E \left [ \left | \L^{j}(\Gamma(F,-\L^{-1}F)-1) \right |^{2}\right ]^{\frac{1}{2}} \nonumber\\
    && \qquad + B^{(p,q)}_{d,j} \,  \E \left [\| I_\sigma^{-1} \Fvect \|^{8d-8+4\mu_1^{(p,q)}} \right ]^{\frac{1}{4}} \, \E\left [ |\overline{\Gamma}_{p,q}|^{4}\right ]^{\frac{1}{4}} \E \left [ \left | \L^{j}(\Gamma(F,-\L^{-1}F)-1) \right |^{2}\right ]^{\frac{1}{2}} \, ,
\end{eqnarray}
where we used the Cauchy-Schwarz inequality several times. Now, notice that, by convexity,
\begin{align*}
    \| I_\sigma^{-1} \Fvect \|^{8d-8+4\mu_1^{(p,q)}} & = \left ( \sum_{\ell = 1}^d \sigma_\ell^{-2} F_\ell^2 \right )^{4d-4+2\mu_1^{(p,q)}} \leq d^{4d-5+2\mu_1^{(p,q)}} \sum_{\ell = 1}^d (\sigma_\ell^{-2} F_\ell^2)^{4d-4+2\mu_1^{(p,q)}} \, .
\end{align*}
Hence, by sub-additivity of the power function with exponent less than $1$,
\begin{align*}
    \E \left [\| I_\sigma^{-1} \Fvect \|^{8d-8+4\mu_1^{(p,q)}} \right ]^{\frac{1}{4}} &\leq d^{\frac{4d-5+2\mu_1^{(p,q)}}{4}} \sum_{\ell = 1}^d \E \left [(\sigma_\ell^{-2} F_\ell^2)^{4d-4+2\mu_1^{(p,q)}} \right ]^\frac{1}{4}\\
    & \leq d^{4d-5+2\mu_1^{(p,q)}} d \underset{1\leq \ell \leq d}{\max} \E\left [ (\sigma_\ell^{-2} F_\ell^2)^{4d-4+2\mu_1^{(p,q)}} \right ]^\frac{1}{4} \, .
\end{align*}
Then, since $\Gamma(F,-\L^{-1} F) \in \WW_{\leq 2d-2}$, by Lemma \ref{bounds-L} we have 
\begin{align*}
    \E \left [ \left | \L^{j}(\Gamma(F,-\L^{-1}F)-1) \right |^{2}\right ]^{\frac{1}{2}} \leq (2d-2)^{j} \, \E\left [ |\Gamma(F,-\L^{-1}F)-1|^2 \right ]^\frac{1}{2} = (2d-2)^{j} \bigl ( \Var\Gamma(F,-\L^{-1}F) \bigr )^{\frac{1}{2}} \, .
\end{align*}
Moreover, using the hypercontractivity estimates \eqref{eq:hypercontractivity} (applied to the relevant chaoses), we have 
\begin{align}\label{aux:bounding-1}
    \E\left [ |\overline{\Gamma}_{p,q}|^4 \right ]^{\frac{1}{4}} \leq 3^{p+q-2}\, \E \left [ |\overline{\Gamma}_{p,q}|^2  \right ]^{\frac{1}{2}} \, ,\quad \E\left [ (\sigma_\ell^{-2} F_\ell^2)^{4d-4+2\mu_1^{(p,q)}} \right ] \leq (4d-5+2\mu_1^{(p,q)})^{\frac{\ell}{2}(4d-4+2\mu_1^{(p,q)})} \, ,
\end{align}
since $\E[\sigma_\ell^{-2} F_\ell^2 ] = 1$. In particular, we get
\begin{align}\label{aux:bounding-2}
    \E \left [\| I_\sigma^{-1} \Fvect \|^{8d-8+4\mu_1^{(p,q)}} \right ]^{\frac{1}{4}} & \leq d^{4d-5+2\mu_1^{(p,q)}} (4d-5+2\mu_1^{(p,q)})^{\frac{d}{8}(4d-4+2\mu_1^{(p,q)})} d \, .
\end{align}
Combining \eqref{aux:bounding-3} with \eqref{aux:bounding-1} and \eqref{aux:bounding-2}, we deduce the result wanted.
\end{proof}
\noindent Then, summing over all $p,q$ leads to the following bound.
\begin{proposition}\label{proposition:EE_stein_discrepancy_d_1}
There exists $C_{d} >0$ such that, for all $F \in \WW_{\leq d}$ with $\E[F]=0$ and $\E[F^2]=1$, and all $h\in \mathscr{C}^0_b(\R)$ with $\|h\|_\infty \le 1$,
\begin{align*}
    \bigl|\E \bigl[ \operatorname{Proj}_{\geq 2d-1}(\phi'_h)(F)\,(\Ga(F,-\L^{-1} F)-1) \bigr] \bigr|
    \le C_d\, \left ( \sum_{1\leq p \leq q \leq d} \Var\Gamma(F_p,F_q) \right )^{\frac{1}{2}} \bigl(\Var \Ga(F,\L^{-1} F)\bigr)^{\frac{1}{2}}.
\end{align*}
\end{proposition}

\noindent The preceding estimate can also be expressed solely in terms of $\sum_{p\leq q} \Var\Gamma(F_p,F_q)$, although we do not use this formulation.

\begin{lemma}\label{lemma:var_gamma_Lmoins1}
Let $F$ be an element of $\WW_{\leq d}$ such that $\E[F^2]=1$ Then there exists $C_d>0$ such that 
\begin{align}\label{ineq:borne.gamma.moins1}
\Var \Gamma(F, \L^{-1} F )  \leq C_d  \sum_{1\leq p \leq q \leq d} \Var\Gamma(F_p,F_q)\,.
\end{align}
\end{lemma}
\begin{proof}
By linearity, one has
\begin{align*}
    \Gamma(F,-\L^{-1} F) = \sum_{p=1}^d \frac{1}{p} \, \Gamma(F_p) + \sum_{p<q} \left (\frac{1}{p}+\frac{1}{q} \right ) \Gamma(F_p,F_q) \, .
\end{align*}
Hence, by convexity,
\begin{align*}
    \Var \Gamma(F, \L^{-1} F ) & \leq \frac{d(d+1)}{2} \left ( \sum_{p=1}^d \frac{1}{p^2} \Var \Gamma(F_p) + \sum_{p<q} \left ( \frac{1}{p}+\frac{1}{q} \right )^2 \Var\Gamma(F_p,F_q) \right ) \\
    & \leq 2d(d+1) \sum_{1\leq p \leq q \leq d} \Var\Gamma(F_p,F_q) \, .
\end{align*}
\end{proof}

Finally, using \eqref{identity.to.bound} and Proposition \ref{proposition:EE_stein_discrepancy_d_1}, one gets, for all $F\in \WW_{\leq d}$ such that $\E[F]=0$ and $\E[F^2]=1$,
\begin{align}
    \dTV(\prob_F,\boldsymbol{\gamma}_{F,2d}) & \leq C_d \left ( \sum_{1\leq p \leq q \leq d}  \Var \Gamma(F_p, F_q) \right )^{\frac{1}{2}}\bigl ( \Var\Gamma(F,\L^{-1} F) \bigr )^{\frac{1}{2}} \, .
\end{align}
That concludes the proof of Theorem \ref{thm:EE_Wiener_sums} in the case $m=1$.

\section{Higher order expansion}\label{sec-HO}

In this section, we prove \ref{thm:EE_Wiener_sums} for any $m\geq1$. Let $m$ be a positive integer and let $F$ be an element of $\WW_{\leq d}$ such that $\E[F]=0$ and $\E[F^2]=1$. 

As done in \cite{MPS2025.EE}, we are able to generalize the expansions to get higher order terms, the main difficulty being to extend the expansion for continuous functions.

\medskip

As previously done in the case of a single Wiener chaos $F \in \WW_d$, the idea here is to first write out the expansion 
 \begin{eqnarray*}
        &&\E \left [ \phi(\Fvect)(\Gamma(F-\L^{-1}F)-1) \right ] = \sum_{p,q=1}^d \sum_{j=0}^{2d-2} \, \E\left [\Tj^{\sigma, (p,q)}_{j }(\phi)(\Fvect)\overline{\Gamma}_{p,q} \L^{j}(\Gamma(F,-\L^{-1}F)-1) \right ],
    \end{eqnarray*}
for $\phi$ smooth enough. Then the idea is to write the same expansion \eqref{eq:IPP_final} but for the term 
\begin{align*}
    \E\left [\Tj^{\sigma, (p,q)}_{j }(\phi)(\Fvect)\overline{\Gamma}_{p,q} \L^{j}(\Gamma(F,-\L^{-1}F)-1) \right ]
\end{align*}
appearing in the former expansion. As a result, we should obtain an expansion of the form
\begin{align*}
    \E \left [ \phi(\Fvect)(\Gamma(F-\L^{-1}F)-1) \right ] = \sum_{j,k,p_i,q_i} \E \left [ \Tj_{j,k}^{(p_1,q_1,p_2,q_2)}(\phi)(\Fvect) \overline{\Gamma}_{p_2,q_2} \L^k \left [\overline{\Gamma}_{p_1,q_1}\L^j [\Gamma(F,-\L^{-1}F)-1] \right ] \right ] \, ,
\end{align*}
for $\phi$ with a high enough Hermite rank, where $\Tj_{j,k}^{(p_1,q_1,p_2,q_2)}$ is a fraction evaluated at $\LL_\sigma$. Then we can bound the expansion terms using the same arguments as in Section \ref{sec.bounding}, by a power of $\sum_{1\leq p \leq q \leq d} \Var\Gamma(F_p,F_q)$ times $(\Var\Gamma(F,\L^{-1}F) )^{1/2}$.
\subsection{Multivariate expansion}\label{ss:multi_expansion}
In this section, we generalize the expansion \eqref{eq:IPP_final} by swapping $\Gamma(F,-\L^{-1} F)-1$ by an element $W \in \WW_{d_m}$, where $d_m$ denotes the integer  $(2d-2)m$. First, let us introduce some further notations. . Furthermore, denote by $\mathrm{J}_m$ the set of elements $(j_1, \ldots, j_m)$ such that
\begin{align*}
    \forall k \in \llbracket 1, m \rrbracket, \quad j_k \in \llbracket 0,d_k \rrbracket \, .
\end{align*}
We will denote as well by $\mathrm{I}$ the set $\{1,\ldots, d\}^2$. Recall that, in the same spirit as Lemma \ref{prop:Pannulateur}, the polynomial $P_{d_m}(X) = X(X+1)\cdots (X+d_m)$ verifies
\begin{align*}
\forall W \in \WW_{\leq d_m}, \quad P_{d_m}(\L)(W)=0.
\end{align*}
As before, we denote $P_{d_m} = {\sum}^{d_m+1}_{l=1}a^{(m)}_{l}X^l$. 

\medskip

First, we write a more general than the result than Proposition \ref{proposition:IPP_polynomial} by noticing that we replace $\Ga(F,\L^{-1} F)-1$ by $W$ and $P$ by $P_{d_m}$ in the proofs of Lemma \ref{lemma:IPPalpha} and Proposition \ref{proposition:IPP_polynomial} to get the same result.

\begin{proposition}\label{proposition:IPP_polynomial_HO}
    Let $\phi$ be an element of $\Cpold{2d_m+4}{d}$ and let $W$ be an element of $\WW_{\leq d_m}$. Then one has
\begin{align}
   \E[P_{d_m}(\LL_\sigma)(\phi) (\Fvect)W]   = -\sum_{j=0}^{d_m} \sum_{p,q=1}^d \, \E\left [\partial_{pq}^2 P_{d_m,j} \left (\LL_{\sigma} \right )(\phi) (\Fvect) \overline{\Gamma}_{p,q}\L^{j}W \right ],
\end{align}
where $P_{d_m, j}(X) := \sum_{\alpha=1+j}^{d_m+1} a^{(m)}_{\alpha} X^{\alpha-1-j}$.
\end{proposition}

Then, inverting $P_{d_m}(\LL_\sigma)$ on functions $\phi \in \Cpold{2d_m+4}{d}$ such that $\rank_\sigma(\phi) \geq d_m+1$ leads to the next result, generalizing Proposition \ref{proposition:IPP_final_d}.

\begin{proposition}\label{proposition:gen.expansion.d}
    Let $\phi$ be an element of $\Cpold{2d_m+4}{d}$ such that $\rank_\sigma(\phi) \geq d_m+1$ and let $W$ be an element of $\WW_{\leq d_m}$. Then one has
\begin{align}\label{eq:gen.expansion.d}
   \E[\phi(F_1, \ldots, F_d)W]   = \sum_{j=0}^{d_m} \sum_{\pi \in \mathrm{I}} \, \E\left [\Tj_{m,j}^{\sigma, \, \pi} (\phi) (F_1,\ldots,F_d) \overline{\Gamma}_{\pi}\L^{j}W \right ],
\end{align}
where $\overline{\Gamma}_{\pi} := \overline{\Gamma}_{p,q}$, and $
\Tj_{m,j}^{\sigma, \, \pi} (\phi) \coloneq - P_{d_m,j}(\LL_\sigma-p-q) \bigl ( P_{d_m}(\LL_\sigma - p-q) \bigr )^{-1} \partial_{pq}^2 \phi$ for $\pi=(p,q)$.
\end{proposition}

We retrieve a similar factorization for $\Tj_{m,j}^{\sigma, \, \pi}$ than the one in Lemma \ref{lemma:rewrite_of_T_d}, stated in the next Lemma. Its proof is identical to that of Lemma \ref{lemma:rewrite_of_T_d}.

\begin{proposition}\label{proposition:factorisation_Tmi_d}
    Let $j$ be an element of $\llbracket 0, d_m \rrbracket$ and let $\pi = (p,q)$ be an element of $\mathrm{I}$. Then, for every $\phi \in \Cpold{2d_m+4}{d}$ such that $\rank_\sigma(\phi)\geq d_m+1$, we have
    \begin{align}\label{eq:factorisation_Tmi_d}
        \Tj^{\sigma, \, \pi}_{m,j}(\phi) & = -\Q^{\pi}_{m,j} \left ( \LL_{\sigma} \right ) (\LL_\sigma-p-q)^{-1}\partial_{pq}^2 \phi \, ,
    \end{align}
    where 
    \begin{align*}
        \Q_{m,j}^{\pi}(X) := (X-p-q)^{-j} - \underset{\alpha=1}{\overset{j}{\sum}} \, a_\alpha^{(m)} (X-p-q)^{\alpha-j}P_{d_m} ( X-p-q  )^{-1} \, .
    \end{align*} 
\end{proposition}

Since the factorization is stated for $\phi:\R^d \to \R$, it will allow us to use it for composition of operators of the form $\Tj_{m,j}^{\sigma, \, \pi}$, which we would not have been able to do if it were stated for $\phi = [\varphi]^\oplus$.

\medskip

Since we want to iterate the use of \eqref{proposition:IPP_final_extended}, we need to compare the Hermite rank of $\Tj_{m,j}^{\sigma, \pi}(\phi)$ with the Hermite rank of $\phi$. This is the content of the next lemma.
\begin{lemma}\label{lemma:rank_T}
Let $j$ be an element of $\llbracket 0, d_m \rrbracket$ and let $\pi$ be an element of $\mathrm{I}$. Then, for every $\phi\in\Cpold{2d_m+4}{d}$,
\begin{align*}
\rank_\sigma \left ( \Tj^{\sigma,\pi}_{m,j}(\phi) \right ) \geq \rank_\sigma(\phi)-2 \, .
\end{align*}
\end{lemma}
\noindent The proof of Lemma \ref{lemma:rank_T} follows from adapting Steps 2 and 3 in the proof of \cite[Proposition 6.13]{MPS2025.EE} to the multidimensional case.

\subsection{Higher order expansion}
In this section, we derive the general expansion of \begin{align*}
    \E[\varphi(F)(\Gamma(F,-\L^{-1} F)-1)]
\end{align*} 
as a sum of terms controlled by
\begin{align*}
\left ( \sum_{1\leq p \leq q \leq d} \Var\Gamma(F_p,F_q) \right )^{\frac{m}{2}} \bigl ( \Var\Gamma(F,\L^{-1} F) \bigr )^{\frac{1}{2}}\,,
\end{align*}
provided that $\varphi$ is sufficiently smmoth and has a sufficiently high Hermite rank. To this end, we apply Proposition \ref{proposition:gen.expansion.d} iteratively. Starting from $\E[\phi(\Fvect)(\Gamma(F,-\L^{-1} F)-1)]$, we write
\begin{eqnarray*}
       &&\E \left [ \phi(\Fvect)(\Gamma(F-\L^{-1}F)-1) \right ] = \sum_{p,q=1}^d \sum_{j=0}^{2d-2} \, \E\left [\Tj^{\sigma, (p,q)}_{j }(\phi)(\Fvect)\overline{\Gamma}_{p,q} \L^{j}(\Gamma(F,-\L^{-1}F)-1) \right ].
\end{eqnarray*}
Then, for each term
\begin{align*}
\E\left [\Tj^{\sigma, (p,q)}_{j }(\phi)(\Fvect)\overline{\Gamma}_{p,q} \L^{j}(\Gamma(F,-\L^{-1}F)-1) \right ]\,,
\end{align*}
we apply \eqref{eq:gen.expansion.d}, and repeat the procedure at each step on the terms generated by this expansion. Each iteration requires appropriate conditions on both the regularity and the Hermite rank of $\phi$. The proof follows closely the proof of \cite[Proposition 6.13]{MPS2025.EE} and leads to the following result.
\begin{proposition}\label{proposition:IPP_final_HO}
    Let $\phi$ be an element of $\Cpold{2d_m+4}{d}$ such that $\rank_\sigma(\phi) \geq d_m+2m-1$. Then one has
\begin{align}\label{eq:IPP_final_sum_HO}
   \E[\phi(F_1,\ldots, F_d)(\Gamma(F,-\L^{-1}F)-1)]   =   \sum_{\pi \in \mathrm{I}^m} \sum_{j \in \J_{m}}\, \E\left [\boldsymbol{\Tj}_{m,j}^{\sigma, \, \pi} (\phi) (F_1,\ldots, F_d)W_{j}^{\pi} \right ],
\end{align}
where $\boldsymbol{\Tj}_{m,\,j}^{\sigma, \, \pi}(\phi) := \Tj_{m,j_m}^{\sigma, \pi_m} \cdots \Tj_{1,j_1}^{\sigma, \pi_1} (\phi)$ and 
\begin{align*}
W_{j}^{\pi} := \overline{\Gamma}_{p_{m},q_{m}} \L^{j_m} \left [ \cdots \left [ \overline{\Gamma}_{p_{1},q_{1}}  \L^{j_{1}} \left [ \Gamma(F,-\L^{-1} F)-1 \right ] \right ] \cdots  \right ],
\end{align*}
if $\pi = ((p_1,q_1), \ldots, (p_{m}, q_{m}))$.
\end{proposition}

We now express $\boldsymbol{\Tj}_{m,\,j}^{\sigma, \, \pi}$ on functions $\phi$ of the form $[\varphi]^\oplus$ with $\varphi:\R \to \R$, using the factorizations \eqref{eq:factorisation_Tmi_d}. As before, we will write $\boldsymbol{\Tj}_{m,\,j}^{\sigma, \, \pi}(\varphi)$ for $\boldsymbol{\Tj}_{m,\,j}^{\sigma, \, \pi}[\varphi]^\oplus$ and identify that way $\boldsymbol{\Tj}_{m,\,j}^{\sigma, \, \pi}$ as an operator on function from $\R$ to $\R$. 

\begin{proposition}\label{prop:facto}
    Let $j$ be an element of $\J_m$ and let $\pi$ be an element of $\mathrm{I}^m$. Then, for every $\varphi \in \Cpol{2d_m+4}$ such that $\rank(\varphi)\geq d_m+2m-1$, we have
    \begin{align}
        \boldsymbol{\Tj}^{\sigma, \, \pi}_{m,j}(\varphi) & := (-1)^m\mathbf{Q}_{m,j}^{\pi} \left ( \LL_{\sigma} \right ) \S_\sigma^{\pi,m}(\varphi) \, ,
    \end{align}
    where 
    \begin{align*}
        \mathbf{Q}_{m,j}^{\pi}(X) := \prod_{k=1}^m  \Q_{k,j_k}^{\pi_k} \left ( X - \sum_{\ell =k+1}^{m} |\pi_{\ell}| \right ) \, ,
    \end{align*}
    and
    \begin{align*}
        \S_\sigma^{\pi,m} (\varphi) := \left (\prod_{k=1}^m \left (\LL_\sigma-\sum_{\ell = 0}^{k-1} |\pi_{m-\ell}| \right ) \right )^{-1}[\varphi^{(2m)}]^\oplus \, .
    \end{align*}
\end{proposition}
For easier expressions, in the following we set
\begin{align}\label{notations-nu-rho}
    \nu_{\pi,k} \coloneq \sum_{\ell = k+1}^m |\pi_\ell|\, , \quad \rho_{\pi,k} \coloneq \sum_{\ell = 0}^{k-1} |\pi_{m-\ell}|\,.
\end{align}

\subsection{Extending operators}

In this section, we aim to extend the operator $\boldsymbol{\Tj}^{\sigma, \, \pi}_{m,j}$ to functions $\varphi \in \Cpol{0}$ with Hermite rank at least $d_m+2m-1$. In view of Proposition \ref{prop:facto}, it is enough to extend the operator $\S_{\sigma}^{\pi,m}$, initially defined on $\Cpol{2m}$ to an operator acting on functions in $\Cpol{0}$. In Section \ref{ss:studyofTjk} we already proved the case $m=1$. We showed that
\begin{eqnarray*}
&& \S_\sigma^{1,(p,q)}(\varphi)(x) = \\
&& \frac{1}{\E[\Gamma(F)]} \left ( [\varphi]^\oplus(x) + (p+q)\RR_\sigma(p+q)[\varphi]^\oplus(x) + \RR_{\sigma}^{V_2}(p+q)[\varphi]^\oplus(x) + \sum_{b=1}^d \frac{x_b}{\sigma_b} \RR_{\sigma,b}^{V_{1,b}}(p+q)[\varphi]^\oplus(x) \right ),
\end{eqnarray*}
for every $x\in\R^d$, where, for $V$ polynomial,
\begin{align*}
    \RR_\sigma^{V}(p+q)[\varphi]^\oplus(x) := -\int_0^{+\infty} e^{-(p+q)t} \, \Pt_{t}^{\sigma, V}[\varphi]^\oplus(x) \, \dint t \, ,
\end{align*}
and, for $b \in \llbracket 1,d\rrbracket$,
\begin{align*}
\RR_{\sigma,b}^{V}(p+q)[\varphi]^\oplus(x) \coloneq -\int_0^{+\infty} \frac{e^{-(p+q+b)t}}{\left ( 1-e^{-2bt} \right )^{1/2}} \Pt_t^{\sigma,V}[\varphi]^\oplus(x) \, \dint t \, .
\end{align*}
The main idea was to rewrite $[\varphi'']^\oplus$ as 
\begin{align}\label{trick-on-varphi}
 [\varphi'']^{\oplus}(x_1,\ldots,x_d) &= m_{\Gamma}^{-1} \left (  \LL_d [\varphi]^\oplus(x_1,\ldots,x_d)+\sum_{p=1}^d px_p [\varphi']^\oplus(x_1,\ldots,x_d) \right ),
\end{align}
where $m_\Gamma := \E \left [\Gamma(F)\right ]$, and then use the property of the semigroup $\Pt_t^{\sigma}$ in order to make the operator $\LL_\sigma$ and the first derivative of $\varphi$ disappear. Namely, we obtained the following result.
\begin{lemma}\label{lemma:trick-on-semigroup}
    Let $\varphi$ be an element of $\Cpol{2}$. Then, for every $x\in \R^d$ and $t\geq 0$,
    \begin{align}\label{eq:trick-on-semigroup}
    \Pt_t^\sigma[\varphi'']^\oplus(x) = \frac{1}{m_\Gamma} \left ( \partial_t \Pt^\sigma_t [\varphi]^\oplus( x)+ \Pt_t^{\sigma,V_{2}} [\varphi]^\oplus(x) + \sum_{b=1}^d \frac{x_b}{\sigma_b} \frac{e^{-bt}}{\left (1-e^{-2bt} \right )^{1/2}} \Pt_t^{\sigma, V_{1,b}} [\varphi]^\oplus(x) \right) \, ,
\end{align}
where $m_\Gamma = \E[\Gamma(F)]$.
\end{lemma}
Iterating \eqref{trick-on-varphi} for each term appearing in \eqref{eq:trick-on-semigroup} would lead to a lengthy and technically involved argument, mainly because the semigroups now carry polynomial weights, which substantially changes the analysis. We therefore apply \eqref{eq:trick-on-semigroup} only to the term 
\begin{align*}
    \partial_t \Pt^\sigma_t[\varphi]^\oplus(x) \, ,
\end{align*}
after an integration by parts with respect to $t$, and handle the remaining terms separately. To this end, we use the following representation of the operator $\S_\sigma^{\pi,m}$, which follows directly from Proposition~\ref{proposition:resolvent_d}.

\begin{lemma}\label{lemma:other_form_S_d}
    Let $\pi$ be an element of $\mathrm{I}^m$ and let $\varphi$ be an element of $\Cpol{2m}$. Then
    \begin{align}\label{eq:other_form_S_d}
       \forall x \in \R^d, \quad  \S_\sigma^{\pi,m}(\varphi)(x)  = (-1)^m \int_0^{+\infty} K_m(t) \Pt_t^\sigma [\varphi^{(2m)}]^\oplus(x) \, \dint t \, ,
    \end{align}
    where $K_m(t) = e_{\rho_{\pi,k}}^+ \ast \cdots \ast e_{\rho_m}^+(t)$, $e_{\rho_{\pi,k}}^+(t) := e^{-\rho_{\pi,k} t}\I_{\R_+}(t)$ and $\rho_{\pi,k}(t) := \sum_{\ell = 0}^{k-1} |\pi_{m-\ell}|$.
\end{lemma}
\begin{proof}
    Let $x$ be an element of $\R^d$. Using Proposition \ref{proposition:resolvent_d} and the semigroup property, one has
    \begin{align*}
        \S_\sigma^{\pi,m}(x) & = \left (\prod_{k=1}^m \left (\LL_\sigma-\rho_{\pi,k} \right ) \right )^{-1}[\varphi^{(2m)}]^\oplus(x) = (-1)^m \int_{\R_+^m} e^{-\sum_{k=1}^m \rho_{\pi,k} t_k} \, \left (\Pt_{t_1}^\sigma \cdots \Pt_{t_m}^\sigma \right )[\varphi^{(2m)}]^\oplus(x) \, \dint t_1 \cdots \dint t_m \\
        & = (-1)^m \int_{\R_+^m} e^{-\sum_{k=1}^m \rho_{\pi,k} t_k} \, \ \Pt_{t_1 +\cdots + t_m}^\sigma [\varphi^{(2m)}]^\oplus(x) \, \dint t_1 \cdots \dint t_m \, .
    \end{align*}
    Then using the change of variable 
    \begin{align*}
        t = t_1 + \cdots + t_m, \qquad \forall k \in \{1,\ldots, m-1\}, \quad s_k = t_k \, ,
    \end{align*}
    we obtain
    \begin{align*}
        \S_\sigma^{\pi,m}(x) & = (-1)^m \int_{\Delta^{m-1}(t)}\int_{0}^{+\infty} e^{-\sum_{k=1}^{m-1} \rho_{\pi,k}s_k}e^{-\rho_m \left (t-\sum_{k=1}^{m-1} s_k \right )} \, \ \Pt_{t}^\sigma [\varphi^{(2m)}]^\oplus(x) \, \dint t \, \dint s_1 \cdots \dint s_{m-1} \, ,
    \end{align*}
    where $\Delta^{m-1}(t) := \{ s \in \R_+^{m-1} \; / \; \sum_{k=1}^{m-1} s_k \leq t \}$. Now, setting
    \begin{align*}
        K_m(t) & := \int_{\Delta^{m-1}(t)} e^{-\sum_{k=1}^{m-1} \rho_{\pi,k}s_k}e^{-\rho_m \left (t-\sum_{k=1}^{m-1} s_k \right )} \, \dint s_1 \cdots \dint s_{m-1} \\
        & =  \int_{\Delta^{m-1}(t)} \left (\prod_{k=1}^{m-1} e_{\rho_{\pi,k}}^+(s_k) \right )e_{\rho_m}^+ \left (t-\sum_{k=1}^{m-1} s_k \right ) \, \dint s_1 \cdots \dint s_{m-1} = e_{\rho_{\pi,k}}^+ \ast \cdots \ast e_{\rho_m}^+(t) \, ,
    \end{align*}
    we get
    \begin{align*}
        \S_\sigma^{\pi,m}(\varphi)(x)  = (-1)^m \int_0^{+\infty} K_m(t) \Pt_t^\sigma [\varphi^{(2m)}]^\oplus(x) \, \dint t \, ,
    \end{align*}
    as desired.
\end{proof}
Before extending the operator $\S_\sigma^{\pi,m}$, we establish several properties of $K_m$ and its derivatives. We also derive a standard formula for the derivatives of the semigroup $\Pt_t^{\sigma,V}[\varphi]^\oplus$, analogous to \eqref{eq:derive_PtV}, that exploits the particular structure of the function $[\varphi]^\oplus$.
\begin{lemma}\label{lemma:asympt_Km_0}
   The function $K_m$ is infinitely differentiable on $]0,+\infty[$. Furthermore, for every $j\in \llbracket 0, m-1 \rrbracket$,
    \begin{align*}
        K_m^{(j)}(t) = \frac{t^{m-1-j}}{(m-1-j)!} + O(t^{m-j}) \, , \quad \text{ as } t \to 0 \, .
    \end{align*}
    In particular, there exist $B_{d,m,j} >0$ such that
    \begin{align*}
        \forall t \in ]0,1], \quad | K_m^{(j)}(t)| \leq B_{d,m,j} t^{m-1-j} \, .
    \end{align*}
\end{lemma}
\begin{proof}
    Let us first prove the case $j=0$. After a change of variable $t u_k=s_k$, we get
    \begin{align}\label{eq:formule-Km}
        K_m(t) & = t^{m-1} \int_{\Delta^{m-1}(1)} e^{-t\sum_{k=1}^{m-1} \rho_{\pi,k}u_k}e^{-\rho_{\pi,m} t\left (1-\sum_{k=1}^{m-1} u_k \right )} \, \dint u_1 \cdots \dint u_{m-1} \nonumber \\
        & = t^{m-1} \int_{\Delta^{m-1}(1)} e^{-tA(u)} \, \dint u_1 \cdots \dint u_{m-1} \, ,
    \end{align}
    where $A(u) := \sum_{k=1}^{m-1} \rho_{\pi,k}u_k+\rho_{\pi,m} \left (1-\sum_{k=1}^{m-1} u_k \right )$. Since $A(u)$ is a convex combination of the $\rho_{\pi,k} \geq 0$, we have $0\leq A(u) \leq \underset{1\leq k\leq m}{\max} \rho_{\pi,k}$, and by a standard inequality, we get
    \begin{align*}
        \forall t \in \R_+, \quad 1-e^{-tA(u)} \leq t\underset{1\leq k\leq m}{\max} \rho_{\pi,k} = t \rho_{\pi,m} \, .
    \end{align*}
    Noticing that $|\Delta^{m-1}(1)| = \tfrac{1}{(m-1)!}$, it follows 
    \begin{align*}
        \left | K_m(t)-\frac{t^{m-1}}{(m-1)!} \right | \leq t^{m-1} \int_{\Delta^{m-1}(1)} |e^{-tA(u)}-1| \, \dint u_1 \cdots \dint u_{m-1} \,  \leq \frac{ t^m}{(m-1)!} \rho_{\pi,m} \, ,
    \end{align*}
    hence the result for $j=0$. Now, let $j$ be an element of $\llbracket 1, m-1\rrbracket$. Denote by $f_m(t)$ the integral $\int_{\Delta^{m-1}(1)} e^{-tA(u)} \, \dint u_1 \cdots \dint u_{m-1}$. Since $\Delta^{m-1}(1)$ is compact and $A$ is bounded, $f_m$ is infinitely differentiable on $]0,+\infty[$ and differentiating under the integral sign leads to
    \begin{align*}
        f_m^{(j)}(t) = (-1)^j \int_{\Delta^{m-1}(1)} A(u)^j e^{-tA(u)} \, \dint u_1 \cdots \dint u_{m-1} \, .
    \end{align*}
    In particular, the derivatives of $f_m$ are all bounded on $]0,1]$, and $K_m$ is also infinitely differentiable on $]0,+\infty[$. By Leibniz's formula, one has
    \begin{align*}
        K_m^{(j)}(t) & = \sum_{\ell=0}^j \binom{j}{\ell} \frac{(m-1)!}{(m-1-j+\ell)!} t^{m-1-j+\ell} f_m^{(\ell)}(t) \\
        & = \frac{(m-1)!}{(m-1-j)!}t^{m-1-j}f_m(t) +t^{m-1-j+1}\sum_{\ell=1}^j \binom{j}{\ell} \frac{(m-1)!}{(m-1-j+\ell)!} t^{\ell-1} f_m^{(\ell)}(t) \, .
    \end{align*}
    Now, since $f_m(0) = \tfrac{1}{(m-1)!}$, we deduce the result in general.
\end{proof}

\begin{lemma}\label{lemma:asympt_Km_infty}
    Let $j$ be an element of $\llbracket 0, m-1 \rrbracket$. Then there exists $C_{m,j}>0$ such that
    \begin{align*}
        \forall t \geq 1, \quad |K_m^{(j)}(t)| \leq C_{m,j}t^{m-1} e^{-\rho_{\pi,1} t} \, .
    \end{align*}
\end{lemma}
\begin{proof}
    Recall from the previous proof that 
    \begin{align*}
        K_m^{(j)}(t) & = \sum_{\ell=0}^j \binom{j}{\ell} \frac{(m-1)!}{(m-1-j+\ell)!} t^{m-1-j+\ell} f_m^{(\ell)}(t) \, ,
    \end{align*}
    where $f_m^{(j)}(t) = (-1)^j \int_{\Delta^{m-1}(1)} A(u)^j e^{-tA(u)} \, \dint u_1 \cdots \dint u_{m-1}$. For all $u \in \Delta^{m-1}(1)$, one has 
    \begin{align*}
        A(u) = \sum_{k=1}^{m-1} \rho_{\pi,k} u_k + \rho_{\pi,m} \left ( 1 - \sum_{k=1}^{m-1} u_k \right )\geq \left ( \sum_{k=1}^{m-1} u_k + 1 -\sum_{k=1}^{m-1} u_k \right )\underset{1\leq k\leq m}{\min} \rho_{\pi,k} = \rho_{\pi,1} \, ,
    \end{align*}
    hence 
    \begin{align*}
        |f_m^{(\ell)}(t)| \leq e^{-\rho_{\pi,k}} \int_{\Delta^{m-1}(1)} A(u)^\ell \, \dint u \leq \frac{\rho_{\pi,m}^\ell  }{(m-1)!}e^{-\rho_{\pi,1} t} \, .
    \end{align*}
  The triangular inequality yields then
    \begin{align*}
       |K_m^{(j)}(t)| \leq \sum_{\ell=0}^j \binom{j}{\ell} \frac{(m-1)!}{(m-1-j+\ell)!} t^{m-1-j+\ell} |f_m^{(\ell)}(t)| \leq \left ( \sum_{\ell=0}^j \frac{\binom{j}{\ell} \rho_{\pi,m}^{\ell} }{(m-1-j+\ell)!} \right ) t^{m-1} e^{-\rho_{\pi,1} t}\, ,
    \end{align*}
    for every $t \geq 1$.
\end{proof}

Before stating the next lemma, let's settle some notations first. Using Mehler's formula \eqref{eq:mehler_formula_d}, we write
\begin{align*}
    \Pt_t^{\sigma, V} [\varphi]^\oplus(x) = \E \left [ V(Z)\varphi (\beta_t(x) + \sqrt{v_t}u_t \cdot Z )\right ] \, ,
\end{align*}
where 
\begin{align*}
    \beta_t(x) :=\sum_{p=1}^d e^{-pt}x_p\, , \quad v_t := \sum_{p=1}^d \sigma_p^2(1-e^{-2pt}) \, , \quad  u_t := v_t^{-\frac{1}{2}}(\sigma_1(1-e^{-2t})^{\frac{1}{2}}, \ldots, \sigma_d (1-e^{-2dt})^{\frac{1}{2}}) \, .
\end{align*}
In particular, one has $\|u_t\| = 1$.  

\begin{lemma}\label{lemma:formula_PtV_deriv}
    Let $\varphi$ be an element of $\Cpol{2m}$ and let $V$ be a nonzero polynomial. For every $j$ element of $\llbracket 0, 2m \rrbracket$, one has
    \begin{align*}
        \forall t >0, \; \forall x \in \R^d, \quad \Pt_t^{\sigma, V} [\varphi^{(j)}]^\oplus(x) = v_t^{-\frac{j}{2}} \Pt_t^{\sigma, A_{j,t}(V)}[\varphi]^\oplus(x) \, ,
    \end{align*}
    where 
    \begin{align*}
        A_{j,t}(V)(z) := \sum_{\ell = 0}^j (-1)^{j-\ell} \binom{j}{\ell} H_\ell(u_t \cdot z) \partial_{u_t}^{j-\ell}V(z) \, .
    \end{align*}
\end{lemma}
\begin{proof}
Consider the function $g$ defined by
    \begin{align*}
        g(z) := \varphi (\beta_t(x) + \sqrt{v_t}u_t \cdot z ) \, .
    \end{align*}
    Then $g$ is in $\Cpold{2m}{d}$ and 
    \begin{align*}
        \partial_{u_t} g(z) = \sum_{p=1}^d \sqrt{v_t} (u_t)_p^2 \varphi'( \beta_t(x)+\sqrt{v_t}u_t \cdot z) = \sqrt{v_t} \, \varphi'( \beta_t(x)+\sqrt{v_t}u_t \cdot z) \,.
    \end{align*}
    Hence, more generally,
    \begin{align*}
        \partial_{u_t}^j g(z) =(\sqrt{v_t})^j \, \varphi^{(j)}( \beta_t(x)+\sqrt{v_t}u_t \cdot z)\,.
    \end{align*}
    Furthermore, using the Gaussian integration by parts, we get
    \begin{align*}
        \E[V(Z) \partial_{u_t}f(Z)] & = \sum_{p=1}^{d} (u_t)_p\, \E[V(Z)\partial_pf(Z)] = \sum_{p=1}^{d} (u_t)_p\, \E[(Z_pV(Z)-\partial_p V(Z) )f(Z)] \\
        & =\E[((u_t \cdot Z)V(Z)-\partial_{u_t} V(Z) )f(Z)]\, .
    \end{align*}
    We deduce that 
    \begin{align*}
         \Pt_t^{\sigma, V} [\varphi^{(j)}]^\oplus(x) = v_t^{-\frac{j}{2}} \, \E[A_{j,t}(V)(Z) \varphi (\beta_t(x) + \sqrt{v_t}u_t \cdot Z )] \, ,
    \end{align*}
    where $A_{j,t}(V)$ is obtained by differentiating $j$ times $V(z)\gamma_d(z)$ in the direction of $u_t$ ($\gamma_d$ here stands for the density of the standard Gaussian distribution on $\R^d$).
\end{proof}

We leverage now the previous lemmas to extend the operator $\S_\sigma^{\pi,m}$ to an operator acting on $\Cpol{0}$, with explicit pointwise polynomial control of the image.

\begin{proposition}\label{proposition:prolongement_S_d_m}
    Let $\pi$ be an element of $\mathrm{I}^m$. The operator $\S_\sigma^{\pi,m}$ can be extended to an operator on functions $\varphi \in \Cpol{0}$. Furthermore, if $\varphi \in \Cpol{0}$ verifies 
\begin{align}\label{ineq:estimates_polynomial_varphi_oplus}
\forall x \in \R^d,\quad \left | [\varphi]^\oplus(x)\right | \leq C \left ( 1 +\| I_\sigma^{-1} x \|^l \right ),
\end{align}
for some $C>0$ and $l \in \N$, then $\S_\sigma^{\pi,m}(\varphi)\in \Cpold{0}{d}$ and there exist $A_{l,m}^\pi>0$ such that
\begin{align}\label{ineq:estimate_S_d_m}
\forall x \in \R^d, \quad \left | \S_\sigma^{\pi,m}(\varphi)(x) \right | \leq CA_{l,m}^\pi \left ( 1 + \|I_\sigma^{-1} x \|^{l+1} \right ).
\end{align}
\end{proposition}

\begin{proof}
Let $\varphi$ be an element of $\Cpol{2m}$, let $x$ be an element of $\R^d$. Using first Lemma \ref{lemma:other_form_S_d} and then Lemma \ref{lemma:trick-on-semigroup} to rewrite the semigroup in the integral, we get
\begin{eqnarray*}
     && \S_\sigma^{\pi,m}(\varphi)(x) = (-1)^m \int_0^{+\infty} K_m(t) \Pt_t^\sigma [\varphi^{(2m)}]^\oplus(x) \, \dint t \\
    && = \frac{(-1)^m}{m_\Gamma}\int_0^{+\infty} K_m(t) \left ( \partial_t \Pt^\sigma_t [\psi]^\oplus( x) + \Pt_t^{\sigma,V_{2}} [\psi]^\oplus(x) + \sum_{b=1}^d \frac{x_b}{\sigma_b} \frac{e^{-bt}}{\left (1-e^{-2bt} \right )^{1/2}} \Pt_t^{\sigma, V_{1,b}} [\psi]^\oplus(x) \right ) \, \dint t \, ,
\end{eqnarray*}
where $\psi := \varphi^{(2m-2)}$. Then, using Lemma \ref{lemma:formula_PtV_deriv}, one has 
\begin{eqnarray*}
     && \mathcal{A}_t^{(m-1)}(\varphi)(x) := \Pt_t^{\sigma,V_{2}} [\psi]^\oplus(x) + \sum_{b=1}^d \frac{x_b}{\sigma_b} \frac{e^{-bt}}{\left (1-e^{-2bt} \right )^{1/2}} \Pt_t^{\sigma, V_{1,b}} [\psi]^\oplus(x) \\
     && =  v_t^{-m+1} \Pt_t^{\sigma, A_{2m-2,t}(V_2)}[\varphi]^\oplus(x) + \sum_{b=1}^d \frac{x_b}{\sigma_b} \frac{e^{-bt}}{\left (1-e^{-2bt} \right )^{1/2}} v_t^{-m+1} \Pt_t^{\sigma, A_{2m-2,t}(V_{1,b})}[\varphi]^\oplus(x) \, .
\end{eqnarray*}
Using Lemmas \ref{lemma:asympt_Km_0}, \ref{lemma:asympt_Km_infty} and \ref{lemma:formula_PtV_deriv}, we get that $t \mapsto K_m(t)\mathcal{A}_t^{(m-1)}(\varphi)(x)$ is integrable with respect to $t$ on $[0,+\infty[$, since $v_t \sim ct$, as $t\to 0$, for some $c>0$. Then, by expanding and integrating by parts, we get
\begin{align*}
     \S_\sigma^{\pi,m}(\varphi)(x) & = \frac{(-1)^{m+1}}{m_\Gamma}\int_0^{+\infty} K'_m(t) \Pt^\sigma_t [\psi]^\oplus( x) \, \dint t + \frac{(-1)^m}{m_\Gamma}\int_0^{+\infty} K_m(t) \mathcal{A}_t^{(m-1)}(\varphi)(x) \, \dint t \, ,
\end{align*}
where the additional terms from the integration by parts vanish in virtue of Lemmas \ref{lemma:asympt_Km_0} and \ref{lemma:asympt_Km_infty}. Then, iterating $(m-2)$ times the same argument on the integral $\int_0^{+\infty} K'_m(t) \Pt^\sigma_t [\psi]^\oplus( x) \, \dint t$, we get
\begin{align}\label{eq:prolong-S-d}
     \S_\sigma^{\pi,m}(\varphi)(x) &= \frac{1}{m_\Gamma^{m}} \left ( [\varphi]^\oplus(x) + \int_0^{+\infty} K_m^{(m)}(t)\Pt_t^{\sigma}[\varphi]^\oplus(x) \, \dint t \right ) \\
     & \qquad + \sum_{j=0}^{m-1} \frac{(-1)^{m+j}}{m_\Gamma^{j+1}}\int_0^{+\infty} K_m^{(j)}(t) \mathcal{A}_{t}^{(m-1-j)}(\varphi)(x) \, \dint t \, , \nonumber
\end{align}
where we get an extra $[\varphi]^\oplus$ from the last integration by parts, since $K_{m}^{(m-1)}(0)=1$, and where
\begin{align*}
    \mathcal{A}_t^{(k)}(\varphi)(x) := v_t^{-k} \Pt_t^{\sigma, A_{2k,t}(V_2)}[\varphi]^\oplus(x) + \sum_{b=1}^d \frac{x_b}{\sigma_b} \frac{e^{-bt}}{\left (1-e^{-2bt} \right )^{1/2}} v_t^{-k} \Pt_t^{\sigma, A_{2k,t}(V_{1,b})}[\varphi]^\oplus(x) \,.
\end{align*}
With the sale arguments than before, the function $t\mapsto K_m^{j}(t)\mathcal{A}^{m-1-j}_t(\varphi)(x)$ is integrable on $[0,+\infty[$, for all $j$. Since the right-hand side of \eqref{eq:prolong-S-d} is a linear operator in $\varphi$, defined for $\varphi \in \Cpol{0}$ (the necessary intergability conditions follows from the usual arguments), it yields a linear extension of $\S_\sigma^{\pi,m}$ on $\Cpol{0}$.

The polynomial estimates \eqref{ineq:estimate_S_d_m} follows from Lemma \ref{lemma:estimates_PtV_pol} applied on $\Pt_t^\sigma [\varphi]^\oplus$ and on each $\mathcal{A}_t^{(k)}(\varphi)$, using the fact that each function $t\mapsto K_m^{(j)}(t)$ is integrable on $[0,+\infty[$. With the same estimations, combined with the continuity under the integral sign theorem, we obtain that $\S_\sigma^{\pi, m}(\varphi)$ More precisely, recall the bound from Lemma \ref{lemma:estimates_PtV_pol}: for a nonzero polynomial $V$, one has 
\begin{align*}
    |\Pt_t^{\sigma, V} [\varphi]^\oplus(x)| \leq B_{l,V} \left (1+\|I_\sigma^{-1}x \|^l \right )\,, \quad B_{l,V} =2^{l-1} \,\E \left [ \left | V(I_\sigma^{-1}N_\sigma ) \right | \left ( 1 + \|I_\sigma^{-1}N_\sigma\|^l \right ) \right ] \, .
\end{align*}
Since the polynomials $V$ occurring in the expression of $\mathcal{A}_t^{(k)}(\varphi)$ depend on $t$, the bound $B_{l,V}$ might, at first sight, also appear to depend on $t$. The argument below shows, however, that it can be chosen uniformly in $t$. By the triangular inequality, we have
\begin{align*}
     |A_{j,t}(V)(z)| \leq  \sum_{\ell = 0}^j \binom{j}{\ell} |H_j(u_t \cdot z)| |\partial_{u_t}^{j-\ell}V(z)| \,.
\end{align*}
On one side, for a polynomial $Q = \sum_{k=0}^D q_kX^k$, the Cauchy-Schwarz inequality applied to the monomials reaps
\begin{align*}
    |Q(u_t \cdot z)| \leq \sum_{k=1}^D |q_k| |u_t \cdot z|^k )\leq \sum_{k=1}^D |q_k| \|u_t\|^k \|z\|^k = \sum_{k=1}^D |q_k|  \|z\|^k \, .
\end{align*}
Similarly, $\partial_{u_t}^{j-\ell}V(z)$ is a sum of monomials whose coefficients involve powers of the coordinates of $u_t$. Since $\|u_t\|=1$, all these coordinates are bounded in absolute value by $1$. Therefore, $|\partial_{u_t}^{j-\ell}V(z)|$ admits a polynomial bound in $\|z\|$ that is uniform in $u_t$, and hence independent of $\sigma$. It follows that $|A_{j,t}(V)(z)|$ is bounded by a polynomial with positive coefficients evaluated at $\|z\|$. Moreover, as noted above, $I_\sigma^{-1}N_\sigma$ is independent of $\sigma$. Consequently, the constant $B_{l, A_{j,t}(V)}$ depends only on $l$ and $d$.
\end{proof}

\begin{proposition}\label{proposition:continuity_S_d}
Let $\pi$ be an element of $\mathrm{I}^m$. Let $(\varphi_n)_{n\in\N}$ be a sequence of elements of $\Cpol{0}$ such that 
\begin{align*}
\forall n \in \N, \; \forall x \in \R^d, \quad \left |[\varphi_n(x)]^\oplus \right | \leq C \left (1+\|I_\sigma^{-1} x \|^l \right ),
\end{align*}
for some $C>0$ and $l\in \N$. Suppose that the sequence $(\varphi_n)_{n\in\N}$ converges pointwise to a continuous function $\varphi$. Then the sequence $\left ( \S_\sigma^{\pi,m}(\varphi_n) \right )_{n\in\N}$ converges pointwise to $\S_\sigma^{\pi,m}(\varphi)$.
\end{proposition}
\begin{proof}
    The proof goes verbatim as the proof of Proposition \ref{proposition:continuity_S_d_1} dealing with the case $m=1$. First, we write the representation
    \begin{align}\label{aux:continuity-S-d-m}
        \S_\sigma^{\pi,m}(\varphi_n)(x) &= \frac{1}{m_\Gamma^{m}} \left ( [\varphi_n]^\oplus(x) + \int_0^{+\infty} K_m^{(m)}(t)\Pt_t^{\sigma}[\varphi_n]^\oplus(x) \, \dint t \right ) \\
     & \qquad + \sum_{j=0}^{m-1} \frac{(-1)^{m+j}}{m_\Gamma^{j+1}}\int_0^{+\infty} K_m^{(j)}(t) \mathcal{A}_{t}^{(m-1-j)}(\varphi_n)(x) \, \dint t \, . \nonumber
    \end{align}
    Secondly, we apply carefully Lebesgue's dominated convergence theorem to each integral appearing in \eqref{aux:continuity-S-d-m}.
\end{proof}

\noindent The next results follow as corollaries of the precedent Proposition \ref{proposition:prolongement_S_d_m}, combined with Lemma \ref{lemma:estimates_resolvent_polynomial_d}.

\begin{proposition}\label{proposition:continuity_Tmi_d}
Let $j$ be an element of $\mathrm{J}_m$ and let $\pi$ be an element of $\mathrm{I}^m$. Let $(\varphi_n)_{n\in\N}$ be a sequence of elements of $\Cpol{0}$ such that $\forall n \in \N, \;\rank(\varphi_n)\geq d_m+2m-1$ and 
\begin{align*}
\forall n \in \N, \; \forall x \in \R^d, \quad \left |[\varphi_n(x)]^\oplus \right | \leq C \left (1+\|I_\sigma^{-1} x \|^l \right ),
\end{align*}
for some $C>0$ and $l\in \N$. Suppose that the sequence $(\varphi_n)_{n\in\N}$ converges pointwise to a continuous function $\varphi$. Then the sequence $\left ( \boldsymbol{\Tj}^{\sigma,\,\pi}_{m,j}(\varphi_n) \right )_{n\in\N}$ converges pointwise to $\boldsymbol{\Tj}^{\sigma,\,\pi}_{m,j}(\varphi)$.
\end{proposition}

For the next proposition, we adapt and generalize the proof of Proposition \ref{proposition:continuity_Tjk_d}: no additional argument is needed.
\begin{proposition}\label{proposition:estimates_Tmi_d}
Let $j$ be an element of $\J_m$ and let $\pi$ be an element of $\mathrm{I}^m$. Let $\varphi$ be an element of $\Cpol{0}$ such that $\rank(\varphi) \geq d_m+2m-1$, and
\begin{align*}
\forall x \in \R^d, \quad  \left | [\varphi]^\oplus(x)\right | \leq C \left ( 1 +\| I_\sigma^{-1} x \|^l \right ),
\end{align*}
for some $C>0$, $l \in \N$. Then there exist $A_{l,m,j}^\pi>0$ such that
\begin{align}\label{ineq:estimate_Tmi_d}
\forall x \in \R^d, \quad \left | \boldsymbol{\Tj}^{\sigma,\,\pi}_{m,j}(\varphi)(x) \right | \leq CA_{l,m,j}^\pi \left ( 1 + \|I_\sigma^{-1} x \|^{l+\mu_m} \right )\, ,
\end{align}
where $\mu_m := 1 + \sum_{k=1}^m \mu_{k,\pi_k}$ and $\mu_{k,\pi_k} := \sum_{\ell = \nu_{\pi,k} }^{d_k} \ell - \nu_{\pi,k}+1$.
\end{proposition}

Propositions \ref{proposition:continuity_Tmi_d} and \ref{proposition:estimates_Tmi_d}  allow to extend \eqref{eq:IPP_final_sum_HO} to continuous functions.

\begin{proposition}\label{proposition:IPP_final_extended_HO}
    Let $\varphi$ be an element of $\Cpol{0}$ such that $\rank(\varphi) \geq d_m+2m-1$. Then one has
\begin{align}\label{eq:expansion_HO}
   \E[\varphi(F)(\Gamma(F,-\L^{-1}F)-1)]   =  \sum_{j \in \J_{m}} \sum_{\pi \in \mathrm{I}^m} \, \E\left [\boldsymbol{\Tj}_{m,j}^{\sigma, \, \pi} (\varphi) (F_1,\ldots, F_d)W_{j}^{\pi} \right ].
\end{align}
\end{proposition}

The proof of Proposition \ref{proposition:IPP_final_extended_HO} follows exactly the same arguments as the proof of Proposition \ref{proposition:IPP_final_extended} and is hence omitted.

\subsection{Estimation of the remainders}

In this section, we provide a bound for each term
\begin{align*}
    \E\left [\boldsymbol{\Tj}_{m,j}^{\sigma, \, \pi} (\varphi) (F_1,\ldots, F_d)W_{j}^{\pi} \right ]
\end{align*}
appearing in the expansion \eqref{eq:expansion_HO}. We use the same strategy and arguments used in Section \ref{sec.bounding}. Let us start with the following lemma. 

\begin{lemma}\label{lemma:bound.W}
Let $\pi$ be an element of $\mathrm{I}^m$ and let $j$ be an element of $\J_m$. Then there exist a constant $C_{d,m}>0$ such that
\begin{align*}
\E \left [ \left ( W_j^{\pi} \right )^2 \right ] \leq C_{d,m} \left ( \prod_{k=1}^m\, \Var\Gamma(F_{p_k},F_{q_k}) \right )  \Var\Gamma(F,\L^{-1} F) \, ,
\end{align*}
where $\pi = ((p_1,q_1),\ldots, (p_m,q_m))$.
\end{lemma}

\begin{proof}
    Since we have
    \begin{align*}
W_{j}^{\pi} := \overline{\Gamma}_{p_{m},q_{m}} \L^{j_m} \left [ \cdots \left [ \overline{\Gamma}_{p_{1},q_{1}}  \L^{j_{1}} \left [ \Gamma(F,-\L^{-1} F)-1 \right ] \right ] \cdots  \right ],
\end{align*}
using the Cauchy-Schwarz inequality, we get
\begin{align}\label{ineq:aux_Lbound}
    \E \left [ \left ( W_j^{\pi} \right )^2 \right ] \leq \E\left [\overline{\Gamma}_{p_{m},q_{m}}^4 \right ] ^{\frac{1}{2}}\E\left [ \left ( \L^{j_m} \left [ \cdots \left [ \overline{\Gamma}_{p_{1},q_{1}}  \L^{j_{1}} \left [ \Gamma(F,-\L^{-1} F)-1 \right ] \right ] \cdots  \right ] \right )^4 \right ]^{\frac{1}{2}}\, .
\end{align}
Then, using the hypercontractivity property \eqref{eq:hypercontractivity}, one has
\begin{align*}
    \E\left [\overline{\Gamma}_{p_{m},q_{m}}^4 \right ]^{\frac{1}{2}} \leq 3^{p_m+q_m-1} \E\left [\overline{\Gamma}_{p_{m},q_{m}}^2 \right ] \, ,
\end{align*}
and 
\begin{eqnarray*}
    && \E\left [ \left ( \L^{j_m} \left [ \cdots \left [ \overline{\Gamma}_{p_{1},q_{1}}  \L^{j_{1}} \left [ \Gamma(F,-\L^{-1} F)-1 \right ] \right ] \cdots  \right ] \right )^4 \right ]^{\frac{1}{2}} \\
    && \leq 3^{d_m} \,  \E\left [ \left ( \L^{j_m} \left [ \cdots \left [ \overline{\Gamma}_{p_{1},q_{1}}  \L^{j_{1}} \left [ \Gamma(F,-\L^{-1} F)-1 \right ] \right ] \cdots  \right ] \right )^2 \right ] .
\end{eqnarray*}
Now, since 
\begin{align*}
    \overline{\Gamma}_{p_{m-1},q_{m-1}}\L^{j_{m-1}} \left [ \cdots \left [ \overline{\Gamma}_{p_{1},q_{1}}  \L^{j_{1}} \left [ \Gamma(F,-\L^{-1} F)-1 \right ] \right ] \cdots  \right ] \in \WW_{\leq d_{m}} \, ,
\end{align*}
applying Lemma \ref{bounds-L} in \eqref{ineq:aux_Lbound} yields
\begin{align*}
    \E \left [ \left ( W_j^{\pi} \right )^2 \right ] \leq C_{d,m} \, \E\left [\overline{\Gamma}_{p_{m},q_{m}}^2 \right ] \E \left [ \left ( \overline{\Gamma}_{p_{m-1},q_{m-1}}\L^{j_{m-1}} \left [ \cdots \left [ \overline{\Gamma}_{p_{1},q_{1}}  \L^{j_{1}} \left [ \Gamma(F,-\L^{-1} F)-1 \right ] \right ] \cdots  \right ]\right )^2\right ] \, ,
\end{align*}
where $C_{d,m} =  3^{d_m+p_m+q_m-1}d_m^{2j_m}$. Iterating this process, we get the wanted assertion.
\end{proof}
Thus, using Lemma \ref{lemma:bound.W} and with the same arguments used for Lemma \ref{lemma:pre-boundation} and Proposition \ref{proposition:EE_stein_discrepancy_d_1}, we derive the following final estimation.

\begin{proposition}\label{proposition:EE_stein_discrepancy_d_HO}
There exists $C_{d,m} >0$ such that, for all $F \in \WW_{\leq d}$ with $\E[F]=0$ and $\E[F^2]=1$, and all $h\in \mathscr{C}^0_b(\R)$ with $\|h\|_\infty \le 1$,
\begin{align*}
    \bigl|\E \bigl[ \operatorname{Proj}_{\geq 2dm-1}(\phi'_h)(F)\,(\Ga(F,-\L^{-1} F)-1) \bigr] \bigr|
    \le C_{d,m}\, \left ( \sum_{1\leq p\leq q \leq d} \Var\Gamma(F_p,F_q) \right )^{\frac{m}{2}} \bigl(\Var \Ga(F,\L^{-1} F)\bigr)^{\frac{1}{2}}.
\end{align*}
\end{proposition}
\noindent This concludes the proof of Theorem \ref{thm:EE_Wiener_sums}.

\section{Refinements and limitations of the remainder bounds}\label{sec.shenanigans}

In this section, we explore the optimality of the added factor $\sum_{p\leq q} \Var\Gamma(F_p,F_q)$. As mentioned in the introduction, existing results on Edgeworth expansions for smooth functionals of Gaussian fields (see \cite{NP2009a}) indicate that the convergence rate is governed by
\begin{align*}
    \Var\Gamma(F,\L^{-1} F) \, .
\end{align*}
It is therefore natural to expect the corresponding Edgeworth expansion to be expressed in powers of this quantity. This is precisely what is established in \cite{MPS2025.EE} when $F$ belongs to a single Wiener chaos. Another non-asymptotic result is obtained in \cite{Celli2026}, where the total variation distance between the law of a neural network and its Edgeworth approximation is bounded by a power of $\tfrac{1}{n}$, consistently with the convergence rate identified in \cite{FHMNP2025}.

With that in mind, in a first subsection, we control 
\begin{align}
    \Lambda(F) =  \sum_{1\leq p\leq q \leq d} \Var\Gamma(F_p,F_q)  \, , 
\end{align}
by a power of 
\begin{align}
    \Delta(F) := \bigl ( \Var\Gamma(F,\L^{-1} F) \bigr )^{\frac{1}{2}}\, ,
\end{align}
hence yielding a genuine analogue of \cite[Theorem 1.2]{MPS2025.EE} in the case of finite sums of chaoses. Then, in the second subsection, we discuss this bound, as well as the lack of equivalence between the two quantities. We also provide alternate non-optimal ways to bound $\Lambda(F)$.

\subsection{A general comparison}
In this section, we prove the following control of $\Lambda(F)$ by a power of $\Delta(F)$, for elements $F$ lying in a finite sum of chaoses.
\begin{proposition}\label{prop:bound-by-var}
    Let $d\geq 2$ be an integer. There exists a constant $C_d>0$ such that, for every centered $F\in \WW_{\leq d}$ with $\E[F^2]=1$,
    \begin{align}\label{ineq:bound-by-var}
    \Lambda(F) \leq C_d \,\Delta(F)^{\frac{1}{d-1}}\,.
\end{align}
\end{proposition}
Before proving Proposition \ref{prop:bound-by-var}, we establish two preliminary lemmas. The first sharpens the usual monotonicity inequality for $L^p$-norms.

\begin{lemma}\label{lemma:sub-bound-var}
    Let $n$ be a positive integer. Then there exists a constant $C_n>0$ such that 
    \begin{align*}
        \forall Q \in \WW_{\leq n}, \; \forall r \geq 2, \quad \|Q\|_r \geq c_d r^{\frac{1}{2}} \bigl ( \Var(Q) \bigr )^{\frac{1}{2}} \, .
    \end{align*}
\end{lemma}
\begin{proof}
    Suppose that $Q \in \WW_{\leq n}$ is a polynomial with a finite number of variables, evaluated in a sequence of iid standard Gaussian random variables. An application of \cite[Theorem 3.3]{ODZ2016} (see also \cite[Theorem 1.3]{AW2015}) to $f = Q/\bigl(\Var(Q)\bigr)^{1/2}$, with $t=\sqrt{r}$, yields an absolute constant $C>0$ such that
    \begin{align*}
        \prob \left ( |f|> r^{\frac{1}{2}}\right ) \geq e^{-Crn^2} \, .
    \end{align*}
    Therefore we have 
    \begin{align*}
        \|f\|_r  = \E\left [ |f|^r \right ]^{\frac{1}{r}} \geq \E\left [ |f|^r\I_{\{|f| \geq r^{1/2} \} } \right ]^{\frac{1}{r}} \geq  r^{\frac{1}{2}} \prob \left ( |f|> r^{\frac{1}{2}}\right )^{\frac{1}{r}} \geq e^{-Cn^2} r^{\frac{1}{2}}\, ,
    \end{align*}
    that is 
    \begin{align}\label{aux:for-polynomial}
        \|Q\|_r \geq C_n r^{\frac{1}{2}}\bigl ( \Var(Q) \bigr )^{\frac{1}{2}}
    \end{align}
    where $C_n := e^{-Cn^2}$. 
    
    In the general case, we approximate $Q \in \WW_{\leq n}$ with polynomials of degree at most $n$ and finitely many variables evaluated in a sequence of iid Gaussian random variables. Then taking the limit in \eqref{aux:for-polynomial} leads to the result. For more details on the approximating argument, we refer the reader to \cite[Sections 3.1.3--3.1.4]{HMP2024}.
\end{proof}

The next lemma provides an estimate for the resolvent operator of $\L$ evaluated at $\Gamma(F,-\L^{-1}F)-1$.
\begin{lemma}\label{lemma:aux}
    Let $F$ be an centered element of $\WW_{\leq d}$ such that $\E[F^2] = 1$. Consider the resolvent operator
    \begin{align*}
        R := (\Id-\tfrac{1}{2}\L)^{-1} = \int_0^{+\infty} e^{-u} P_{\frac{u}{2}} \, \dint u\,,
    \end{align*}
    and set $ Z_F := R\bigl ( \Gamma(F,-\L^{-1}F)-1\bigr )$.
    Then one has
    \begin{align}\label{ineq:borne-ZF}
       \forall r \geq 2, \quad \|Z_F\|_r \leq (r-1)^{d-1} \Delta(F)\,.
    \end{align}
    \end{lemma}
    \begin{proof}
        One has 
        \begin{align*}
            \frac{1}{2}\L \Gamma(\L^{-1} F) = -\Gamma(F,-\L^{-1}F)+\Gamma_2(\L^{-1} F) \, ,
        \end{align*}
        where $\Gamma_2$ is the iterated gradient, defined by
        \begin{align*}
            \Gamma_2(F,G) := \frac{1}{2} \bigl ( \L \Gamma(F,G) -\Gamma(F,\L G) - \Gamma(\L F, G) \bigr )\,.
        \end{align*}
        Furthermore, one has
        \begin{align*}
            \Gamma_2(G) = \Gamma(G)+\|D^2G\|_{ HS}^2 \, ,
        \end{align*}
        for every $G \in \WW_{\leq d}$ (for more details, see e.g. \cite[Exercise 2.11.3]{NP2012}). Hence
        \begin{align*}
            \frac{1}{2}\L \Gamma(\L^{-1} F) = -\Gamma(F,-\L^{-1}F)+\Gamma (\L^{-1} F) + \|D^2(\L^{-1}F)\|_{ HS}^2\,,
        \end{align*}
        that is 
        \begin{align}\label{eq:ready-for-R}
           1+ \Gamma(F,-\L^{-1}F)-1 = \left ( \Id -\frac{1}{2} \L \right )\Gamma(\L^{-1}F) + \|D^2(\L^{-1}F) \|_{HS}^2 \, .
        \end{align}
        Since $R$ is a nonnegative operator and preserves constants, applying $R$ to \eqref{eq:ready-for-R} we reap
    \begin{align*}
        1+Z_F =\Gamma(\L^{-1}F)+R\,\|D^2 \L^{-1}F\|_{HS}^2 \, ,
    \end{align*}
    and $0 \leq \Gamma(\L^{-1}F) \leq 1+Z_F$. Moreover, one has
    \begin{align*}
        \E[Z_F]=0\, , \quad Z_F \in \WW_{\leq 2d-2}\,, \quad \|Z_F\|_2 \leq \Delta(F) \, , 
    \end{align*}
    by property of contraction of $(P_t)_{t\geq 0}$.
    Therefore, by the hypercontractivity property \eqref{eq:hypercontractivity}, we deduce the result.
    \end{proof}

We are now ready to prove Proposition \ref{prop:bound-by-var}.

\begin{proof}[Proof of Proposition \ref{prop:bound-by-var}]
    Let $F$ be an centered element of $\WW_{\leq d}$. Suppose that $\E[F^2] =1$. We divide the proof in several steps. 
    
\noindent \underline{Step 1}: For $s,t \geq 0$, define 
\begin{align*}
    U_{s,t} := \Gamma( P_s\L^{-1}F,P_t \L^{-1} F ) = \sum_{p,q=1}^d \frac{e^{-ps}e^{-qt}}{pq} \Gamma(F_p,F_q) \, .
\end{align*}
Denote by $M$ the Vandermonde matrix $(e^{-pj})_{1\leq j,p\leq d}$, so that we have
\begin{align*}
    \mathcal{U} = M \mathcal{B} M^{\top} \, ,
\end{align*}
where $\mathcal{U}_{i,j} = U_{i,j}$ and $\mathcal{B}_{p,q} = \tfrac{\Gamma(F_p,F_q)}{pq}$. Since $M$ is invertible, we deduce that each $\Gamma(F_p,F_q)$ is a linear combination of coordinates of $\mathcal{U}$, with coefficients only depending on $d$. Hence, there exists $C_d >0$ such that 
\begin{align}\label{aux:intermediary-bound}
    \Lambda(F) \leq C_d \underset{1\leq i,j \leq d}{\max} \Var( U_{i,j}) \, .
\end{align}
\noindent \underline{Step 2}: Fix $s,t >0$. Using Lemma \ref{lemma:sub-bound-var}, we get
\begin{align}\label{aux:borne-var-U}
   \forall r \geq 2, \quad \bigl ( \Var(U_{s,t}) \bigr )^{\frac{1}{2}} \leq C_d r^{-\frac{1}{2}} \|U_{s,t}\|_r \,.
\end{align}
Now, using the Cauchy--Schwarz inequality for $\Gamma$ (see, e.g., \cite{BGL2013}), we obtain \begin{align*} 
|U_{s,t}| & = | \Gamma( P_s \L^{-1} F, P_t \L^{-1} F )| \leq \bigl ( \Gamma(P_s \L^{-1} F) \bigr )^{\frac{1}{2}} \bigl ( \Gamma(P_t \L^{-1} F) \bigr )^{\frac{1}{2}} \\
& \leq \frac{1}{2} \left ( \Gamma(P_s \L^{-1} F)+\Gamma (P_t \L^{-1} F) \right ) \leq \frac{1}{2} \left ( P_s \Gamma(\L^{-1} F)+P_t\Gamma( \L^{-1} F) \right )\, , 
\end{align*} 
where the last inequality follows from \begin{align*} 
\Gamma(P_s \L^{-1} F) \leq e^{-2s} P_s \Gamma(\L^{-1} F) \leq P_s \Gamma(\L^{-1} F) \,. \end{align*} 
Then, using the positivity of $P_s$ together with Lemma \ref{lemma:aux}, we obtain 
\begin{align*} 
P_s \Gamma(\L^{-1} F)+P_t\Gamma( \L^{-1} F) \leq 2 + P_s Z_F + P_t Z_F \, , 
\end{align*} 
where $Z_F$ is defined in Lemma \ref{lemma:aux}. Finally, by positivity of the integrands, the contraction property of $P_s$, and the estimate \eqref{ineq:borne-ZF}, we obtain
\begin{align}\label{aux:borne-norm-U}
     \|U_{s,t}\|_r \leq \left \|1 +\frac{1}{2} \left (  P_s Z_F + P_t Z_F \right ) \right \|_r \leq 1 +r^{d-1} \Delta(F) \, .
\end{align}
\noindent \underline{Step 3}: Combining \eqref{aux:borne-var-U} and \eqref{aux:borne-norm-U}, we get 
\begin{align*}
   \forall r \geq 2,\quad  \bigl ( \Var(U_{s,t}) \bigr )^{\frac{1}{2}} \leq C_d r^{-\frac{1}{2}}\left (1+r^{d-1} \Delta(F) \right ) \, .
\end{align*}
Then, optimizing in $r$, we get $r = \Delta(F)^{-\frac{1}{d-1}}$. Now, if $\Delta(F) \leq 1$, then set $r = 2\Delta(F)^{-\frac{1}{d-1}} \geq 2$ and we get
\begin{align}\label{aux:intermediary-bound-bis}
    \bigl ( \Var(U_{s,t}) \bigr )^{\frac{1}{2}} \leq C_d\, \Delta(F)^{\frac{1}{2d-2}} \, ,
\end{align}
where $C_d$ have changed. If $\Delta(F) >1$ then since $\E[F^2]=1$, we already have the conclusion
\begin{align*}
    \Lambda(F)\leq A_d \leq A_d \,\Delta(F)^{\frac{1}{2d-2}} \, ,
\end{align*}
for some constant $A_d >0$. Combining \eqref{aux:intermediary-bound} and \eqref{aux:intermediary-bound-bis} concludes the proof.
\end{proof}

\subsection{Two chaoses with orders of opposite parity}
Consider the following model: let $p,q$ be positive integers such that $p\leq q$ and $p+q$ is odd. Then set 
\begin{align}\label{eq:sum-pq}
    F:=F_p+F_q\,,
\end{align}
with $F_p\in\WW_p$ and $F_q\in\WW_q$. Suppose furthermore that $\E[F^2]=1$. This setting was studied in \cite{BCKBS2025,DDMN2026}, where the authors bound the total variation distance between the law of $F$ and the standard normal distribution by a constant times the square root of the fourth cumulant of $F$. This establishes a fourth moment theorem for this class of random variables, although such a theorem does not hold for arbitrary finite sums of Wiener chaoses.

\begin{proposition}\label{prop:comparison-pq}
Let $1\le p<q$ be integers such that $p+q$ is odd. Let $F_p\in\WW_p$ and $F_q\in\WW_q$, and assume that $F=F_p+F_q$ satisfies $\E[F^2]=1$. Then
\begin{equation}\label{eq:comparison-pq}
\Lambda(F)\le q^2\Delta(F)^2 \, .
\end{equation}
\end{proposition}

To prove Proposition \ref{prop:comparison-pq}, we need two lemmas. The first one below is a direct consequence of \cite[Theorem~1.6]{MNPP2016} in the case $d=2$.

\begin{lemma}\label{lemma:cov-squares}
Let $p,q$ be nonnegative integers, $U\in\WW_p$ and $V\in\WW_q$.
Then
\begin{equation}\label{eq:cov-square-positivity}
 \operatorname{Cov}(U^2,V^2) = \E[U^2V^2]-\E[U^2]\E[V^2]\geq 0 \, .
\end{equation}
\end{lemma}

The second one is the same result but a version for the carré-du-champ operator $\Gamma$, consequence of Lemma \ref{lemma:cov-squares}.

\begin{lemma}\label{lemma:cov-gamma}
Let $p,q$ be nonnegative integers, $U\in\WW_p$ and $V\in\WW_q$.
Then
\begin{equation}\label{eq:cov-gamma-positivity}
 \operatorname{Cov}(\Gamma(U),\Gamma(V)) \geq 0 \, .
\end{equation}
\end{lemma}
\begin{proof}
If $p=1$ or $q=1$, one of $\Gamma(U)$ and $\Gamma(V)$ is deterministic, and the conclusion is immediate. Hence we assume that $p,q\geq2$. 
Let $(e_i)_{i\geq1}$ be an orthonormal basis of $\mathfrak{H}$ (defined in Section \ref{sec.contraction}), and set
\begin{align*}
    D_iU:=\langle DU,e_i\rangle_{\mathfrak H}\, ,\qquad D_iV:=\langle DV,e_i\rangle_{\mathfrak H} \, .
\end{align*}
Since $D_iU\in\WW_{p-1}$ and $D_jV\in\WW_{q-1}$, Lemma \ref{lemma:cov-squares} yields
\begin{align*}
    \operatorname{Cov}\left ((D_iU)^2,(D_jV)^2\right) \geq 0 \, ,
\end{align*}
for every $i,j$ positive integers. For $n\geq 1$, let us define the partial sums
\begin{align*}
    A_n:=\sum_{i=1}^n (D_iU)^2 \, , \qquad B_n:=\sum_{j=1}^n (D_jV)^2 \, .
\end{align*}
By bilinearity of the covariance, one has
\begin{align*}
    \operatorname{Cov}(A_n,B_n) = \sum_{i,j=1}^n \operatorname{Cov}\left ((D_iU)^2,(D_jV)^2\right) \geq 0 \, .
\end{align*}
By Parseval's identity, we have
\begin{align*}
    A_n \xrightarrow[n\to +\infty]{\text{a.s.}} \Gamma(U) \, , \qquad  B_n \xrightarrow[n\to +\infty]{\text{a.s.}} \Gamma(V) \, .
\end{align*}
Moreover, both $\Gamma(U)$ and $\Gamma(V)$ are $L^2$, since they belong to finite sums of Wiener chaoses. Thus, using the fact that
\begin{align*}
    \forall n \in \N^*, \quad 0\leq A_n \leq \Gamma(U)\, , \quad  0\leq B_n \leq \Gamma(V) \, ,
\end{align*}
the dominated convergence theorem gives
\begin{align*}
    A_n \xrightarrow[n\to +\infty]{L^2} \Gamma(U) \, , \qquad  B_n \xrightarrow[n\to +\infty]{L^2} \Gamma(V) \, .
\end{align*}
By continuity of the covariance under $L^2$-convergence, we conclude that
\begin{align*}
    \operatorname{Cov}(\Gamma(U),\Gamma(V)) = \lim_{n\to\infty}\operatorname{Cov}(A_n,B_n) \geq 0 \, .
\end{align*}
\end{proof}
Now, we have everything to prove Proposition \ref{prop:comparison-pq}.
\begin{proof}[Proof of Proposition \ref{prop:comparison-pq}]
By symmetry and bilinearity of $\Gamma$, one has
\begin{equation*}
 \Gamma(F,-L^{-1}F) = \frac{1}{p}\Gamma(F_p) + \frac{1}{q}\Gamma(F_q)
   + \left ( \frac{1}{p} + \frac{1}{q} \right ) \Gamma(F_p,F_q) \, ,
\end{equation*}
hence, since $\E[F^2]=1$, we get
\begin{align}\label{eq:Gamma-L}
    \Gamma(F,-L^{-1}F)-1 = \frac{1}{p}A + \frac{1}{q}B + \left ( \frac{1}{p} + \frac{1}{q} \right ) \Gamma(F_p,F_q) \, ,
\end{align}
where $A := \Gamma(F_p)-p\sigma_p^2$ and $B:= \Gamma(F_q)-q\sigma_q^2$.
Let us write $F_a = I_a(f_a)$, with $f_a \in \mathfrak{H}^{\odot a}$, for every $a \in \{p,q\}$. Using the product formula \eqref{eq:product-formula}, and the formula
\begin{align*}
    \Gamma(F_a,F_b) = \frac{1}{2}(\L+(a+b)\Id)(F_aF_b) \, ,
\end{align*}
one gets 
\begin{align*}
    \Gamma(F_a,F_b) = ab\sum_{r=1}^{a\wedge b} (r-1)! \binom{a-1}{r-1} \binom{b-1}{r-1} I_{a+b-2r} \left( f_a\widetilde{\otimes}_r f_b \right) \, ,
\end{align*}
for every $a,b$ in $\{p,q\}$.
Since the orders of the chaoses in the decomposition of $\Gamma(F_p)$ (resp. $\Gamma(F_q)$) are even, we have, after re-centering,
\begin{align*}
    A \in \underset{k=1}{\overset{p-1}{\bigoplus}} \, \WW_{2k} \, , \qquad  B \in \underset{k=1}{\overset{q-1}{\bigoplus}} \, \WW_{2k} \, .
\end{align*}
On the other hand, since $p+q$ is odd, every order $p+q-2r$ in the expansion of $\Gamma(F_p,F_q)$ is odd, that is
\begin{align*}
    \Gamma(F_p,F_q) \in \underset{\substack{1\leq k \leq p+q-2 \\ k \text{ odd}}}{\bigoplus} \, \WW_k \, .
\end{align*}
In particular, by orthogonality, we deduce that
\begin{align*}
    \E[A \Gamma(F_p,F_q)] = \E[B \Gamma(F_p,F_q)] =0\,.
\end{align*}
Therefore, squaring \eqref{eq:Gamma-L} and taking the expectation, we get
\begin{align*}
    \Delta(F)^2 = \frac{1}{p^2}\Var\Gamma(F_p) + \frac{1}{q^2}\Var\Gamma(F_q) + \frac{2}{pq} \operatorname{Cov}(\Gamma(F_p),\Gamma(F_q)) + \left ( \frac{1}{p} +\frac{1}{q} \right )^2 \Var\Gamma(F_p,F_q) \, .
\end{align*}
By Lemma \ref{lemma:cov-gamma}, one has $\operatorname{Cov}(\Gamma(F_p),\Gamma(F_q)) \geq 0$. Consequently, since $p<q$, we get
\begin{align*}
    \Delta(F)^2 \geq \frac{1}{q^2} \left ( \Var\Gamma(F_p) + \Var\Gamma(F_q) + + \Var\Gamma(F_p,F_q) \right ) = \frac{1}{q^2}\Lambda(F) \, .
\end{align*}
\end{proof}

Combining Proposition \ref{prop:comparison-pq} with Theorem \ref{thm:EE_Wiener_sums}, we get the following corollary.

\begin{corollary}
Let $m$ be a positive integer. Then there exists $C_{d,m}>0$ such that for all $F$ of the form \eqref{eq:sum-pq} with $\E[F]=0$ and $\E[F^2]=1$,
\begin{align*}
\dTV(\prob_F,\boldsymbol{\gamma}_{F,M}) & \leq C_{q,m} \, \bigl ( \Var\Gamma(F,\L^{-1} F) \bigr )^{\frac{m+1}{2}}\, ,
\end{align*}
where $M = 2qm$.
\end{corollary}

\subsection{Examples and counterexamples}

Here, we provide below an example of non-optimality of $\Lambda(F)$ and compare it to $\Delta(F)^2$.
\begin{example}[Cancellations between chaotic components]
The bound \eqref{ineq:borne.gamma.moins1} is, in general, sub-optimal. Consider the setting of Example \ref{example:un} and define
\begin{align*}
    F_n := F_{n,2} + F_{n,4} \, ,
\end{align*}
where 
\begin{align*}
    F_{n,2} := \sqrt{1-\frac{4}{9n}} \widetilde{F}_{n,2} \, , \quad F_{n,4} := -\frac{2}{3\sqrt n} \widetilde{F}_{n,4} \, .
\end{align*}
Then $\E[F_n^2] = 1$, for all $n$. Using the upcoming computations in Example \ref{example:un}, we compute
\begin{align*}
    \Lambda(F_n) \sim \frac{104}{9n} \, , \qquad \Delta(F_n)^2 \sim \frac{392}{81n^2} \, , 
\end{align*}
thus showing that \eqref{ineq:borne.gamma.moins1} is not optimal, as well as \eqref{ineq:bound-by-var}. Furthermore, it shows that the result of Proposition \ref{eq:comparison-pq} do not hold on general for an arbitrary sum of two chaoses.
\end{example}

Another way of bounding the error in the Edgeworth approximation is to obtain a bound involving only the variances of $\Gamma(F_p)$ and $\Gamma(F_q)$. For that, we apply Lemma \ref{lemma:estimates_variance_pq}, which yields
\begin{align}\label{ineq:bound1}
\Var\Gamma(F_p,F_q) \; \lesssim_{p,q}  \;\Var\bigl(\Gamma(F_p)\bigr) + \Var\bigl(\Gamma(F_q)\bigr) + \sigma_p^2 \Var\bigl(\Gamma(F_q)\bigr)^{1/2}\,,
\end{align}
for $p<q$. We then use the more convenient bound
\begin{align}\label{ineq:bound2}
\Var\Gamma(F_p,F_q) \; \lesssim_{p,q}  \; \sqrt{\Var\Gamma(F_p) + \Var\Gamma(F_q)} \, .
\end{align}
Neither bound is asymptotically sharp in general, as the following examples illustrate.

\begin{example}[Loss through contraction estimates $1$]\label{example:un}
Consider the sequence of random variables $(F_n)$ defined as follow. Let $X_1, \ldots, X_n$ be iid standard Gaussian random variables, and denoted by $R_n$ the sum $\sum_{i=1}^n X_i^2$. Then define
\begin{align*}
\widetilde{F}_{n,2} \coloneq \frac{R_n-n}{\sqrt{2n}} \in \WW_2\,, \quad
 \widetilde{F}_{n,4}  \coloneq \frac{R_n^2-2(n+2)R_n+n(n+2)}{\sqrt{8n(n+2)}} \in \WW_4\,,
\end{align*}
($\widetilde{F}_{n,4}$ is the normalised projection of $R_n^2$ onto $\WW_4$) and finally define
\begin{align*}
F_n & \coloneq  F_{n,2}+F_{n,4} \, ,
\end{align*}
where $F_{n,2} := \sqrt{1-\tfrac{1}{4n}} \widetilde{F}_{n,2}$ and $F_{n,4} := -\tfrac{1}{2\sqrt{n}} \widetilde{F}_{n,4}$. 
One has $\forall n \in \N^*, \; \E[F_n^2] = 1$. Furthermore, 
\begin{align*}
\Gamma(F_{n,4}) = \frac{2R_n}{4n^2(n+2)} \left ( R_n-n-2 \right )^2\,, \quad \Gamma( F_{n,2} ) = \left (1-\frac{1}{4n} \right )\frac{2}{n}R_n\,,
\end{align*}
and 
\begin{align*}
\Gamma( F_{n,2}, F_{n,4}) = -\frac{\left ( 1-\tfrac{1}{4n} \right )^{\frac{1}{2}}}{n\sqrt{n}\sqrt{n+2}}R_n(R_n-n-2)\,.
\end{align*}
Thus,
\begin{align*}
\Var\Gamma(F_{n,2}) = \left ( 1- \frac{1}{4n} \right )^{2}\frac{8}{n}\,, \quad \Var\Gamma(F_{n,4}) = \frac{2(n^2+31n+106)}{n^3(n+2)}\,,
\end{align*}
and
\begin{align*}
     \Var\Gamma(F_{n,2},F_{n,4}) = \frac{2(n+6)}{n^2}\left ( 1- \frac{1}{4n} \right ) \, .
\end{align*}
Consequently, the bound 
\begin{align*}
    \Var \Gamma(F_{n,2},F_{n,4}) \; \lesssim \; \sqrt{ \Var\Gamma(F_{n,2}) + \Var\Gamma(F_{n,4})} 
\end{align*}
has right-hand side of order $n^{-\frac{1}{2}}$ whereas its left-hand side has order $n^{-1}$.
\end{example}

\begin{example}[Loss through contraction estimates $2$]\label{example:deux}
In the same setting as in Example \ref{example:un}, define
\begin{align*}
F_{n,2} :=\frac{R_n-n}{\sqrt{2n}} \, , \qquad F_{n,3} :=\frac{1}{\sqrt{6n}}\sum_{i=1}^n(X_i^3-3X_i) \, .
\end{align*}
Then $F_{n,2}\in\WW_2$, $F_{n,3}\in \WW_3$, and
\begin{align*}
    \mathbb{E}[F_{n,2}^2]=\mathbb E[F_{n,3}^2]=1 \,.
\end{align*}
On one hand,
\begin{align*}
    \Gamma(F_{n,2},F_{n,3})
=\frac{\sqrt{3}}{n}\sum_{i=1}^n(X_i^3-X_i)\, .
\end{align*}
On the other hand,
\begin{align*}
    \Gamma(F_{n,2})=\frac{2}{n}\sum_{i=1}^nX_i^2,
\qquad
\Gamma(F_{n,3})
=\frac{3}{2n}\sum_{i=1}^n(X_i^2-1)^2 \, ,
\end{align*}
We find, after computations,
\begin{align*}
\operatorname{Var} \Gamma(F_{n,2})=\frac{8}{n} \, , \quad \operatorname{Var} \Gamma(F_{n,3}) =\frac{126}{n} \, , \quad \operatorname{Var}\Gamma(F_{n,2},F_{n,3})=\frac{30}{n}
\end{align*}
Consequently, the bound
\begin{align*}
\operatorname{Var}\Gamma(F_{n,2},F_{n,3})
\lesssim
\operatorname{Var}\Gamma(F_{n,2})
+\operatorname{Var}\Gamma(F_{n,3})
+\E[F_{n,2}^2]
\sqrt{\operatorname{Var}\Gamma(F_{n,3})} 
\end{align*}
has right-hand side of order $n^{-\frac{1}{2}}$ whereas its left-hand side has order $n^{-1}$. 
\end{example}
\begin{remark}
Note that the estimates obtained by using the contractions are not available in the broader setting of Markov chaoses, discussed in the introduction.
\end{remark}

\section*{Acknowledgments}
The author warmly thanks his PhD supervisors Guillaume Poly and Yvik Swan for their help and support in this project, as well as Ronan Herry for useful and stimulating discussions. The author is a FRIA grantee of the Fonds de la Recherche Scientifique -- FNRS. This work was supported by ANR LESSBIG.

\section*{Use of generative AI}

The author did not use generative artificial intelligence in the development of the mathematical arguments, proofs, or results presented in this article. AI-assisted tools were used only to support bibliographic searches and to improve the English exposition. All cited references were independently checked by the author, who takes full responsibility for the content of the manuscript.
\printbibliography

\end{document}